\documentstyle{amsppt}
\input amstex
\nologo\magnification=1200
\overfullrule=0pt
\hsize = 5.5 true in\hoffset=.5in
\vsize = 9 true in
\input BoxedEPS
\SetTexturesEPSFSpecial 
\loadeusm

\def\fn{{\frak n}}

\def\fnl{{\fn_{\scriptscriptstyle L}}}

\def\fnm{{\fn_{\scriptscriptstyle M}}}
\def\fnn{{\fn_{\scriptscriptstyle N}}}

\def\fnd#1{{\fn_{#1}}}
\def\fndm{{\fnd{\!\scriptscriptstyle M}\!}}

\def\fS*{{\frak S_*}}
\def\fs{{\frak s}}

\def\fsj{{\fs_j}}
\def\fsd#1{{\fs_{#1}}}
\def\fsdh{{\fsd{H}}}

\def\ds{{\dot s}}

\def\dstp{{\ds_{T'}}}

\def\sA{{\eusm A}}

\def\sE{{\eusm E}}

\def\sL{{\eusm L}}

\def\sM{{\eusm M}}

\def\sMdone{{\sM_1}}
\def\sMdl{{\sM_\la}}
\def\sMdz{{\sM_0}}
\def\sMdzone{{\sM_{0,1}}}
\def\sMdzoneuth{{\sM_{0,1}^\th}}
\def\sMdzoneu#1{{\sM_{0,1}^{#1}}}
\def\tM{{\widetilde {\sM}}}

\def\tMdld{{\widetilde {\sM}_{\la, d}}}
\def\tMd{{\tM_d}}
\def\sN{{\eusm N}}

\def\sU{{\eusm U}}

\def\sR{{\eusm R}}

\def\sRdone{{\sR_1}}

\def\sRdzone{{\sR_{0,1}}}
\def\sRdl{{\sR_\la}}

\def\C{{\Bbb C}}

\def\R{{\Bbb R}}
\def\T{{\Bbb T}}

\def\Z{{\Bbb Z}}
\def\Zp{{\Z_p}}

\def\Q{{\Bbb Q}}
\def\N{{\Bbb N}}

\def\bbjk{{b(j,k)}}

\def\bbij{{b(i,j)}}
\def\bbji{{b(j,i)}}

\def\bbonetwo{{b(1, 2)}}

\def\bbtwoone{{b(2,1)}}

\def\bbkj{{b(k,j)}}

\def\bbii{{b(i,i)}}
\def\a{{\alpha}}

\def\da{{\dot\a}}

\def\b{{\beta}}

\def\f{{\varphi}}
\def\p{{\psi}}
\def\z{{\zeta}}
\def\Ad{{\text{\rm Ad}}}

\def\edi{{e_i}}
\def\edonen{{e_{1,N}}}
\def\edtwon{{e_{2,N}}}
\def\edin{{e_{i,N}}}
\def\edjn{{e_{j,N}}}
\def\edkn{{e_{k,N}}}

\def\edim{{e_{i,M}}}

\def\edj{{e_j}}

\def\edk{{e_k}}

\def\edii{{e_{i,i}}}

\def\ed#1{{e_{#1}}}

\def\tedi{{\tilde e_i}}

\def\tedj{{\tilde e_j}}

\def\tedk{{\tilde e_k}}

\def\tilg{{\tilde g}}

\def\tilh{{\tilde h}}

\def\tp{{\tilde p}}

\def\tq{{\tilde  q}}

\def\tu{{\tilde u}}

\def\txti{{\text{\rm i}}}
\def\txk{{\text{\rm k}}}
\def\txm{{\text{\rm m}}}

\def\txmd#1{{\txm_{#1}}}
\def\txmdi{{\txmd{i}}}

\def\txmthth1{{\txm_{\thth1}}}
\def\txn{{\text{\rm n}}}

\def\txs{{\text{\rm s}}}

\def\Qm{{Q_\txm}}

\def\thds{{\th_s}}

\def\thdstp{{\th_{\dstp}}}

\def\thdtprime{{\th_{T'}}}

\def\thdtprimeun{{\th_{T'}^n}}
\def\thdm{{\th_\txm}}

\def\th{{\theta}}
\def\tht{{\theta_T}}
\def\tht'{{\th_{T'}}}
\def\la{{\lambda}}

\def\lada{{\la_a}}
\def\ladb{{\la_b}}

\def\ladauijk{{\la_a^{i,j,k}}}

\def\ladaz{{\lad{a, 0}}}

\def\ladab{{\la_{a, b}}}
\def\ladasa{{\la_{\scriptscriptstyle{\AS a}}}}

\def\lad#1{{\la_{#1}}}

\def\ladmu{{\lad{\mu}}}

\def\mum{{\mu_M}}

\def\muda{{\mu_a}\!}
\def\mudb{{\mu_b}\!}

\def\mudauijk{{\mu_a^{i,j,k}}\!}

\def\mudainv{{\mu_a^{-1}\!}}

\def\Aut{{\text{\rm Aut}}}
\def\Autf'{{\Aut_\f'}}
\def\Autp'{{\Aut_\p'}}
\def\Cntp'{{\Cnt_\p'}}
\def\Cnt{{\text{\rm Cnt}}}
\def\Ob{{\text{\rm Ob}}}
\def\Obm{{\Ob_{\text{\rm m}}}}
\def\cnt{{\Cnt}}
\def\cntr{{\Cnt_{\text{\rm r}}}}
\def\Hom{{\text{\rm Hom}}}

\def\Int{{\text{\rm Int}}}
\def\Intp'{{\Int_\p'}}

\def\Ker{{\text{\rm Ker}}}

\def\mod{{\text{\ \rm mod\ }}}
\def\modd{{\text{\!\rm\mod\!\!}}}

\def\out{{\text{\rm out}}}

\def\Autf{{\Aut_\f}}

\def\Autf'{{\Aut_\f'}}
\def\Proj{{\text{\rm Proj}}}

\def\Sp{{\text{\rm Sp}}}
\def\Tr{{\text{\rm Tr}}}
\def\two{{\rm I\!I}}
\def\twoone{{\rm I\!I$_1$}}
\def\threee{{\text{\rm I\!I\!I}}}
\def\threeone{{\rm I\!I\!I$_1$}}
\def\threel{{\rm I\!I\!I$_{\lambda}$}}
\def\three0{{\rm I\!I\!I$_0$}}
\def\threez{{\three0}}
\def\twoinf{{\rm I\!I$_{\infty}$}}

\def\tsU{{\widetilde {\sU}}}
\def\tsUdz{{\tsU_0}}

\def\ta{{\widetilde {\a}}}

\def\d{{\delta}}
\def\bd{{\Bar\d}}

\def\Ddone{{D_1}}

\def\Ddi{{D_i}}
\def\Ddij{{D_{i,j}}}

\def\Edij{{E_{i,j}}}

\def\tmu{{\tilde \mu}}
\def\tmuda{{\tilde\mu_a}}

\def\La{{\Lambda}}
\def\g{{\gamma}}

\def\ggdz{{g_0}}
\def\ggdone{{g_1}}
\def\dotggdone{{\dot g_1}}

\def\ggdtwo{{g_2}}

\def\dotggdtwo{{\dot g_2}}
\def\ggdthree{{g_3}}

\def\dotggdthree{{\dot g_3}}

\def\ggdonedotsthree{{\ggdone, \ggdtwo, \ggdthree}}

\def\hdone{{h_1}}
\def\hdtwo{{h_2}}

\def\hdoneinv{{h_1^{-1}}}

\def\ggdi{{g_i}}

\def\ggd#1{{g_{#1}}}

\def\ggdj{{g_j}}

\def\ggdjone{{g_{j+1}}}

\def\ggdkmone{{g_{k-1}}}
\def\ggdk{{g_k}}

\def\ggdkone{{g_{k+1}}}
\def\ggdktwo{{g_{k+2}}}

\def\ggdnmone{{g_{n-1}}}
\def\ggdn{{g_n}}
\def\ggdnone{{g_{n+1}}}

\def\gdotsgn{{\ggdone, \ggdtwo, \cdots, \ggdn}}
\def\gdotsgnmone{{\ggdone, \ggdtwo, \cdots, \ggdnmone}}
\def\gdotsgtwon{{\ggdtwo, \ggdthree, \cdots,
\ggdn}}

\def\gdotsgtwonmone{{\ggdtwo, \ggdthree, \cdots,
\ggdnmone}}
\def\ggdnmone{{g_{n-1}}}
\def\gdotsgthrn{{\ggdthree,  \cdots,
\ggdn}}

\def\gsigdotsn{{g_{\sig(1)}, g_{\sig(2)}, \cdots,
g_{\sig(n)}}}

\def\gdotsgnmone{{\ggdone, \ggdtwo, g_3, \cdots,
\ggdnmone}}

\def\ggdsigone{{g_{\sig(1)}}}

\def\ggdsigone{{g_{\sig(1)}}}
\def\ggdsigtwo{{g_{\sig(2)}}}

\def\ggdsign{{g_{\sig(n)}}}

\def\md#1{{m_{#1}}}

\def\mdonetwo{{\md{12}}}
\def\mdonethr{{\md{13}}}
\def\mdonefour{{\md{14}}}

\def\mdtwothr{{\md{23}}}
\def\mdtwofour{{\md{24}}}

\def\mdjk{{m_{jk}}}

\def\hdone{{h_1}}
\def\hdtwo{{h_2}}

\def\kdone{{k_1}}

\def\kdj{{k_j}}
\def\kdjone{{k_{j+1}}}
\def\elldone{{\ell_1}}

\def\elldi{{\ell_i}}

\def\elldk{{\ell_k}}

\def\tilg{{\tilde g}}

\def\wdz{{w_0}}
\def\dotwdz{{\dot w_0}}
\def\wdone{{w_1}}
\def\dotwdone{{\dot w_1}}
\def\zdz{{z_0}}

\def\tilh{{\tilde h}}

\def\part{{\partial}}

\def\partg{{\part_{\scriptscriptstyle{G}}}}
\def\tpartg{{\tilde\part_{\scriptscriptstyle{G}}}}

\def\partl{{\part_L}}
\def\partdone{{\part_1}}
\def\partdtw{{\part_2}}
\def\partm{{\part_{\scriptscriptstyle  M}}}
\def\partn{{\part_{\scriptscriptstyle  N}}}

\def\partgm{{\part_{G_\txm}}}
\def\parth{{\part_{\scriptscriptstyle  H}}}
\def\parthm{{\part_{\scriptscriptstyle  H_\txm}}}
\def\partq{{\part_{\scriptscriptstyle  Q}}}
\def\partqm{{\part_\Qm}}

\def\log{{\text{\rm log}}}

\def\sig{{\sigma}}

\def\sigp{{\sigma^{\psi}}}

\def\sigps{{\sigma_s^{\psi}}}
\def\sigpdT{{\sig_T^\p}}

\def\siguz{{\sig^0}}
\def\sigdz{{\sig_0}}

\def\sigdguz{{\sig_g^0}}
\def\sigulmu{{\sig^{\la, \mu}}}

\def\sigdkkone{{\sig_{k, k+1}}}

\def\fsfntr{{faithful semi-finite normal trace}}
\def\wt{{semi-finite normal weight}}
\def\fwt{{faithful \wt}}
\def\botimes{{\Bar \otimes}}
\def\id{{\text{\rm id}}}

\def\r0{{\sR_0}}

\def\r01{{\sR_{0,1}}}

\def\Map{{\text{\rm Map}}}

\def\botimes{{\overline \otimes}}
\def\wt{{\text\allowlinebreak{semi-finite normal weight}}}

\def\fwt{{faithful \wt}}

\def\vna{{\text\allowlinebreak{von Neumann algebra}}}
\def\vnas{{\text\allowlinebreak{von Neumann algebras}}}

\def\tB{{\text{\rm B}}}

\def\tC{{\text{\rm C}}}

\def\tCta#1{{\tC_\ta^{#1}}}

\def\tCtaun{{\tCta{n}}}
\def\tH{{\text{\rm H}}}

\def\tHau#1{{\tH_\a^{#1}}}
\def\tHtau#1{{\tH_\ta^{#1}}}
\def\tHaun{{\tHau{n}}}

\def\tHaunmone{{\tHau{n-1}}}
\def\tHtaunmone{{\tHtau{n-1}}}

\def\tZ{{\text{\rm Z}}}
\def\tZu#1{{\text{\rm Z}^{#1}}}
\def\tZunlr#1{{\tZu{n}\lr{#1}}}
\def\tZunmonelr#1{{\tZu{n-1}\lr{#1}}}

\def\tZa#1{{\tZ_\a^{#1}}}

\def\tZaun{{\tZa{n}}}
\def\tZaunmone{{\tZa{n-1}}}

\def\QED{{\hfill$\heartsuit$}}

\def\inv{{^{-1}}}

\document

\def\X{{\text{\rm Xext}}}

\def\res{{\text{\rm res}}}
\def\Res{{\text{\rm Res}}}

\def\inf{{\text{\rm inf}}}

\def\fsm{{\fs_\txm}}

\def\dfsm{{\dot \fs_\txm}}
\def\dfs{{\dot\fs}}

\def\sh{{\fs\!_{\scriptscriptstyle H}\!}}
\def\shs{{\fs\!_{\scriptscriptstyle H}^*\!}}

\def\sj{{\fs_j}}

\def\pig{{\pi\!_{\scriptscriptstyle G}\!}}

\def\piginv{{\pi\!_{\scriptscriptstyle G}^{-1}\!}}
\def\pigs{{\pi\!_{\scriptscriptstyle G}^*\!}}

\def\piq{{\pi_{\scriptscriptstyle Q}}}
\def\piqm{{\pi_{\scriptscriptstyle \Qdm}}}
\def\piqs{{\pi_{\scriptscriptstyle Q}^*}}

\def\pim{{\pi_{\txm}}}

\def\dpim{{\dot\pi_{\txm}}}
\def\pims{{\pi_{\txm}^*}}

\def\etapp{{\eta_p}}
\def\etappi{{\eta_\pdi}}
\def\etappj{{\eta_\pdj}}
\def\etappk{{\eta_\pdk}}

\def\etat{{{\eta_{T}}}}

\def\etat'{{\eta_{T'}}}

\def\etad#1{{\eta_{#1}}}

\def\etada{{\etad{a}}}

\def\_#1{{_{#1}}}
\def\^#1{{^{#1}}}

\def\(#1){{^{({#1})}}}

\def\scirc{{\lower-.3ex\hbox{{$\scriptscriptstyle\circ$}}}}

\def\bracettt#1{{\left\{#1\right\}_T}}
\def\bracett'#1{{\left\{\!#1\right\}_{T'}}}
\def\bracet'm#1{{\bracett'{\txm\left({#1}\right)}}}
\def\bracep#1{{\left\{#1\right\}_p}}

\def\bracedpi#1{{\left\{{#1}\right\}_\pdi}}
\def\bracedpj#1{{\left\{{#1}\right\}_\pdj}}
\def\bracedpk#1{{\left\{#1\right\}_\pdk}}
\def\bracketd#1#2{{\left[#1\right]_{#2}}}

\def\bracketp#1{{\left[#1\right]_{p}}}

\def\bracketpdi#1{{\left[#1\right]_\pdi}}
\def\bracketpdj#1{{\left[#1\right]_\pdj}}
\def\bracketpdk#1{{\left[#1\right]_\pdk}}

\def\brackett#1{{\left[{#1}\right]_{T}}}
\def\brackett'#1{{\left[{#1}\right]_{T'}}}
\def\rtz{{\R/T\Z}}
\def\rt'z{{\R/T'\Z}}

\def\cdag{{\text\allowlinebreak{countable discrete amenable group}}}
\def\cdabg{{\text\allowlinebreak{countable discrete abelian group}}}
\def\cdfabg{{\text\allowlinebreak{countable discrete free abelian group}}}

\def\hjr{{{\text{\rm HJR}-exact sequence}}}
\def\mhjr{{\text modified \hjr}}

\def\tcone{{\tC^1}}

\def\tctw{{\tC^2}}
\def\tcthr{{\tC^3}}

\def\tzda{{\tZ_a}}
\def\tzdb{{\tZ_b}}

\def\tzdwha{{\tZ_\whaa}}

\def\tcn{{\tC^n}}

\def\tcnone{{\tC^{n+1}}}

\def\tbmsout{{\tB_{\txm, \fs}^\out}}

\def\tbth1{{\tB_\th^1}}

\def\tbmsout{{\tB_{\txm, \fs}^\out}}
\def\tcone{{\tC^1}}
\def\tctw{{\tC^2}}
\def\tzth1{{\tZ_\th^1}}
\def\thth1{{\tH_\th^1}}
\def\ththr{{\tH^3}}

\def\thasout{{\tH_{\a, \fs}^\out}}

\def\thmsout{{\tH_{\txm, \fs}^\out}}

\def\tzthr{{\tZ^3}}

\def\tbthr{{\tB^3}}

\def\tcthr{{\tC^3}}
\def\thtw{{\tH^2}}

\def\tztw{{\tZ^2}}

\def\tzout{{\tZ^\out}}

\def\tzmsout{{\tZ_{\txm, \fs}^\out}}

\def\tbtw{{\tB^2}}

\def\pr{{\text{\rm pr}}}
\def\prdone{{\pr_1}}
\def\prdtwo{{\pr_2}}

\def\pijj'{{\pi_{J, J'}}}
\def\pijj's{{\pi_{J, J'}^*}}

\def\tmu{{\tilde \mu}}

\def\mudtxa{{\mu_{\text {\rm a}}}}

\def\linfrt'z{{L^\infty(\R/T'\Z)}}

\def\rtz{{\R/T\Z}}
\def\rt'z{{\R/T'\Z}}

\def\End{{\text{\rm End}}}
\def\txd{{\text{\rm d}}}

\def\dhjr{{\d_{\scriptscriptstyle{\text{\rm HJR}}}}}
\def\Inf{{\text{\rm Inf}}}

\def\Qm{{Q_\txm}}

\def\dfn={{\overset{\text{\rm def}}\to=}}
\def\BBig(-{{\Big(\!\!-}}
\def\BBiglangle-{{\Big\langle\!\!-}}

\def\:{{\text{\rm:}}}
\def\;{{\text{\rm;}}}
\def\M{{\text{\rm M}}}
\def\SL{{\text{\rm SL}}}

\def\twopi{{2\pi}}
\def\twopii{{\twopi\txti}}

\def\sumkton{{\sum_{k=1}^n}}

\def\sumkztonone{{\sum_{k=0}^{n+1}}}

\def\sumdionetop{{\sum_{i=1}^p}}

\def\sumdjonetonmone{{\sum_{j=1}^{n-1}}}
\def\sumdjtwotonmone{{\sum_{j=2}^{n-1}}}

\def\sumdij{{\sum_{i,j}}}

\def\sumdijk{{\sum_{i< j< k}}}

\def\sumddijk{{\displaystyle\sumd{i,j<k}}}

\def\sumd#1{{\sum_{#1}}}

\def\moneuj{{(-1)^{j}}}

\def\moneuk{{(-1)^{k}}}

\def\AS{{\text{\rm AS}}}

\def\moneuk{{(-1)^k}}
\def\moneun{{(-1)^n}}
\def\moneunone{{(-1)^{n+1}}}

\def\sign{{\text{\rm sign}}}

\def\pidz{{\pi_0}}

\def\pidzs{{\pi_0^*}}

\def\pidnone{{\pi_{n+1}}}

\def\pidk{{\pi_k}}

\def\sigus{{\sig^*}}

\def\adz{{a_0}}

\def\adone{{a_1}}

\def\adoneuk{{a_1^k}}

\def\adtwo{{a_2}}

\def\adi{{a_i}}

\def\adj{{a_j}}

\def\adk{{a_k}}

\def\adoneone{{a_{11}}}
\def\adonetwo{{a_{12}}}

\def\adtwoone{{a_{21}}}
\def\adtwotwo{{a_{22}}}

\def\bda{{b_a}}

\def\bdj{{b_j}}

\def\Bd#1{{B_{#1}}}

\def\Bdij{{\Bd{ij}}}
\def\Bdii{{\Bd{ii}}}
\def\Bdji{{\Bd{ji}}}

\def\Bdjk{{\Bd{jk}}}
\def\Bdkj{{\Bd{kj}}}
\def\Bdkk{{\Bd{kk}}}
\def\Bdik{{\Bd{ik}}}
\def\Bdki{{\Bd{ki}}}
\def\bd#1{{\bd_{#1}}}

\def\bd#1{{b_{#1}}}

\def\bdone{{b_1}}

\def\bdtwo{{b_2}}

\def\bdi{{b_i}}

\def\cda{{c_a}}

\def\cdb{{c_b}}
\def\cdbuij{{c_b^{i,j}}}
\def\cdbuii{{c_b^{i,i}}}

\def\cdtxa{{c_{\text {\rm a}}}}

\def\cdtxs{{c_{\text {\rm s}}}}

\def\edij{{e_{i,j}}}

\def\ed#1{{e_{#1}}}
\def\edkell{{\ed{k,\ell}}}
\def\edjk{{\ed{j,k}}}
\def\edkj{{\ed{k,j}}}
\def\edik{{\ed{i,k}}}

\def\mdz{{m_0}}

\def\mdj{{m_j}}

\def\mdij{{m_{i,j}}}
\def\mdji{{m_{j,i}}}

\def\elldone{{\ell_1}}

\def\elldk{{\ell_k}}

\def\kdone{{k_1}}
\def\kdn{{k_n}}
\def\kdnmone{{k_{n-1}}}

\def\kdtwo{{k_2}}

\def\kdn{{k_n}}

\def\rdone{{r_1}}

\def\rdij{{r_{i,j}}}
\def\rdji{{r_{j,i}}}

\def\r'dij{{r'_{i,j}}}
\def\r'dji{{r'_{j,i}}}
\def\pdz{{p_0}}
\def\pdone{{p_1}}

\def\pdtwo{{p_2}}

\def\pdthr{{p_3}}

\def\pdk{{p_k}}
\def\pdr{{p_r}}

\def\qdone{{q_1}}
\def\qdtwo{{q_2}}
\def\qdthr{{q_3}}
\def\tqdone{{\tilde q_1}}
\def\tqdtwo{{\tilde q_2}}
\def\tqdthr{{\tilde q_3}}
\def\qdi{{q_i}}
\def\qdj{{q_j}}

\def\rdone{{r_1}}

\def\rdi{{r_i}}

\def\pdi{{p_i}}

\def\pdj{{p_j}}

\def\pdr{{p_r}}
\def\edz{{e_0}}
\def\edone{{e_1}}

\def\edtwo{{e_2}}

\def\edthr{{e_3}}

\def\edfour{{e_4}}

\def\tedz{{\tilde e_0}}
\def\tedone{{\tilde e_1}}
\def\tedtwo{{\tilde e_2}}

\def\fdz{{f_0}}
\def\fdone{{f_1}}
\def\fdtwo{{f_2}}

\def\sdtxo{{s_{\text{\rm o}}}}

\def\sdone{{s_1}}
\def\sdtwo{{s_2}}
\def\sdthr{{s_3}}
\def\sdi{{s_i}}

\def\udi{{u_i}}

\def\udij{u_{i,j}}
\def\udji{{u_{j,i}}}

\def\vdi{{v_i}}

\def\zd#1{{z_{#1}}}
\def\zdi{{\zd{i}}}
\def\dotzdi{{\dot z_i}}

\def\zd#1{{z_{#1}}}

\def\vdone{{v_1}}

\def\bar xdone{{{\Bar x}_1}}

\def\xdij{{x_{ij}}}

\def\xdij{{x_{i,j}}}

\def\ydij{{y_{i,j}}}

\def\zdz{{z_0}}
\def\zdzuk{{z_0^k}}
\def\zdzu#1{{z_0^{#1}}}
\def\zdzukdn{{\zdzu{\kdn}}}
\def\zdzukdnmone{{\zdzu{\kdnmone}}}

\def\zdzukdone{{z_0^{k_1}}}
\def\zdzukdtwo{{z_0^{k_2}}}
\def\zdzukdthr{{z_0^{k_3}}}
\def\zdzukdj{{z_0^\kdj}}
\def\zdzukdjone{{z_0^\kdjone}}
\def\zdzukdnmone{{z_0^\kdnmone}}

\def\zdzun{{z_0^n}}

\def\zdzukone{{z_0^\kdone}}

\def\zdzukn{{z_0^{k_n}}}

\def\zdone{{z_1}}

\def\zdtwo{{z_2}}

\def\zdk{{z_k}}
\def\zdm{{z_m}}

\def\zdr{{z_r}}

\def\bzdz\inv{{{\Bar z}_0^{-1}}}

\def\udone{{u_1}}

\def\txs{{\text{\rm s}}}

\def\Zinf{{\Z^{<\N}}}
\def\sdij{{s_{i,j}}}
\def\sdji{{s_{j,i}}}

\def\lcm#1#2{{\left(\frac {#1#2}{\gcd({#1},{#2})}\right)}}
\def\lcm#1#2#3{{\left(\frac
{#1#2#3}{\gcd({#1},{#2}, {#3})}\right)}}

\def\zdzu#1{{z_0^{#1}}}

\def\zd#1{{z_{#1}}}

\def\elldi{{\ell_i}}

\def\lad#1#2{{\la_{#1}^{#2}}}

\def\lad#1{{\la_{#1}}}

\def\detdijk{{\det\!{}_{ijk}}}
\def\detdjik{{\det\!{}_{jik}}}
\def\detdkij{{\det\!{}_{kij}}}

\def\explr#1{{\exp\!\lr{#1}}}

\def\explrtwopii#1{{\exp\!\left(\twopii\lr{#1}\right)}}

\def\explrtwopiim#1{{\exp\!\left(-\twopii{#1}\right)}}

\def\lr#1{{\left({#1}\right)}}
\def\lrbrace#1{{\left\{#1\right\}}}
\def\lrbracket#1{{\left[{#1}\right]}}
\def\lrbracketd#1#2{{\left[{#1}\right]_{#2}}}

\def\sumddijk{{\sum_{i,j<k}}}
\def\sumd#1{{\sum_{#1}}}

\def\sumdij{{\sumd{i,j}}}

\def\Qdm{{Q_{\txm}}}

\def\Gdm{{G_{\txm}}}
\def\Hdm{{H_{\txm}}}

\def\ot#1#2{{{#1}\otimes{#2}}}

\def\aaijk{{a(i,j,k)}}

\def\aajik{{a(j,i,k)}}

\def\aakij{{a(k,i,j)}}
\def\aakik{{a(k,i,k)}}

\def\aaiik{{a(i,i,k)}}

\def\bbij{{b(i,j)}}
\def\bbiz{{b(i,0)}}
\def\bbzi{{b(0,i)}}
\def\bbzj{{b(0,j)}}
\def\bbjz{{b(j,0)}}

\def\cck{{c(k)}}
\def\ccij{{c(i,j)}}

\def\ccii{{c(i,i)}}

\def\ccjk{{c(j,k)}}

\def\ffdijk{{f_{i,j,k}}}

\def\ffdkij{{f_{k,i,j}}}

\def\ffdjik{{f_{j,i,k}}}

\def\ot{{\otimes}}

\def\Bdetaz{{B_{\eta, \z}}}

\def\Asaijk{{(\AS a)(i,j,k)}}

\def\ffdc{{f_c}}

\def\Ndz{{\N_0}}
\def\Ndzutw{{\N_0^2}}
\def\Ndzuthr{{\N_0^3}}

\def\xxdaijk{{X_a(i,j,k)}}
\def\xxdasaijk{{X_{\scriptscriptstyle{\AS a}}(i,j,k)}}
\def\xxdasa{{X_{\scriptscriptstyle{\AS a}}}}

\def\xxdaik{{X_a(i,k)}}
\def\yydaijk{{Y_a(i,j,k)}}

\def\yydaik{{Y_a(i,k)}}

\def\zzdaijk{{Z_a(i,j,k)}}

\def\zzdaik{{Z_a(i,k)}}
\def\xxda{{X_a}}
\def\yyda{{Y_a}}
\def\zzda{{Z_a}}
\def\uudaijk{{U_a(i,j,k)}}

\def\uudaik{{U_a(i,k)}}
\def\vvdaijk{{V_a(i,j,k)}}

\def\vvdaik{{V_a(i,k)}}

\def\uuda{{U_a}}
\def\vvda{{V_a}}

\def\lcm{{\text{\rm  lcm}}}

\def\tzzdbij{{\tZ_b(i,j)}}
\def\tzzdbii{{\tZ_b(i,i)}}

\def\tzzda{{\tZ_a}}

\def\tzzdaijk{{\tZ_a(i,j,k)}}
\def\tzzdaik{{\tZ_a(i,k)}}

\def\tbbbda{{\tB_a}}

\def\tbbbdaijk{{\tB_a(i,j,k)}}

\def\Lada{{\La_a}}
\def\Ladb{{\La_b}}
\def\Ladbij{{\La_b(i,j)}}
\def\Ladbii{{\La_b(i,i)}}

\def\tbbdb{{\tB_b}}
\def\tbbdbii{{\tB_b(i,i)}}

\def\tbbdbij{{\tB_b(i,j)}}

\def\tbbda{{\tB_a}}
\def\tbbdb{{\tB_b}}

\def\tbbda{{\tB_a}}

\def\nudb{{\nu_b}}

\def\dotadi{{\dot a_i}}

\def\Deldz{{\Delta_0}}
\def\edietan#1#2{{\etappi\!\lr{{\bracketpdi{\edi(#1)}};
{\bracketpdi{\edi(#2)}}}}}
\def\edjetan#1#2{{\etappj\!
\lr{{\bracketpdj{\edj(#1)}};
{\bracketpdj{\edj(#2)}}}}}
\def\edketan#1#2{{\etappk\lr{{\bracketpdk{\edk(#1)}};
{\bracketpdk{\edk(#2)}}}}}

\def\fsdm{{\fs_\txm}}

\def\zdj{{z_j}}
\def\dottxmdi{{\dot\txm_i}}

\def\gcdpijk{{\gcd(\pdi, \pdj,\pdk)}}

\def\RuDelt{{\R^\Delta}}

\def\ndij{{n_{i,j}}}

\def\mudauijk{{\mu_a^{ijk}}}
\def\mudauik{{\mu_a^{ik}}}
\def\ddik{{\d_{ik}}}

\def\ddjell{{\d_{j\ell}}}
\def\lrbrace#1{{\left\{#1\right\}}}

\def\whaa{{\widehat a}}
\def\whaaijk{{\widehat a(i,j, k)}}
\def\whaajik{{\widehat a(j,i, k)}}
\def\whaakij{{\widehat a(k,i,j)}}
\def\whaaiik{{\widehat a(i,i, k)}}
\def\whaakik{{\widehat a(k,i, k)}}

\def\whaaijj{{\widehat a(i,j,j)}}
\def\pidt{{\pi_T}}
\def\pidt'{{\pi_{T'}}}
\def\prdone{{\pr_1}}

\def\chidz{{\chi_{\scriptscriptstyle{0}}}}
\def\tzdb{{\tZ_b}}

\def\tzdbij{{\tZ_b(i,j)}}

\def\pidij{{\pi_{ij}}}

\def\pidii{{\pi_{ii}}}

\def\pidiui{{\pi_i^i}}
\def\pidijui{{\pi_{i,j}^i}}
\def\pidijuj{{\pi_{i,j}^j}}

\def\pidkiui{{\pi_{ki}^i}}

\def\pidk{{\pi_{k}}}

\def\\dotpidk{{\dot\pi_{k}}}

\def\gcdlr#1{{\gcd\!\lr{#1}}}

\def\ladbij{{\la_b(i,j)}}
\def\ladbki{{\la_b(k,i)}}

\def\wdij{{w_{i,j}}}
\def\wdji{{w_{j,i}}}
\def\zdaijk{{z_a(i,j,k)}}

\def\sRdz{{\sR_0}}
\def\sRdzone{{\sR_{0,1}}}

\def\sdt{{s_{T}}}
\def\sdt'{{s_{T'}}}
\def\pdtxo{{p_{\text{\rm o}}}}

\def\Gallr#1{{\text{\rm Gal}\lr{#1}}}
\def\eum#1{{e^{-{#1}}}}
\def\tup{{T'}}
\def\simdc{{\sim_c}}
\def\appfindimf{{\text\allowlinebreak{approximately finite dimensional factor}}}
\def\appfindiml{{\text\allowlinebreak{approximately finite dimentsional}}}
\def\rmhjrsq{{\text\allowlinebreak{reduced modified HJR-sequence}}}
\def\Rmhjrsq{{\text\allowlinebreak{Reduced Modified HJR-Sequence}}}
\def\tadzdz{{\ta_\zdz}}

\def\afdf{{\text\allowlinebreak{approximately finite dimensional factor}}}
\def\Zda{{Z_a}}
\def\Zdb{{Z_b}}
\def\Bda{{B_a}}
\def\Bdb{{B_b}}
\def\ddij{{D(i,j)}}
\def\fsfnt{{\text\allowlinebreak{faithful semi finite normal trace}}}
\def\pidtp{{\pi_{T'}}}
\def\Mapgsa{{\Map(G, \sA)}}
\def\Mapdzgsa{{\Map_0(G, \sA)}}
\def\tsA{{\widetilde \sA}}
\def\tc{{\tilde c}}
\def\bc{{\Bar c}}


\nologo
\TagsOnRight 
\loadeusm                    
\document 
\topmatter
\title{Outer Actions of a Discrete Amenable Group\\ on Approximately Finite Dimensional\\ Factors \threee, \\
The Type \threel\ Case, $\pmb{0<\la<1}$,\\ Examples }
\endtitle
\author{Yoshikazu Katayama and Masamichi Takesaki}
\endauthor
\dedicatory{Dedicated to the Memory of Masahiro Nakamura }
\enddedicatory
\address{YK: Division of Mathematical Sciences,
Osaka Kyoiku University,
Asahigaoka 4-698-1 Kashihara Osaka, 582-8582 Japan.}
\endaddress
\email{ katayama\@cc.osaka-kyoiku.ac.jp}
\endemail
\address{MT: Department of Mathematics, UCLA, PO Box 951555, Los Angeles, California 90095-1555.  Mail Address: 3-10-39 Nankohdai, Izumi-ku, Sendai, 981-8003 Japan}
\endaddress
\email{mt\@math.ucla.edu}
\endemail
\leftheadtext{Outer Actions \threee}
\rightheadtext{Y. Katayama and M. Takesaki}

\abstract
{In this last article of the series on outer actions of a countable dicrete amenable group on AFD factors, we analyze  outer actions of a countable discrete free abelian group on an AFD factor of type $\text{I\!I\!I}_\lambda$,  $0<\lambda< 1$, and compute outer conjugacy invariants. As a byproduct, we discover the asymmetrization technique for coboundary condition on a $\Bbb T$-valued cocycle of a torsion free abelian group, which might have been known by the group cohomologists. As the asymmetrization technique gives us a very handy criteria for coboundaries, we present it here in detail in the second section.}
\endabstract

\toc\noindent
{\S 0. Introduction.}\newline
{\S1. Simple examples and model construction.}\newline
{\S2. Asymmetrization.}\newline
{\S3. Universal resolution for a countable discrete abelian 
group.}\newline
\S4. The characteristic cohomology group $\La(\Hdm, L, M, \T)$.
\S5. The {\rmhjrsq.}\newline
\S6. Concluding remark.
\endtoc

\endtopmatter
\document

\head{\bf \S0. Introduction.}
\endhead
This article concludes the series of our joint work, \cite{KtT1},\cite{ KtT2}, and \cite{KtT3},
on the outer conjugacy
classification of outer actions of a {\cdag} on an approximately finite dimensional, (abbreviated to AFD),  factor by
examining outer actions of a countable discrete abelian group $G$  on an AFD factor $\sRdl$ of type {\threel},
$0<\la<1$. The cocycle conjugacy classification theory of actions of a  \cdag\ on an AFD factor had been completed through the work  of many mathematicians over three decades, \cite{Cnn3, Cnn4, Cnn6, Cnn7, Jn, JT, Ocn,  KtST1,
KtST2, KwST, ST1, ST2, ST3}, prior to the outer conjugacy classification theory.

Unlike the general classification program in operator algebras, the outer conjugacy classification of a {\cdag}
on  $\sRdl$ is almost smooth as shown in the series of previous work, see \cite{KtT3}. Only non-smooth part of the
classification theory stems from the classification of subgroups $N$ of $G$: for instance the classification of subgroups
of a torsion free abelian group of higher rank is non-smooth. We refer the work of Sutherland concerning Borel parameterization of polish groups, \cite{St2}. When  the modular automorphism part $N=\da\inv\lr{\cntr(\sM)}$ of the outer action $\dot\a$ of $G$ on $\sRdl$, is fixed,  the set of invariants becomes
a compact abelian group. It is a rare case in the theory of operator algebras. So we are encouraged to make a concrete
analysis of outer conjugacy class of  a {\cdag}. Of course, without having a concrete date on the group $G$ involved,
we cannot make  a fine analysis. So we take a {\cdfabg} $G$ and study its outer actions on $\sRdl$ and identify the
invariants completely. The justification of this restriction rests on the fact that all outer actions of  a {\cdabg} $A$ can be
viewed as outer actions of $G$ by pulling back the outer action via the quotient map: $G\mapsto A$.  Thanks to all
hard analytic work on the cocycle conjugacy classification in the past,  cited in the reference, our work is very algebraic
and indeed done by cohomological computations. 

We will begin first by relating the discrete core of $\sRdl$ and the core of an AFD factor $\sRdone$ of type {\threeone}.  
This analysis will enable us to have a simple model construction with given invariants, which is presented here in Section 1.
Single automorphisms and a pair of commuting automorphisms of $\sRdl$ are studied first. Then we will work on the
asymmetrization of a cocycle of a {\cdabg} which will provide a powerful tool for analysis of the third cohomology group $\ththr(G,
\T)$. The general theory of group cohomology is available to us today, for example see \cite{Brw}. But we have to work with individual
cocycles to analyze outer actions. So we have to have a tool to work with a cocycle directly beyond the computation  of the cohomology
group. For example, we have to identify which data of a given cocycle contributes to the modular automorphism part of the action in
question. Thus we will work on the cohomology group based on a very primitive method of chasing cocycles, through which we
discover the asymmetrization technique which provides us a quite handy criterion for the coboundary condition on a
cocycle of a torsion free abelian group. In our previous work, \cite{KtT1, KtT2, KtT3}, the outer conjugacy classification of a {\cdag} outer actions
were studied by a resolution of the relevant third cocycle. In the abelian case, it is shown that there is a universal
resolution group which takes care of all third cocycles at once which simplifies greatly the investigation of outer actions
of a \cdabg. The  {\rmhjrsq} will provide us a tool to chase the cocycles along with the asymmetrization technique. The first step of studying outer actions of {\cdabg} $G$ on a factor $\sM$ of type \threel, $0<\la<1$, is to find a {\cdag} $H$ and
a surjective homomorphism $\pig: H\mapsto G$ so that the pull back $\pigs(c)$ is a coboundary,  the process called the
resolution of a cocycle $c\in\tzthr(G, \T)$. Then the outer action $\da$ is  identified with a lifting $\shs(\a)$ of an action
$\a$ of $H$ through a cross-section $\sh: G\mapsto H$ of the homomorphism $\pig$. Luckily, a {\cdabg} $G$ admits a
universal resolution $\lrbrace{H, \pig}$, a group $H$ and a surjective homomorphism $\pig: H\mapsto G$ such that
$\pigs\lr{\ththr(H, \T)}=\{1\}$. The group $H$ is constructed via relatively simple process from a {\cdfabg} $G$. This
makes us possible to reduce the study of an outer action $\da$ of $G$ to that of an action $\a$ of $H$. Now, the action $\a$ of $H$
does not  lift to the discrete core $\tMd$ if $\modd(\a)\neq 1$. So we construct a central extension $\Hdm$ of $H$:
$$\CD
0@>>>\Z@>n\to \zdzun>>\Hdm@>>>H@>>>1
\endCD
$$
and work with the characteristic cohomology group $\La(\Hdm, L, M, \T)$ where the normal subgroup $L$ stands
for the inverse image
$L=\piginv(N)$ with $N=\da\inv\lr{\cntr(\sM)}$. Thus we are going to investigate the {\rmhjrsq}:
$$\CD
\thtw(H, \T)@>\Res>>\La(\Hdm, L, M, \T)@>\d>>\thmsout(G, N, \T)@>\Inf >>\ththr(H, \T)\\
@|@V\pims VV@V\partgm VV@|\\
\thtw(H, \T)@>\res>>\La(H, M, \T)@>\dhjr>>\ththr(G, \T)
@>\pigs>>\ththr(H, \T)
\endCD
$$
Here $\fs$ is a fixed cross-section of the quotient map: $G\mapsto Q=G/N$. 
The groups appeared on the above exact sequences are all compact abelian groups and are indeed computable as shown in this paper.

We refer \cite{Brw, EMc, McWh, Hb, Jn} for the general cohomology theory of abstract groups and \cite{St1} for the cohomology theory related to \vnas. We refer  \cite{Tk1, Tk2, Tk3} for the general theory of \vnas. Concerning the discrete core of a factor of type \threel, we refer \cite{Cnn1, Cnn2, CT, FT1 and FT2}. 

This work was originated from the authors' visit to the Erwin Shr\"odinger Institute, Vienna,  and the University of Rome, La Sapienza, in the spring of 2005 and further developed throughout the subsequent years. The second named author visited the Erwin Shr\"odinger Institute in the fall of 2008 again where the final touch on the joint work was given. The authors are greatly indebted to these institute, in particular to Professors Klaus Schmidt and Sergio Doplicher who made our collaboration possible and pleasant. We would like to record here our sincere appreciation to their support and hospitality.

\head{\bf \S1 Simple Examples and Model Construction.}
\endhead

\subhead\nofrills{Factors  of Type \threel\ and Type \threeone, and Their Cores:}
\endsubhead\quad 
We begin by the following folklore theorem in the structure theory of factors of type \threee.
\proclaim{Theorem 1.1} Let $\lrbrace{\sMdzone, \tau, \th}$ be a factor of type {\twoinf} equipped with {\fsfnt} $\tau$ and trace scaling automorphism $\th$ by $\la, 0<\la<1$, i.e., $\tau\scirc\th=\la\tau$ and let $\sM=\sM_{0,1}^\th$ be the fixed point subalgebra of $\sMdzone$ by $\th$. Then we have the following statements\:
\roster
\item"i)" The {\vna} $\sM$ is a factor of type \threel\;
\item"ii)" The triplet $\lrbrace{\sMdzone, \tau, \th}$ is conjugate to the discrete core of $\sM$\;
\item"iii)" For an automorphism $\a\in \Aut\lr{\sMdzone}$, the following statements are equivalent\:
\itemitem {\rm a)} $\a(\sM)=\sM;$
\itemitem{\rm b)} $\a\scirc \th=\th\scirc \a$.
\item"iv)" Let $\Aut(\sMdzone, \sM)$ be the group of automorphisms of $\Aut\lr{\sMdzone}$ leaving $\sM$ globally invariant. Then we have the following exact sequence\:
$$\CD
0@>>>\Z@>n\to \th^n>>\Aut(\sMdzone, \sM)@>\a\to \a|_\sM>>\Aut(\sM)@>>>1;
\endCD
$$
\item"v)" The subgroup $\lrbrace{\th^n: n\in \Z}$ is the Galois group of the pair $\lrbrace{\sMdzone, \sM}$ in the sense that
$$
\lrbrace{\th^n: n\in \Z}=\lrbrace{a\in \Aut(\sMdzone): \a(x)=x, x\in\sM}.
$$
\item"vi)" If $\a\in \Aut\lr{\sMdzone, \sM}$, then the modulus
$\modd_\sMdzone(\a)$ as a member of $\Aut\lr{\sMdzone}$ gives the
modulus $\modd_\sM(\a)$ of the restriction $\a|_\sM\in\Aut(\sM)$ in the following way\:
$$
\modd_\sM(\a)=\pidtp\lr{\modd_\sMdzone(\a)}\in\rt'z,
$$
where  
$$\gathered
T'=-\log\la, \quad T=\frac \twopi {T'},\\
\pidtp: s\in\R\mapsto \dstp=s+T'\Z\in\rt'z.
\endgathered
$$
\endroster
\endproclaim
\demo\nofrills{Proof.}\quad The statements (i) and (ii) are known in the general structure theory of a factor of type \threee, see \cite{Tk2, Chapter XII, \S2 and \S6}.

v) We prove the statement (v) first.  Let $\p$ be a generalized trace of $\sM$, i.e., a {\fwt} on $\sM$ such that $\p(1)=+\infty$ and $\sigpdT=\id$. 
Then the covariant system $\lrbrace{\sMdzone, \th}$ is conjugate to the dual system
$\lrbrace{\sM\rtimes_\sigp \rtz, \Z, \widehat \sigp}$. So we identify them, so that $\sMdzone$ admits a periodic one parameter 
unitary group $\lrbrace{u^\p(s): s\in \R}$:
$$\gathered
u^\p(T)=1,\quad  \th(u^\p(s))=\la^{is}u^\p(s),\quad \text{and}\quad
\sigps=\Ad(u^\p(s))|_\sM, \quad s\in\R.
\endgathered
$$
Furthermore, the one parameter unitary group  $\lrbrace{u^\p(s): s\in \R}$ together with $\sU(\sM)$ generates the normalizer $\tsUdz(\sM)=\lrbrace{v\in\sU\lr{\sMdzone}: v\sM v^*=\sM}$, giving the semi-direct product decomposition $\tsUdz=\sU(\sM)\rtimes_\sigp \rtz$. Suppose that $\a\in\Aut\lr{\sMdzone}$ leaves $\sM$ pointwise fixed. We then show that $x\in\sM$ and
$u^\p(s)^*\a\lr{u^\p(s)}, s \in \R$, commute:
$$\aligned
x&u^\p(s)^*\a\lr{u^\p(s)}=u^\p(s)^*u^\p(s)xu^\p(s)^*\a\lr{u^\p(s)}\\
&=u^\p(s)^*\sigps(x)\a\lr{u^\p(s)}=
u^\p(s)^*\a\lr{\sigps\lr{x}u^\p(s)}\\
&=u^\p(s)^*\a\lr{u^\p(s)x}=u^\p(s)^*\a\lr{u^\p(s)}x,
\endaligned
$$
so that $u^\p(s)^*\a\lr{u^\p(s)}=\sMdzone\cap\sM'=\C$. Hence there exists a scalar $\mu(s)\in \T$ such that
$$
\a\lr{u^\p(s)}=\mu(s)u^\p(s), \quad s\in\R.
$$
Since $u^\p(T)=1$, we have $\mu(T)=1$. Since $\mu(s+t)=\mu(s)\mu(t), s, t\in\R$, we have
$$
\mu(s)=\la^{\txti ns}, s \in \R,\quad \text{for some }\quad n\in\Z.
$$
Since $\sM$ together with $\lrbrace{u^\p(s): s\in\R}$ generate the whole algebra $\sMdzone$, we conclude that $\a=\th^n$. 
This shows (v).

iii) Suppose that $\a\in\Aut\lr{\sMdzone}$ leave $\sM$ globally invariant. Let $\b=\a_\sM=\a|_\sM$ be 
the automorphism of $\sM$ obtained as the restriction of $\a$ to $\sM$. Then the uniqueness of a generalized trace on $\sM$ gives a 
scalar $s\in\R$ and a unitary $v\in\sU(\sM)$ such that 
$$
e^{-s}\p=\p\scirc \lr{\Ad (v)\scirc \b}.
$$
This means that 
$$
\modd(\b)=\modd(\Ad(v)\scirc \b)=\dstp=s+T'\Z\in
 \rt'z,
$$
and that  $\sigp$ and $\Ad(v)\scirc \b$ commute. Hence it is possible to extend $\Ad(v)\scirc \b$ to the automorphism 
$\g\in\Aut\lrbrace{\sMdzone}$ such that
$$
\g\lr{u^\p(t)}=u^\p(t), \quad t\in \R, \quad \g(x)=\Ad(v)\scirc \b(x), \quad x \in \sM.
$$
Now we compare $\a$ and $\g$ on $\sM$: 
$$
\g(x)=\Ad(v)\scirc \b(x)=\Ad(v)\scirc \a(x), \quad x\in\sM.
$$
From (v) it follows that $\a$ is of the form:
$$
\a=\th^n\scirc \Ad(v^*)\scirc \g.
$$
for some $n\in\Z$. Since $\th$ commutes with both $\g$ and $\Ad(v)$, $\a$ and $\th$ commute. Hence the implication $(a)\Rightarrow (b)$ follows. 
The reversed implication: $(b)\Rightarrow (a)$ is trivial. So the proof of (iii) is complete. 

iv)  This follows from (iii) and (v).

Let $\lrbrace{\tM, \R, \tau, \th}$ be the non-commutative flow of weights on $\sM$ so that the covariant system $\lrbrace{\sMdzone,
\Z, \th}$ is identified with $\lrbrace{\sM\vee \lrbrace{\p^{}}, \th_{T'}}$

vi) Fix a member $\a\in \Aut\lr{\sMdzone, \sM}$ and let
$\txm(\a)=\modd(\a)\in\R$ so that
$$
\tau\scirc \a=e^{-\txm(\a)}\tau.
$$
Consider the crossed product 
$$
\tM=\sMdzone\rtimes_\th \Z\cong \sM\botimes\sL\lr{\ell^2(\Z)}
$$
and the generalize trace $\f=\tau\scirc \sE$ on $\tM$:
$$
\f(x)=\tau\scirc \sE(x)=\tau\lr{\int_\rtz\hat\th_s(x)\txd s},
\quad x\in\tM_+.
$$
With $U\in\sU(\tM)$ the unitary corresponding to the crossed product $\sMdzone\rtimes_\th \Z$, we extend $\a$ to $\ta\in\Aut\lr{\tM}$ by:
$$
\ta(x)=\a(x), \quad x\in \sMdzone, \quad \ta(U)=U.
$$
Then we have for each $x\in\tM_+$
$$\aligned
\f\scirc \ta(x)&=\tau\lr{\int_\rtz\hat\th_s\lr{\ta(x)}\txd s}=\tau\lr{\a\lr{\int_\rtz\hat\th_s(x)\txd s}}\\
&=e^{-\txm(\a)}\tau\lr{\int_\rtz\hat\th_s(x)\txd s}\\
&=e^{-\txm(\a)}\f(x).
\endaligned
$$
Hence we get
$$
\modd(\ta)=\lrbracket{\txm(\a)}_{T'}=\txm(\a)+T'\Z\in\rt'z.
\tag1.1
$$
Since the covariant systems $\lrbrace{\sM, \a}$ and $\lrbrace{\tM, \ta}$ are cocycle conjugate, we have
$\modd(\ta)=\modd(\a)$. 
This completes the proof.
\QED

\enddemo

Now, we denote by $\sRdz$ an approximately finite dimensional factor of type \twoone\ throughout the paper.

A factor $\sMdone$ of type \threeone\ generates one parameter family $\lrbrace{\sMdl: 0<\la\leq 1}$ of factors of type \threel, 
who share the same discrete core $\sMdzone$. So let $\sMdone$ be a
factor of type \threeone, and $\lrbrace{\sMdzone, \thds, s\in\R}$  be the
non-commutative flow of weights on $\sMdone$, i.e., $\sMdzone$ is a
factor of type \twoinf\ equipped with a trace scaling one parameter 
automorphism group $\lrbrace{\thds: s\in\R}$ and a \fsfntr\ $\tau$ such
that
$$
\sMdone=\sMdzoneuth, \quad \tau\scirc\thds=e^{-s}\tau, \ s\in\R.
$$
The following is a folklore theorem in the structure theory of type \threee. 
\proclaim{Theorem 1.2} In the above context, 
fixing  $T'>0$, set
$$
\la=\eum{T'}, \quad T=\frac\twopi \tup,
$$
and let $\sMdl$ be
the fixed point subalgebra $\sMdzoneu{\thdtprime}$ of $\sMdzone$ under the automorphism $\thdtprime$. Then the following statements hold\:
\roster
\item"i)" The subalgebra $\sMdl\i \sMdzone$ is a factor of type \threel, whose discrete core is conjugate to the pair $\lrbrace{\sMdzone, \thdtprime}$.
\item"ii)" The triplet $\lrbrace{\sMdzone, \sMdl, \thdtprime}$ is a Galois triplet in the sense\:
$$\gathered
\Gallr{\sMdzone/ \sMdl}=\lrbrace{\thdtprimeun: n\in\Z},\\
\endgathered
$$
where
$$\gathered
\Gallr{\sM/ \sN}=\lrbrace{\a\in\Aut(\sM): \a|_\sN=\id}
\endgathered
$$
for any pair $\sN\i\sM$ of \vnas. Furthermore, we have the following exact sequence\:
$$\CD
1@>>>\lrbrace{\thdtprimeun: n\in\Z}@>>>
\Aut(\sMdl)_\txm@>>>\Aut(\sMdl)@>>>1
\endCD
$$
and 
$$\aligned
\Aut(\sMdl)_\txm&=\lrbrace{\ta\in\Aut\lr{\sMdzone}: \ta(\sMdl)=\sMdl}\\
&=\lrbrace{\ta\in\Aut\lr{\sMdzone}: \ta\scirc \th_{T'}=\th_{T'}\scirc \ta}.
\endaligned
$$

\item"iii)" Another pair $\lrbrace{\sMdl, \sMdone}$ forms a Galois pair\:
$$\gathered
\Gallr{\sMdl/\sMdone}=\lrbrace{\th_\dstp: \dstp=s+T'\Z\in\rt'z, s\in \R},
\endgathered
$$
i.e., an automorphism $\a\in\Aut(\sMdl)$ is of the form
$\a=\thdstp$ for some $\dstp\in\rt'z$ if and only if $\a(x)=x, x\in\sMdone$.
\item"iv)" The modulus of $\thdstp\in\Aut(\sMdl)$ is precisely $-\dstp\in\rt'z$ itself, i.e.,
$$
\modd\lr{\thdstp}=-\dstp\in\rt'z.
$$
\endroster
If any one of $\sMdl, \sMdone$ and $\sMdzone$ is approximately finite dimensional, then all others are approximately finite dimensional and 
 in addition the following statements hold\:
\roster
\item"v)" If  $\a\in\Aut(\sMdl)$ has aperiodic modulus $\txm=\modd(\a)$, i.e., if $k\txm\neq 0$ for every non-zero integer $k\in\Z$, or equivalently if
$$
\frac{\bracett'{\modd(\a)}}{T'}\not\in \Q,
$$
then $\a$ is cocycle conjugate to $\th_{-\txm}$.
\item"vi)" If an automorphism $\a\in\Aut(\sMdl)$ has trivial asymptotic outer period, i.e., $p_a(\a)=0$, then its cocycle conjugacy class is determined by its modulus $\txm=\modd(\a)\in\rt'z$. In fact, the automorphism $\a$ is cocycle conjugate to
the automorphism $\th_{-\txm}\ot \sigdz$ on $\sMdl\cong
\sMdl\botimes \sRdz$, where $\sigdz\in\Aut(\sRdz)$ is any aperiodic automorphism of the \appfindimf\ $\sRdz$. If $\txm\neq0$, then we have $\th_\txm\sim \th_\txm\ot\sigdz$.
\endroster
\endproclaim
\demo\nofrills{Proof.} \quad We present a proof  the statements  (v) and (vi). Choose an automorphism $\a\in\Aut(\sMdl)$ such that $\txm=\modd(\a)$ is aperiodic.  Let $\sRdz$ be an approximately finite dimensional (to be abbreviated to AFD afterward) factor of type \twoone\ realized as the infinite tensor product of two by two matrix algebras 
$$
\sRdz=\prod_{n\in\Z}{}^\otimes \lrbrace{M_n, \tau_n}
$$
relative to the normalized traces $\tau_n=\Tr/2$ on $M_n=\M(2,\C)$. Let $\sigdz$ be the Bernouille shift automorphism of $\sRdz$, i.e., the automorphism determined by the following:
$$
\sigdz\lr{\prod_{n\in\Z}{}^\otimes x_n}=\prod_{n\in\Z}{}^\otimes x_{n+1}.
$$
Then thanks to the grand theorem of Connes, \cite{Cnn6, Tk3, page 267}, $\a$ and $\a\otimes \sigdz$ are cocycle conjugate under the identification of $\sMdl$ and $\sMdl\ \botimes\ \sRdz$ because the asymptotic outer period $p_a(\a)$ of $\a$ is zero, $p_a(\a)=0$. The same is true for $\th_\txm$, i.e., $\thdm\ \simdc\ \thdm\otimes\sigdz$, where ``$\simdc$" means the outer conjugacy. Since $\modd(\a_1\otimes\a_2)=\modd(\a_1)+\modd(\a_2)$ on $\sMdl\ \botimes\ \sMdl\cong \sMdl$, we have
$$
\a\ \simdc\ \a\otimes\sigdz\ \simdc\ \a\otimes\thdm\otimes \th_{-\txm}\ \simdc\
\sigdz\otimes \th_{-\txm}\ \simdc\  \th_{-\txm}.
$$
This proves the statement (v).

vi) Suppose that $p\in\N$ is the period of $\txm\in \rt'z$, i.e., the smallest
non-negative integer with $p\txm=0$. We assume that $p\neq 0$. Let 
$$
\lrbrace{\edjk: 1\leq j,k\leq p}
$$
be the standard matrix units of the
$p\times p$-matrix algebra $\M(p; \C)$, and for each $n\in\N$
set $M_n=\M(p; \C), n\in\N$, and also consider the diagonal unitary 
$$
u_n=\sumdionetop \explrtwopii{\frac{i-1}p}\edii\in \text{\rm U}(p;\C)
\i M_n
$$
of order $p$, i.e., $u_n^p=1$. Now we identify the AFD factor
$\sRdz$ with the infinite tensor product:
$$
\sRdz=\prod_{n\in\N}{}^\otimes \lrbrace{M_n, \tau_n}, \quad
\tau_n=\frac 1n \Tr
$$
and let
$$
\sig_p=\prod_{n\in\N}{}^\ot \Ad(u_n)\in \Aut(\sRdz)\in\Aut(\sRdz).
$$
Then the automorphism $\sig_p$ has the properties:
$$\gathered
\sig_p^k\not\in \Int(\sRdz),\ \text{for}\ k=1, \cdots, p-1, \ \text{and}\
\sig_p^p=\id\\
\th_\txm\sim_c \th_\txm\ot \sig_p \quad \text{on}\quad \sMdl\cong 
\sMdl\botimes\sRdz,\\
\th_\txm\ot\th_{-\txm}\sim_c \id\ot\id\ot\sig_p\quad \text{on}\quad
\sMdl\botimes\sMdl\cong\sMdl\botimes\sMdl\botimes\sRdz,\\
\sigdz\ot\sig_p\sim_c \sigdz\quad \text{on}\quad
\sRdz\botimes\sRdz\cong\sRdz.
\endgathered
$$
If $\a\in \Aut(\sMdl)$ has the trivial asymptotic outer period $p_a(\a)=0$, then
the automorphism $\a$ has the properties:
$$\gathered
\a\simdc \sig_p\ot\a\quad \text{on}\quad \sMdl\cong 
\sRdz\botimes\sMdl, \\
\th_{\txm} \ot \a \sim_c \id\ot \sigdz\quad \text{on}\quad
\sMdl\botimes
\sMdl
\cong \sMdl\botimes \sRdz,\\
\th_{-\txm}\ot\sigdz\ \simdc\ \th_{-\txm}\ot\thdm \ot
\a\ \simdc\ \sig_p\ot\a\ \simdc\
\a
\endgathered
$$
under the isomorphisms:
$$
\sMdl\botimes \sRdz\cong \sMdl\botimes \sMdl\botimes \sMdl
\cong \sRdz\botimes\sMdl\cong\sMdl.
$$
This completes the proof.
\QED
\enddemo
Thus if $\modd(\a)$ is aperiodic, or $p_a(\a)=0$, then the grand
theorem of Connes \cite{Cnn6, Tk3, page 270}, identifies the cocycle
conjugacy class of $\a\in\Aut(\sMdl)$. But if $\modd(\a)$ has non trivial
period, and $\pdone=p_a(\a)\neq 0$, then the cocycle conjugacy class
of $\a$ involves algebraic invariants. For example, one has to consider
the extension of $\a$ to the discrete core $\tMdld$ on which $\a$ alone
cannot act. In fact, one has to consider a larger group $\Z^2$ than the
integer group $\Z$. So we continue to the next paragraph.

\subhead\nofrills{Invariants for Single Automorphisms: \quad}
\endsubhead
We consider a single automorphism of a factor $\sM$ of type \threel, which can be viewed as an action of 
the integer additive group $\Z$. As the integer group $\Z$ appears in many different roles, we denote it by $G=\Z$. Let $\adone$ be 
the generator of the group $G$ so that $G=\Z\adone$. Sometime, we view $G$ as a multiplicative group in which case $G$ becomes
$G=\lrbrace{\adoneuk: k\in\Z}$. Since the integer group is cohomologically trivial, i.e., $\thtw(G, \T)=\ththr(G, \T)=\lrbrace{1}$, there
is  no distinction between the cocycle conjugacy problem and the outer conjugacy problem of actions of $G$. Namely, an outer action
$\dot\a$ of $G$ comes always from an action $\a$ of $G$ and outer conjugacy of the outer action $\dot\a$ of $G$ is the same as the
cocycle conjugacy of the action $\a$ of $G$. Hence the obstruction $\Ob(\dot\a)$ of $\dot\a$ and the characteristic invariant $\chi(\a)$ of
$\a$ is handily identified. The same is true for the modular obstruction $\Obm(\dot\a)$ and the modular characteristic invariant
$\chi_\txm(\a)$. 

As the single automorphism cocycle conjugacy classification wasn't handled properly in our previous work, \cite{KtST1, KtST2}, and more importantly the presentation of a single  automorphism on a factor of type \threel\ in the book of the second named author \cite{Tk3} contains a minor mistake, we present it here in some detail.  

Since the case that the modulus $\txm=\modd(\a)$ is aperiodic, then the last theorem takes care of the cocycle conjugacy of $\a$, i.e., it must be cocycle conjugate to $\th_{-\txm}$. So we handle only the case that $\bracett'{\modd(\a)}$ is rational multiple of $T'$.

Suppose $\a\inv(\cntr(\sM))=\Z\bdone$ and $\bdone=\pdone \adone, \pdone\in\N$.

Choose a pair $\pdone, \qdone\in\N$  of positive integers
$\qdone<\pdone$  such that
$$
\txm=\frac\qdone\pdone T'+T'\Z\in\rt'z, \quad 0\leq \qdone<\pdone.
$$

Then we form a group extension: 
$$\gathered
\Gdm=\lrbrace{(g, s)\in G\times \R: g\txm=\dstp=s+T'\Z\in\rt'z},\\
\CD
0@>>>\Z@>k\to\lr{0, kT'}>>\Gdm@>\prdone>>G@>>>0.
\endCD
\endgathered\tag1.2
$$
Set
$$\gathered
\zdz=(0, T'),\quad \zdone=\lr{\adone, \bracett'{\txm}}, \\
\bdone=\pdone\zdone-\qdone \zdz, \quad N=\Z\bdone, \quad
\Qdm=\Gdm/N.
\endgathered\tag1.3
$$
The group $\Gdm$ is equipped with a distinguished homomorphism $\txk_\txm=\prdtwo$ to $\R$:
$$
\txk_\txm(g, s)=s\in\R, \quad (g, s)\in\Gdm.
\tag1.4
$$

Let $\piq: g\in \Gdm\mapsto \dot g\in \Qdm$ be the quotient map 
and further set
$$\gathered
\Ddone=\gcdlr{\pdone, \qdone}, \quad
\text{and}\quad \rdone=\frac\pdone\Ddone, \quad \sdone
=\frac\qdone\Ddone,
\endgathered\tag1.5
$$
and find a pair $\udone, \vdone\in\Z$ of integers such that
$$
1=\rdone\udone-\sdone\vdone, \quad \text{equivalently}\quad \Ddone=\pdone\udone-\qdone\vdone,
$$
which can be done through the Euclid algorithm. In the event that $\qdone=0$, the modulus $\txm$ is trivial, i.e.,
$\txm=0$ and $\Gdm=G\oplus \Z$. 
\proclaim{Theorem 1.3 (Invariants for a Single Automorphism with Periodic Modulus)}\quad In the case that $\pdone$ and $\qdone$ are both non-zero, we have the following statements with $\Ddone=\gcdlr{\pdone,\qdone}$\:
\roster
\item"i)" The pair $\lrbrace{\zdz, \zdone}$ is a free basis of $\Gdm$ so that every element $g\in\Gdm$ is written uniquely in the form\:
$$
g=\edz(g)\zdz+\edone(g)\zdone.
$$
\item"ii)" The group $\Gdm$ admits another free basis $\lrbrace{\wdz, \wdone}$ such that
$$
\bdone=\Ddone\wdone,
$$
and therefore
$$\gathered
N=\Ddone\Z\wdone, \quad \Qdm=\Z\dotwdz\oplus\Z\dotwdone,\\
 \Ddone\dotwdone=0\ \text{ in }\ \Qdm\cong \Z\oplus \Z_\Ddone,
\endgathered
$$
where the dotted notations indicate their images in the quotient group $\Qdm$.
\item"iii)" The character group $\widehat\Qdm$ of $\Qdm$ and the characteristic cohomology group $\La(\Gdm, N, \T)$ is identified under the correspondence\:
$$
\la_\chi(n\bdone; g)=\chi(\piq( g))^n, \quad g\in \Gdm,\ \chi\in\widehat\Qdm.
\tag1.6
$$
\item"iv)" The character group $\widehat\Qdm$ is given by the exact sequence\:
$$\CD
0@>>>\Z^2@>>>\R\oplus \lr{\frac 1\Ddone\Z}
@>\exp(\twopii\ \cdot)>>\T\oplus \Z_\Ddone=\widehat\Qdm@>>>0,
\endCD
$$
 which describes the characteristic cohomology group $\La(\Gdm, N, \T)$\:
$$
\La(\Gdm, N, \T)\cong \T\oplus \Z_\Ddone.
\tag1.7
$$
\endroster
If $\chi(\zdz)$ is a root of unity, then 
the outer period $\pdtxo(\a)$ of $\a$ is given as the product $\pdone\sdtxo$ with $\sdtxo\in\Z_+$ the smallest non-negative integer $s\in\Z_+$ such that
$1=\chi(\zdz)^s$. If $\chi(\zdz)$ is not a root of unity, then the corresponding automorphism $\a$ is aperiodic, i.e., $\pdtxo(\a)=0$.
\endproclaim
\demo\nofrills{Proof.}\quad i) Since $\prdone(\zdone)=\adone$ and $G$ is a free abelian group, the exact sequence (1.2) splits along with the cross-section: $m\in G\mapsto m\zdone\in \Gdm$.

ii)  We set
$$\gathered
\wdz=\udone\zdz-\vdone \zdone, \quad 
\wdone=-\sdone\zdz+\rdone\zdone.
\endgathered
$$
Since
$$
\zdz=\rdone\wdz+\vdone\wdone,\quad
\zdone=\sdone\wdz+\udone\wdone,
$$
the pair $\lrbrace{\wdz, \wdone}$ is a free basis of $\Gdm$ such that
$$\gathered
\Gdm=\Z\wdz+\Z\wdone,\quad \bdone=\Ddone\wdone, \quad N=\Ddone\Z\wdone,\\
\Qdm=\Gdm/N= \Z\dotwdz\oplus \Z\dotwdone,
\endgathered
$$
as we wanted.

iii) Since $\thtw(N, \T)=\lrbrace{1}$, the second cocycle part of a characteristic cocycle in $\tZ(\Gdm, N, \T)$ is taken to be trivial, 
so that the $\la$-part vanishes on $N$ and therefore it is a character of $\Gdm$ which vanishes on $N$ and factors through the quotient map
$\piq: \Gdm\mapsto \Qdm$. Thus it is of the form:
$$
\la(\bdone; g)=\chi\lr{\piq(g)}, \quad g\in \Gdm, \ \chi\in \widehat\Qdm.
$$ 

iv) It follows from (ii) that the character group $\widehat \Qdm$ is parameterized by $\R\oplus \lr{\frac 1\Ddone\Z}$:
$$
\chi_{x, y}(g)=\explrtwopii{x\fdz(g)+y\fdone(g)}, \quad
g=\fdz(g)\wdz+\fdone(g)\wdone\in\Gdm,
$$
with $(x, y)\in\R\oplus \lr{\frac1\Ddone\Z}$. This gives the exact sequence:
$$\CD
0@>>>\Z^2@>>>\R\oplus \lr{\frac1\Ddone\Z}@>(x,y)\to\chi_{x,y}>>\widehat\Qdm=\T\oplus{\Z_\Ddone}@>>>0
\endCD
$$
The assertion (iv) follows.
This completes the proof.\QED
\enddemo

\subhead\nofrills{Model Construction:\quad}
\endsubhead
Let $G$ be a fixed \cdag\ and $\lrbrace{H, \pig}$ be  a universal resolution group of the third cocycles of $G$, i.e, $\pig: H\mapsto G$ is a surjective homomorphisms such that
$$
\pigs\lr{\tzthr(G, \T)}\i \tbthr(H, \T).
$$
We require $H$ to be a \cdag. Let $M=\Ker(\pig)$. Fix a normal subgroup $N$ of $G$ and set $L=\piginv(N)$. 
With a fixed invariant homomorphism $\txm\in\Hom_G(N, \rt'z)$ such that $\Ker(\txm)\j N$, we use the notation $\txm$ for
$\txm\scirc \pig$ for short and  form a group extension
$\Hdm$:
$$\CD
0@>>>\Z@>>>\Hdm@>\pim>>H@>>>1,
\endCD
$$
where
$$\gathered
\Hdm=\lrbrace{(g, s)\in H\times \R: \txm(g)=\dstp=s+T'\Z\in\rt'z},\\
\pim(g, s)=g\in H, \quad \txk(g, s)=s\in\R, \quad (g, s)\in\Hdm.
\endgathered 
$$
Then we get the following \rmhjrsq:
$$\eightpoint\CD
\cdots@>>>\thtw(H, \T)@>\Res>>\La(\Hdm, L, M, \T)@>\d>>\thmsout(G, N,\T)@>>>1
\endCD
$$
Thus every modular obstruction cocycle $(c, \nu)\in\tzmsout(G, N, \T)$ is of the form: 
$$
(c, \nu)\equiv\d(\la, \mu)\ \mod\ \tbmsout(G, N, \T).
$$
Consequently the construction of an outer action $\dot\a$ of $G$ on an AFD factor $\sMdl$ of type \threel\ with $\Obm(\dot\a)=([c],
\nu)\in\thmsout(G, N, \T)$ is reduced to the construction of an action $\a^{\la,\mu}$ of $\Hdm$ such that 
$$\gathered
\lr{\a^{\la, \mu}}^{-1}\lr{\Int\lr{\sMdl}}\j M, \quad 
\lr{\a^{\la, \mu}}^{-1}\lr{\cnt\lr{\sMdl}}=L,\\
\chi\lr{\a^{\la, \mu}}=\lrbracket{\la, \mu}\in \La(\Hdm, L, M, \T), \\
\mod\!\!\lr{\a_g^{\la, \mu}}=\txm(\pig(g)),\quad g \in \Hdm.
\endgathered
$$
So fix a set of  invariants $(\la, \mu)\in \tZ\lr{\Hdm, L, M, \T}$ and $\txm\in\Hom_G(G, \rt'z)$ such that $\Ker(\txm)\j N$. We are
going to construct the model action $\a^{\la,\mu}$ of $\Hdm$  as follows: 
\roster
\item"Step I:" Let $X$ be a countable but infinite set on which $\Hdm$ acts freely from the left. In the case that $\Hdm$ is an infinite group, then we take $X$ to be $\Hdm$ itself and let $\Hdm$ act on it by the multiplication from the left. So the infinite set $X$ is only needed when $\Hdm$ is a finite group in which case $X$ can be taken to be  the product set $X=\Hdm\times \N$ and $\Hdm$ act on the first component by the left multiplication . 
Let $\lrbrace{M_x, x\in \Hdm}$ be the set of 2 by 2 matrix algebras $\M(2, \C)$ indexed by  elements $x\in X$. 
\item"Step \two:" Form the infinite tensor product 
$$
\sRdz=\prod_{x\in X}{}^{\otimes}\lrbrace{M_x, \tau_x}
$$
relative to the normalized trace
$$
\tau_x\pmatrix\adoneone&\adonetwo\\
\adtwoone&\adtwotwo\endpmatrix=\frac{\adoneone+\adtwotwo}2.
$$
Then let $\siguz$ be the Bernouille action  of   $\Hdm$ on $\sRdz$ which is determined by:
$$
\sigdguz\lr{\prod_{x\in X}\!{}^\otimes a_x}=\prod_{x\in X}\!{}^\otimes a_{gx}.
$$
\item"Step \threee:" Form the twisted partial crossed product of $\sRdz$ by $N$ relative to the second cocycle $\mu\in\tztw(N, \T)$ and the action $\siguz$:
$$
\sMdz= \sRdz\rtimes_{\siguz, \mu}N.
$$
Let $\lrbrace{U(m): m\in N}$ be the projective unitary representation of $N$ to $\sMdz$ corresponding to the twisted crossed product so that
$$\gathered
U(g)U(h)=\mu(g; h)U(gh), \quad g, h\in N;\\
U(g)a U(g)^*=\sigdguz(a), \quad a \in \sRdz, g \in N.
\endgathered
$$
Let $\sigulmu$ be the action of $\Hdm$ on $\sMdz$ determined by:
$$\gathered
\sig^{\la, \mu}_g(U(m))=\la(gmg\inv; g)U(gmg\inv), \quad m\in N, g\in \Hdm;\\
\sig^{\la, \mu}_g(a)=\sigdguz(a), \quad a \in \sRdz, g\in \Hdm.
\endgathered
$$
\item"Step IV:" Let $\sMdzone$ be the AFD factor of type \twoinf\ equipped with trace scaling one parameter automorphism group $\lrbrace{\thds: s\in\R}$ and set 
$$
\sRdzone=\sMdzone\botimes \sMdz.
$$
We then set the action $\ta^{\la, \mu}$ to be the following:
$$
\tilde\a_g^{\la, \mu}=\th_{\txm(g)}\otimes \sig_g^{\la,\mu}\ \text{on}\ 
\sRdzone, \quad g\in \Hdm.
$$
Set
$$
\sR=\lr{\sRdzone}^{\ta_\zdz}.
$$
Since the automorphism $\ta_\zdz=\th_{T'}\ot \sig_\zdz^{\la,\mu}$ scales the trace $\tau$ by $\la=e^{-T'}$, the \vna\ $\sR$ is an AFD factor of type \threel. Finally we define the action $\a^{\la, \mu}$ by the following:
$$
\a_g^{\la,\mu}=\ta_g^{\la,\mu}\Big|_\sR, \quad g\in H,
$$
which makes sense because $\ta_\zdz$ acts trivially on $\sR$.
\endroster
\proclaim{Theorem 1.4 (Model Action)} {\rm i)} The action $\a=\a^{\la, \mu}$ constructed above has the invariants\:
$$\gathered
N=\a^{-1}\lr{\cnt(\sRdl)}, \quad \mod\!\lr{\a_g}=\txm(g), \quad g\in H,\\
\chi(\a)=\lrbracket{\la, \mu}\in \La(\Hdm, L, M, \T), \\
\nu_\a(g)=\bracketd{\frac{ T\text{\rm Log}\lr{\la(g; \zdz)}}{\twopi}}{T}\in\rtz, \quad g\in N.
\endgathered
$$

{\rm ii)} Let $\sh: G\mapsto H$ be a cross-section of the homomorphism $\pig: H\mapsto G$. Then the outer action 
$\a_\sh^{\la, \mu}$ of $G$ has the associated modular obstruction $\d([\la, \mu])=[c^{\la, \mu}, \nu^\la]\in \thmsout(G, N, \T)$.

The construction of {\rm (i)} and {\rm (ii)} exhausts all outer actions of $G$ on the {\afdf} $\sR$ of type {\threel} up to outer conjugacy.
\endproclaim
\demo{Proof} i) Let $\ta$ denote the action $\ta^{\la,\mu}$ of $\Hdm$ on $\sRdzone$. 
Since $\sR$ is the fixed point subalgebra of $\sRdzone$ under the automorphism $\ta_\zdz$, the restriction $\a=\ta|_\sR$ of $\ta$ to $\sR$ factors through the quotient group $H=\Hdm/\lr{\Z\zdz}$. Hence $\a$ is indeed an action of $H$. Since $\sRdzone$ is a factor of type \twoinf\ and
$$\aligned
\tau\scirc \tadzdz&=\tau\scirc \th_{\txm\lr{\zdz}}=e^{-\txm\lr{\zdz}}\tau=e^{-T'}\tau\\
&=\la\tau,
\endaligned
$$
the fixed point subalgebra $\sR$ is a factor of type \threel\ and the pair $\lrbrace{\sRdzone, \tadzdz}$ is the discrete core of the factor 
$\sR$. Since $\sRdzone$ is AFD, $\sR$ is \appfindiml\ by the grand theorem of Connes, \cite{Cnn5}. As $\zdz$ is a central element of
$\Hdm$, $\ta(\Hdm)$
leaves $\sR$ globally invariant and hence its restriction to $\sR$ makes sense. The inner part $\ta(N)$, which is given by the
projective representation $\lrbrace{U(g): g\in N}$, leaves $\sR$ globally invariant, i.e., $U(g), g\in N$, normalizes
$\sR$; thus we have the inclusion $U(N)\i \tsUdz(\sR)$. Hence $N=\a\inv\lr{\Cnt(\sR)}$. As in (1.1), we have
$$
\modd(\a_h)=\txm(h), \quad h\in H.
$$
If $g, \ggdone, \ggdtwo\in N$ and $h\in H$, then 
$$\gathered
\la(g; h)=U^*(g)\ta_h(U(h\inv gh));\\
U(\ggdone)U(\ggdtwo)=\mu(\ggdone;\ggdtwo)U(\ggdone\ggdtwo);\\
\nu_\a(g)=\part_{\tadzdz}\lr{U(g)}=U(g)^*\tadzdz(U(g))=\la(g; \zdz).
\endgathered
$$
Hence $\chi(\ta)=[\la, \mu]\in\La(\Hdm, L, M, \T)$ as required.
 Finally viewing $\nu_\a$ as a homomorphism of $N$ into
$\rtz$, we get $\nu_\a\in \Hom_G(N,\allowmathbreak \rtz)$ as in the assertion of the theorem.

ii) The assertion follows from the construction of $\a^{\la, \mu}$. 
\QED

\enddemo

\subhead\nofrills{Actions and Outer Actions of Two Commuting Automorphisms on an AFD factor 
$\pmb{\sR}$ of type \threel: }
\endsubhead
In this case, we have to consider the free abelian group $G=\Z^2$ of rank
two and its extension $\Gdm\cong \Z^3$ relative to a homomorphism
$\txm: G\mapsto \rt'z$. We fix a subgroup $N$ of $G$, which is going to
represent the inverse image $\a\inv(\cnt(\sMdl))$ of the extended
modular automorphism group. We assume that $N$ is in the diagonal
form, i.e., with a free basis $\lrbrace{\adone, \adtwo}$ of $G$ the
subgroup $N$ is of the form:
$$
N=\pdone\Z\adone+\pdtwo\Z\adtwo.
$$
Of course, one can choose $\pdone$ and $\pdtwo$ in such a way that
$0\leq\pdone\leq \pdtwo$ and $\pdone$ divides $\pdtwo$, but to go
beyond the finite rank case, we don't assume that $\pdone$ is a divisor of
$\pdtwo$, which might make a matter slightly more involved. In the case
that $G=\Z^2$, we have $\ththr(G, \T)=\lrbrace{1}$, so every  outer
action of $G$ comes from an action of $G$. Since $\thtw(G, \T)\cong
\T\neq \lrbrace{1}$, the outer conjugacy class of an action is bigger than
the cocycle conjugacy class. To go further, we recall the reduced
modified \hjr, \cite{KtT3, Theorem 3.11 page 116}:
$$\eightpoint\CD
\thtw(G, \T)@>\Res_\Qdm>>\La(\Gdm, N, \T)@>\d_\Qdm>>\thasout(G,
N,
\T) @>\Inf_\Qdm>>\ththr(G, \T)=\lrbrace{1},
\endCD
$$
where $\Qdm=\Gdm/N$. Here since $\ththr(G, \T)=\lrbrace{1}$, 
we don't have to consider the resolution group $H$ and its subgroup
$M$. To identify the subgroup $N\i G$ as a subgroup of $\Gdm$, we
need a little care. First, set
$$\gathered
\zdz=\lr{0, T'}\in \Gdm,\\ \zdone=\lr{\adone, \frac
{\qdone T'}\pdone}\in\Gdm,
\quad \zdtwo=\lr{\adtwo, \frac{\qdtwo T'}\pdtwo},\\
\bdone=\lr{\pdone\adone, 0}=\pdone\zdone-\qdone\zdz\in\Gdm, \\
\bdtwo=\lr{\pdtwo\adtwo, 0}=\pdtwo\zdtwo-\qdtwo\zdz\in\Gdm,\\
N=\Z\bdone+\Z\bdtwo\i \Gdm=\Z\zdz+\Z\zdone+\Z\zdtwo, \\
\Qdm=\Gdm/N.
\endgathered\tag1.8
$$

This gives the following coordinate system in $\Gdm$ and $N$:
$$\gathered
g=\edonen(g)\bdone+\edtwon(g)\bdtwo\in N, \quad \text{i.e.,}\quad
\edin(g)=\frac{\edi(g)}\pdi,\quad i=1, 2, \\
h=\tedz(h)\zdz+\tedone(h)\zdone+\tedtwo(h)\zdtwo\in \Gdm.
\endgathered\tag1.9
$$

\proclaim{Theorem 1.6 (Invariant)}  Define $\tZ$ and 
$\tB$ by the following\:
$$\left.\aligned
\tZ&=\lrbrace{\gathered b=\lrbrace{\bbij: i=1,2, j=0,1,2}\in \R^6:\\ \pdj\bbij-\qdj\bbiz\in\Z, i=1, 2, \ j=1, 2
\endgathered}, \\
\tB&=\lrbrace{b\in \tZ: \gathered\bbiz, \bbii\in\Z, i=1,2,  \\
\pdtwo\bbonetwo+\pdone\bbtwoone\in\gcdlr{\pdone, \pdtwo}\Z
\endgathered},
\endaligned\right\}\tag1.10
$$
and to each $b\in\tZ$ associate a cochain $\lr{\ladb, \mudb}\in \tZ(\Gdm, N, \T)$ by\:
$$\gathered
\ladb(g; h)=\explrtwopii{\sumd{i=1,2; j=0,1,2}\bbij\edin(g)\tedj(h)},\\
\mudb(\ggdone; \ggdtwo)=1, \quad g, \ggdone, \ggdtwo\in N,\
h\in \Gdm.
\endgathered\tag1.11
$$
Then the cochain
$\lr{\ladb, \mudb}$ is a characteristic cocycle $(\ladb, 1)\in\tZ(\Gdm, N, \T)$. The modular obstruction cocycle $(\cdb, \nudb)=\d(\ladb, 1)\in \tzmsout(G, N, \T)$ corresponding to
$\lr{\ladb, 1}$ takes the form\:
$$\eightpoint\aligned
\cdb&\lr{\dotggdone; \dotggdtwo; \dotggdthree}=\ladb\lr{\fnn(\dotggdtwo;
\dotggdthree); \fs\lr{\dotggdthree}}, \quad \dotggdone, \dotggdtwo,
\dotggdthree\in\Qdm,\\
&=\explrtwopii{\sumd{i=1,2; j=0,1,2}\frac
{\bbij\edietan{\dotggdtwo}{\dotggdthree}\lrbrace{\tedj(\dotggdone)}_\pdj}\pdi},\\
&\hskip.2in
\nudb(g)=\lrbracketd{T\sumd{i=1,2}\bbiz\edin\lr{g}}{T}\in \rtz, \quad g\in N,
\endaligned\tag1.12
$$
where for the notations $\eta_\pdi$ and $\fnn$ we refer the definition in \S$3$, in particular {\rm (3.8)} and {\rm(3.14)}, and
furthermore
$$
\lrbrace{\tedz(\dotggdone)}_\pdz=\tedz(\dotggdone)\in\Z, \quad \dotggdone \in \Qdm.
$$
The  $(i,j)$-component $\tZ(i,j)$ and $\tB(i,j)$ of $\tZ$ and $\tB$ give more precise informations about the cocycles\:
\roster
\item"i)" For $i=1,2$, we have
$$\gathered
\tzzdbii=\lrbrace{z=(x, u)\in \R^2: \pdi x -\qdi u\in\Z},\\
\tbbdbii=\Z\oplus \Z.
\endgathered\tag1.13
$$
The bicharacter $\la_z^{i,i}$ on $N\times \Gdm$ determined by\:
$$\aligned
\la_z^{i,i}(g; h)=\explrtwopii{x\edin(g)\tedi(h)+u\edin(g)\tedz(h)}, 
\endaligned\tag1.14
$$
for each pair $g\in N, h\in \Gdm,$
gives a characteristic cocycle of $\tZ(\Gdm,\allowmathbreak N, \T)$. It is a coboundary if and only
if $z$ is in $\tbbdbii$. The corresponding cohomology class $[\la_z^{i,i}]\in\Ladbii$ has the parameterization\:
$$\aligned
\lrbracket{\la_z^{i,i}}\in \La(i,i)&\sim\lr{\lrbracket{\pdi x -\qdi u}_{\gcd(\pdi, \qdi)}, \lrbracket
{-\vdi x +\udi u}_\Z }\\
&\hskip.2in
\in \Z_{\gcd(\pdi, \qdi)}\oplus
(\R/\Z),
\endaligned\tag1.15
$$
where the integers $\udi, \vdi$ are determined by\:
$$
\pdi\udi-\qdi \vdi = \gcdlr{\pdi, \qdi}, \quad i=1,2,
$$
through the Euclid algorithm. The associated modular obstruction cohomology class $\lr{\lrbracket{c^{i,i}_z, \nu_z^{i,i}}}\in \thmsout(i,i)$ corresponds to the class\:
$$\eightpoint\left.\gathered
\lr{\lrbracket{\pdi x -\qdi u}_{\gcd(\pdi, \qdi)}, \lrbracket
{-\vdi x +\udi u}_\Z }
\in \Z_{\gcd(\pdi, \qdi)}\oplus
(\R/\Z),\\
\nu_z^{i,i}(g)=\lrbracket{Tu\edin(g)}_{T}\in \rtz,
\endgathered\right\} \quad i=1, 2.
$$

{\rm ii)}  With $ (i,j)=(1,2)$,
$$\eightpoint\gathered
\tzzdbij=\lrbrace{(x, u, y, v)\in\R^4: 
\pdj x-\qdj u\in\Z, \pdi y-\qdi v\in \Z};\\
\tbbdbij=
\lrbrace{(x, u, y, v)\in\tzdbij: \gathered
\pdj x+\pdi y\in\gcd(\pdi, \pdj)\Z,\\ u, v\in \Z
\endgathered}.
\endgathered\tag1.16
$$
For each element $z=(x, u, y, v)\in \tzzdbij$, the corresponding bicharacter $\la_z$ on $N\times
\Gdm$\:
$${\eightpoint\aligned
\la_z^{i,j}(g;
h)&=\explrtwopii{x{\edin(g)\tedj(h)}+y{\edjn(g)\tedi(h)}}\\
&\hskip.2in\times
\explrtwopii{u{\edin(g)\tedz(h)}+v{\edjn(g)\tedz(h)}},
\endaligned} \tag1.17
$$
for each pair $g\in N, h\in\Gdm$
is a characteristic cocycle in $\tZ(\Hdm,\allowmathbreak L, M, \T)$. It is a coboundary if and only if
$z\in\tbbdbij$. The cohomology class $[\la_z^{i,j}]\in\Ladbij$ of $\la_z$ corresponds to the parameter class\:
$$\left.\eightpoint\aligned
\lrbracket{z}&=\pmatrix\lrbracket{\mdij\lr{x\rdji+y\rdij}-
\ndij\lr{u\sdji+v\sdij}}_\Z\\
 \lrbracket{\ydij\lr{x\rdji +y\rdij }+\xdij\lr{u\sdji+v\sdij}}_{\Z}\\
\lrbracket{-u\wdij+v\wdji}_\Z\endpmatrix\\
&\hskip.5in
\in \pmatrix \left.\lr{\dfrac1{D(i,j)}\Z}\right/\Z\\\R/\Z\\\R/\Z\endpmatrix,
\endaligned\right\}\tag1.18
$$
where the  integers $D(i,j), \cdots, \wdij$ are those such that\:
$$\left.\gathered
D(i,j)=\gcdlr{\pdi, \pdj, \qdi, \qdj},\\  \Ddij
=\gcdlr{\pdi, \pdj}, \quad \Edij=\gcdlr{\qdi, \qdj},\\
\rdij=\frac\pdi\Ddij, \quad \rdji=\frac\pdj\Ddij\quad
\sdij=\frac \qdi\Edij, \quad \sdji=\frac\qdj\Edij,\\
\mdij=\frac \Ddij {D(i,j)}, \quad\ndij=\frac \Edij {D(i,j)},\\
\qdi\wdij+\qdj\wdji=\Edij, \quad 
\xdij \Ddij+\ydij \Edij=D(i,j).
\endgathered\right\}\tag1.19
$$
The associated modular obstruction class $\lr{\lrbracket{c_z^{i,j}}, \nu_z^{i,j}}\in \thmsout(i,j)$ corresponds to the pair of the classes\:
$$\gathered
\lrbracket{z}\in \pmatrix \left.\lr{\dfrac1{D(i,j)}\Z}\right/\Z\\\R/\Z\\\R/\Z\endpmatrix,\\
\nu_z^{i,j}(g)=\lrbracketd{T\lr{u\frac{\edone(g)}\pdone+v\frac{\edtwo(g)}\pdtwo}}{T}
\in \rtz, \quad g\in N.
\endgathered
$$ 
\endroster
\endproclaim
The proof of this special case is not particularly simpler than the general case, so that we will discuss later 
in the general free abelian group case, see Theorem 4.2.

\head{\bf \S2 Asymmetrization.}
\endhead

Set the notations:
$$
X=\Z_{n+1}=\Z/(n+1)\Z,\quad X_1=X\backslash \{1\}.
$$
The signature of a permutation $\sig$ is the sign of the product:
$$
\sign(\sig)=\sign\left\{\prod_{i<j}(\sig(j)-\sig(i))\right\}.
$$
Let $S$ be the cyclic permutation:
$$
S=(2, 3, \cdots, n, n+1, 1)\in \Pi(X),\tag2.1
$$
whose signature is given by:
$$
\sign(S)=\moneun.\tag2.2
$$
Each element $\sig\in \Pi(X_1)$ is identified with the element of $\Pi(X)$ such that
$$
\sig=(1, \sig(2), \sig(3), \cdots, \sig(n), \sig(n+1))\in \Pi(X).
$$
This identification of an element of $\Pi(X_1)$ with the
corresponding element of $\Pi(X)$ preserves the signature of
$\sig$. 
Then the total permutation group $\Pi(X)$ is the disjoint
union of the translations $\{S^k \Pi(X_1): 0\leq k \leq n\}$, i.e.,
$$
\Pi(X)=\bigcup_{k=0}^nS^k\Pi(X_1), \quad \text{disjoint union}.\tag2.3
$$
{\smc Definition 2.1.} 
The {\it asymmetrization} $\AS\xi$ of $\xi\in \tcn(G, A)$ is defined by the following:
$$
(\AS\xi)(\ggdone, \ggdtwo, \cdots, \ggdn)
=\sum_{\sig\in\Pi(\Z_n)}\sign(\sig)\xi(\ggdsigone, \ggdsigtwo,
\cdots, \ggdsign).\tag2.4
$$
Define $\pidk:
G^{n+1}\mapsto G^{n}$ by the following:
$$\allowdisplaybreaks\aligned
\pidk(&\ggdone, \ggdtwo, \cdots, \ggdn, \ggdnone)\\ 
&=\cases
(\ggdtwo, \ggdthree, \cdots, \ggdn, \ggdnone),\quad k=0;\\
(\ggdone,
\cdots, \ggdkmone, \ggdk\ggdkone, \ggdktwo, \cdots, \ggdnone),\quad 1\leq k \leq
n;\\ (\ggdone, \ggdtwo, \cdots, \ggdn),\quad k=n+1.
\endcases
\endaligned\tag2.5
$$
The boundary operation $d\in \Hom(\Z(G^{n+1}), \Z(G^n))$ is then given by
$$
d= \sumkztonone \moneuk \scirc\pi_k,\tag2.6
$$
and
$$
\part \xi=d^*\xi, \quad \xi\in \tcnone(G, \T).
$$
We view the asymmetrization $\AS$ also as an element of $\End\lr{\Z\lr{G^n}}$ determined by: 
$$\aligned
\AS \lr{\gdotsgn}=\sumd{\sig\in \Pi\lr{\Z_n}}
\sign(\sig)\lr{\gsigdotsn}.
\endaligned
$$

\proclaim{Lemma 2.2} The asymmetrization and the boundary operation are related in the
following way\:
$$
\AS\scirc d=0\quad \text{in }\ \Hom(\Z( G^{n+1}),\Z( G^n)).
$$
\endproclaim
\demo{Proof} Define $Q\in \Hom(\Z(G^{n+1}), \Z(G^n))$ and $R\in  \Hom(\Z(G^{n+1}),
\Z(G^n))$ by:
$$\eightpoint\aligned
Q&=
\sum_{\sig\in\Pi(X_1)}\sum_{j=1}^{n+1}
\Big(\sign\left(S^{j-1}\sig\right)\pidz S^{j-1}\sig
+\moneunone\sign\left(S^{j}\sig\right)\pidnone S^j\sig\Big),
\endaligned
$$
and
$$
Rg=\sum_{\sig\in\Pi(X_1)}\sum_{j=1}^{n+1}\sign(S^j\sig)\sumkton
\moneuk\pidk S^j\sig g, \quad g\in G^{n+1}.
$$
So we have
$$
\AS\scirc d=Q+R.
$$
We know 
$$\allowdisplaybreaks\aligned
&\hskip.5in\pidz S^{j-1}\sig=\pidnone S^j\sig,\quad 1\leq j\leq n;\\
\sign&(S^{j-1}\sig)\pidz S^{j-1}\sig+\moneunone\sign(S^j\sig)\pidnone
S^j\sig\\
&=(-1)^{n(j-1)}\sign(\sig)\pidz S^{j-1}\sig+\moneunone (-1)^{nj}\sign(\sig)\pidnone
S^j\sig\\ &=0.\\
\endaligned
$$
Thus we get
$$
Q=0.
$$

We need the notation $\sigdkkone$ for the flip of $k$ and $k+1$:
$$
\sigdkkone=(1, 2, \cdots, k-2, k-1, k+1, k, k+2, k+3,\cdots, n+1)\in\Pi(X).
$$
Then we get
$$\aligned
\sign(\sigdkkone \rho)\pidk \sigdkkone \rho g+\sign(\rho)\pidk \rho g=0, \quad \rho\in\Pi(X), 1\leq
k\leq n.
\endaligned
$$
Hence we come to the following:
$$\aligned
R&=\sum_{\sig\in\Pi(X_1)}\sum_{j=1}^{n+1}\sign(S^j\sig)\sumkton
\moneuk\pidk S^j\sig\\
&=\sumkton\moneuk\sum_{j=1}^{n+1}\sum_{\sig\in\Pi(X_1)}\sign(S^j\sig)
\pidk S^j\sig\\
&=\sumkton\moneuk \sum_{\rho\in\Pi(X)}\sign(\rho)\pidk\rho\\
&=\sumkton\moneuk\sum_{\rho\in\Pi_0(X)}(\sign(\rho)\pidk\rho+\sign(\sigdkkone
\rho)\pidk \sigdkkone\rho)\\
&=0,
\endaligned
$$
where $\Pi_0(X)$ is the group of even permutations of $X$, i.e., the alternating group.
Therefore we conclude
$$\aligned
\AS \scirc d=Q+R=0.
\endaligned
$$
This completes the proof.
\QED
\enddemo

Let $\Cal A$ be a $G$-module with action $\alpha$.
We recall the dimension shifting theorem and the dimension shift  map $\part$. First we define a new $G$-module $\widetilde{\Cal A}$ as follows:
\roster
\item"i)" Let $\Map(G, \sA)$ be the module $\sA^G$ of all
$\sA$-valued functions on $G$ with pointwise addition.
\item"ii)" View the group $\sA$ as the submodule of $\Map(G, \sA)$ of constant $\sA$-valued functions.
\item"iii)" The action $\a$ of $G$ on $\sA$ is extended to the enlarged additive group
$\Map(G, \sA)$ by:
$$
\lr{\a_h f}(g)=\a_{h}\lr{f(gh)}, \quad f\in \Map(G, \sA),\quad
g, h\in G.
$$
\item"iv)" Form the quotient $G$-module:
$$
\widetilde \sA=\Mapgsa/\sA.
$$
\endroster
Thus we obtain the following equivariant short exact sequence:
$$\CD
0@>>>\sA@>>>\Mapgsa@>>>\widetilde \sA@>>>0.
\endCD\tag2.7
$$
The short exact sequence (2.7) splits in the following way: 
\roster
\item"i)" First, set
$$
j(f)(g)=f(g)-f(e), \quad f\in\Mapgsa, g\in G,
$$
where $e\in G$ is the neutral element of $G$. Then the map $j$ is a homomorphism of $\Mapgsa$ onto the subgroup 
$\Mapdzgsa$ of all $\sA$-valued functions on $G$ vanishing at $e$. Then we get 
$$
\Ker(j)=\sA \i \Mapgsa,
$$ 
so that the map $j$ is viewed as a bijection from $\tsA$ onto $\Mapdzgsa$.

\item"ii)" The map $j$ transforms
the action $\ta$ of $G$ on $\tsA$ to the action, denoted by $\ta$ again, 
on
$\Mapdzgsa$ defined by: 
$$
\lr{\ta_h f}(g)=\a_h\lr{f(gh)}-\a_h\lr{f(h)}, \quad g, h\in G, \quad f\in \Mapdzgsa.
$$
\endroster
With the map $j$, we will identify $\tsA$ and $\Mapdzgsa$. Thus we have a
short exact sequence:
$$\CD
0@>>>\sA@>i>>\Mapgsa@>j>\underset{\fs}\to\longleftarrow>\tsA=\Mapdzgsa@>>>0
\endCD
$$
Let $\fs$ denote the embedding of $\tsA=\Mapdzgsa \hookrightarrow \Mapgsa$,
which is a right inverse of the map $j$. If $\tu \in\tZaunmone\lr{G, \tsA}$, then
$$
0=\partg \tu =j\ \lr{\tpartg \fs(\tu )},
$$
where $\tpartg$ means the coboundary operator in $\tCtaun\lr{G, \Mapgsa}$, 
so that we have $\partg \fs(\tu )\in \tZaun\lr{G, \sA}$. We denote the cohomology
class $\lrbracket{\tpartg \fs(\tu )}\in \tHaun\lr{G, \sA}$ by $\part [\tu ]$ for each
$[\tu ]\in\tHtaunmone\lr{G, \tsA}$.  It is known as the dimension shift theorem that the map $\part$ is
an isomorphism of $\tHtaunmone\lr{G, \tsA}$ onto $\tHaun\lr{G, \sA}$.

{\smc Definition 2.3.} Suppose that 
the group $G$ admits a torsion free central element $\zdz\in G$. A cocycle
$c\in \tZaun\lr{G, \sA}$ is said to be of the {\it standard form} (relative to the
central element $\zdz$) if \roster
\item"i)" For each $\kdone, \cdots \kdn\in\Z$ and $\gdotsgn\in G$,
$$
c(\zdzukone\ggdone, \cdots, \zdzukn\ggdn)=\a_\ggdone\Big(d_c\lr{\kdone;
\ggdtwo,\cdots, \ggdn}\Big)+c\lr{\gdotsgn};
\tag2.8
$$
\item"ii)" The map $k\in\Z\mapsto d_c\lr{k; \gdotsgtwon}\in\sA$ is in
$\tZ_{\a_\zdz}^1\lr{\Z, \sA}$ for each $\gdotsgtwon\in G$, i.e.,
$$\aligned
d_c&\lr{k+\ell; \gdotsgtwon}\\
&\qquad=d_c(k; \gdotsgtwon)+\a_\zdz^k\Big(
d_c\lr{\ell; \gdotsgtwon}\Big).
\endaligned\tag2.9
$$
\item"iii)" For each $k\in\Z$ and $\gdotsgn\in G$, we have
$$\aligned
\lr{\partg d_c}&\lr{k; \gdotsgn}\\
&\qquad=
\a_\zdz^k\Big(c\lr{\gdotsgn}\Big)-c\lr{\gdotsgn}.
\endaligned\tag2.10
$$
\endroster

{\smc Remark 2.4.} The cocycle identity (2.8) can be fulfilled automatically if $d_c$ is chosen in such a way that
$$\gathered
c\lr{\zdz\ggdone, \zdzukdtwo\ggdtwo, \cdots, \zdzukdn\ggdn}=\a_\ggdone\Big(d_c\lr{\gdotsgtwon}\Big)
+c\lr{\gdotsgn}, \\
\lr{\partg d_c}\lr{\gdotsgn}=\a_\zdz\Big(c\lr{\gdotsgn}\Big)
-c\lr{\gdotsgn}.
\endgathered
$$
Because $d_c\lr{k; \gdotsgtwon}$ can be obtained inductively by:
$$\aligned
d_c&\lr{k; \gdotsgtwon}\\
&\qquad=d_c\lr{\gdotsgtwon}
+\a_\zdz\Big(d_c\lr{k-1; \gdotsgtwon}\Big).
\endaligned\tag2.11
$$
In the sequel, we often write $d_c\lr{\gdotsgtwon}$ for the $d$-part of a standard cocycle $c$ without referring to the first variable $k$ in $d_c\lr{k; \gdotsgtwon}$.

\proclaim{Lemma 2.5} In the above context, every cocycle $c\in \tZaun\lr{G,
\sA}$ is cohomologous to a cocycle $\cdtxs$ of the standard form.
\endproclaim
\demo{Proof} For $n=1$, the cocycle identity:
$$
c(\zdzuk g)=\a_g\lr{c\lr{\zdzuk}}+c(g), \quad k\in \Z, g\in G,
$$
shows that with $d_c(k)=c\lr{\zdzuk}$ the cochain $d_c$ satisfies the condition
(i). Now we have
$$\aligned
\a_\zdz^k\Big(c(g)\Big)-c(g)&=c\lr{\zdzuk g}-c\lr{\zdzuk}-c(g)\\
&=c(g)+\a_g\Big(c\lr{\zdzuk}\Big)-c\lr{\zdzuk}-c(g)\\
&=\a_g\Big(d_c\lr{k}\Big)-d_c\lr{k}\\
&=\Big(\partg d_c\Big)\lr{k; g},
\endaligned
$$
which shows the property (ii) for $c$ and $d_c$.

Now assume that our claim is valid for $1, \cdots, n-1$ and for any $G$-module
$\lrbrace{\sA, \a}$.

Choose an equivariant short exact sequence:
$$\CD
0@>>>\sA@>i>>M@>j>\underset{\fs}\to\longleftarrow>\tsA@>>>0
\endCD
$$
such that 
\roster
\item"i)" $\tHaun\lr{G, M}=\{0\}, n\geq 1$
\item"ii)" the cross-section $\fs:\tsA\mapsto M$ is a homomorphism of $\tsA$ into
$M$, but not equivariant, 
\endroster
so that the map
$\partg \fs: \tZaunmone\lr{G, \tsA}\mapsto \tZaun\lr{G, \sA}$ gives rise to
an isomorphism $\part: \tHaunmone\lr{G, \tsA}\mapsto \tHaun\lr{G, \sA}$. For a standard cocycle $\tc\in \tZaunmone\lr{G, \tsA}$, we set, for each
$\zdzukdone\ggdone, \cdots, \zdzukdnmone\ggdnmone\in G$,
$$\aligned
\bc\lr{\zdzukdone\ggdone, \cdots, \zdzukdnmone\ggdnmone}&=\a_\ggdone\Big(\fs\lr{d_\tc\lr{\kdone; \gdotsgtwonmone}}\Big)\\
&\hskip.5in+
\fs\lr{\tc\lr{\gdotsgnmone}}.
\endaligned
$$
Since $j(\bc)=\tc$, we have 
$$
c=\partg \bc\in \tZaun\lr{G, \sA}.
$$
We then compute
$$\allowdisplaybreaks\eightpoint\align
c&\lr{\zdzukdone\ggdone, \cdots, \zdzukdn\ggdn}=(\partg
\Bar c)\lr{\zdzukdone\ggdone, \cdots, \zdzukdn\ggdn}\\
&=\a_{\zdzukdone\ggdone}\Big(\bc\lr{\zdzukdtwo\ggdtwo,
\zdzukdthr\ggdthree, \cdots, \zdzukn\ggdn}\Big)\\
&\hskip.5in+
\sumdjonetonmone\moneuj \bc\lr{\zdzukdone\ggdone, \cdots, \zdzukdj
\ggdj\zdzukdjone\ggdjone, \cdots, \ggdn}\\
&\hskip1in+
\moneun\bc\lr{\zdzukdone\ggdone, \cdots, \zdzukdnmone\ggdnmone}\\
&=\a_{\zdzukdone\ggdone}\Big[\a_\ggdtwo\Big(\fs\lr{d_\tc\lr{\kdtwo; \ggdthree,
\cdots, \ggdn}}\Big)+\bc\lr{\gdotsgtwon}\Big]\\
&\hskip.2in
-\Big[{\a_{\ggdone\ggdtwo}\Big(\fs \lr{d_\tc\lr{\kdone+\kdtwo; \gdotsgthrn}}\Big)+
\bc\lr{\ggdone\ggdtwo, \gdotsgthrn}}\Big]\\
&\hskip.2in+
\sumdjtwotonmone\moneuj \Big[\a_\ggdone\Big(\fs\lr{d_\tc\lr{\kdone; \ggdtwo,
\cdots, \ggdj\ggdjone, \cdots, \ggdn}}\Big) \\
&\hskip1in+
\bc\lr{\ggdone, \cdots, 
\ggdj\ggdjone, \cdots, \ggdn}\Big)\Big]\\
&\hskip.2in+
\moneun\a_\ggdone\Big(\fs\lr{d_\tc\lr{\kdone;
\gdotsgtwonmone}}\Big)\\
&\hskip1in+
\moneun\bc\lr{\gdotsgnmone}\\
&=(\partg\bc)\lr{\gdotsgn}+
\a_{\zdzukdone\ggdone}\Big(\a_\ggdtwo\lr{\fs\lr{d_\tc\lr{\kdtwo;
\ggdthree,
\cdots,
\ggdn}}}\Big)\\
&\hskip.2in+
\a_\ggdone\lr{\a_\zdz^\kdone\Big(\bc\lr{\gdotsgtwon}}-
\bc\lr{\gdotsgtwon}\Big)\\
&\hskip.2in -\a_{\ggdone\ggdtwo} \Big[\fs\lr{d_\tc\lr{\kdone; \gdotsgthrn}}
+\a_\zdz^\kdone\lr{\fs \lr{d_\tc\lr{\kdtwo; \gdotsgthrn}}}\Big]\\
&\hskip.2in+
\sumdjtwotonmone\moneuj \a_\ggdone\lr{\fs\lr{d_\tc\lr{\kdone; \ggdtwo,
\cdots, \ggdj\ggdjone, \cdots, \ggdn}}} \\
&\hskip.2in+
\moneun\Big(\a_\ggdone\lr{\fs\lr{d_\tc\lr{\kdone;
\gdotsgtwonmone}}}\Big)\\
&=(\partg\bc)\lr{\gdotsgn}\\
&\hskip.5in+
\a_\ggdone\Big[\lr{\a_\zdz^\kdone\Big(\bc\lr{\gdotsgtwon}}-
\bc\lr{\gdotsgtwon}\Big)\\
&\hskip1in -\a_{\ggdtwo} \Big(\fs\lr{d_\tc\lr{\kdone; \gdotsgthrn}}\Big)\\
&\hskip.5in+
\sumdjtwotonmone\moneuj \fs\lr{d_\tc\lr{\kdone; \ggdtwo,
\cdots, \ggdj\ggdjone, \cdots, \ggdn}} \\
&\hskip1in+
\moneun\lr{\fs\lr{d_\tc\lr{\kdone;
\gdotsgtwonmone}}}\Big]\\
&=(\partg\bc)\lr{\gdotsgn}\\
&\hskip.5in+
\a_\ggdone\Big[{\a_\zdz^\kdone\Big(\bc\lr{\gdotsgtwon}}-
\bc\lr{\gdotsgtwon}\Big)\\
&\hskip1in -
\partg \lr{\fs\scirc d_\tc}\lr{\kdone; \gdotsgtwon}\Big].
\endalign
$$
Consequently, we get 
$$\aligned
c\lr{\zdzukdone\ggdone, \cdots, \zdzukdn\ggdn}&=
\a_\ggdone\lr{d_c\lr{\kdone; \gdotsgtwon}}+c\lr{\gdotsgn}
\endaligned
$$
with
$$\aligned
c\lr{\gdotsgn}&=(\partg\bc)\lr{\gdotsgn},\\
d_c\lr{\gdotsgtwon}&=\a_\zdz\Big(\bc\lr{\gdotsgtwon}\Big)-
\bc\lr{\gdotsgtwon}\\
&\hskip.4in -\partg \lr{\fs\scirc d_\tc}\lr{\gdotsgtwon}.
\endaligned
$$
We now check the requirement (2.10) for $d_c$ and $c$:
$$\aligned
\a_\zdz&\Big(c\lr{\gdotsgn}\Big)-c\lr{\gdotsgn}\\
&=\a_\zdz\Big(\partg \bc\lr{\gdotsgn}\Big)-\partg \bc\lr{\gdotsgn}\\
&=\partg\Big(\a_\zdz\lr{\bc\lr{\gdotsgn}}-\bc\lr{\gdotsgn}\Big)\\
&=\partg \Big(d_c\lr{\gdotsgtwon}+\partg\fs\scirc d_\tc\lr{\gdotsgtwon}\Big)\\
&=\partg d_c\lr{\gdotsgtwon}.
\endaligned
$$
Thus the cocycle $c$ is standard. This completes the proof.
\QED
\enddemo

We now state the main result on the asymmetrization which extends the work of Olesen-Pedersen -Takesaki, \cite{OPT}:
\proclaim{Theorem 2.6} Let $Q$ be a countable torsion free abelian group. 

{\rm i)} The asymmetrization $\AS$ maps the group $\tZ^n\lr{Q, \T}$ of $\T$-valued $n$-th cocycles  onto 
the compact group $X^n(Q, \T)$  of all asymmetric multi-characters on $n$ variables of $Q$.

{\rm ii)} The following sequence is exact for each $n\in\N$\:
$$\CD
1@>>>\tB^n(Q, \T)@>>>\tZ^n(Q, \T)@>\AS>>X^n(Q, \T)@>>>1.
\endCD
$$
Consequently, $$\tH^n(\Z^m, \T)\cong X^n(\Z^m,
\T)\cong\cases
\T^{\frac{m!}{n!(m-n)!}}&m\geq n\\
0&   m<n\ .
\endcases
$$ More generally, if $Q$ is a countable torsion free abelian
group, then the cohomology group
$\tH^n(Q, \T)$ is naturally isomorphic to the Pontrjagin - Kampen dual of the
$n$-th exterior power
$Q\wedge Q\wedge \cdots \wedge Q$ of $Q$\:
$$
\tH^n(Q, \T)\cong \text{The Pontrjagin - Kampen Dual of }
(Q\wedge Q\wedge\cdots \wedge Q).
$$
{\rm iii)} The group $X^n(Q, \T)$ is a subgroup of $\tZ^n(Q, \T)$ such that
$$\aligned
\tZ^n(Q,& \T)=X^n(Q, \T)\tB^n(Q, \T);\\
X^n(Q, \T)&\cap\tB^n(Q, \T)=\Ker(\text{\rm Power } n!),
\endaligned
$$
and 
$$\AS c=c^{n!}, \quad c\in X^n(Q, \T).
$$
\endproclaim

{\smc Remark 2.7.} If the group $Q$ has torsion, then the theorem fails as
seen in the case that $Q=\Z_p=\Z/p\Z, p\geq 2$,  $\ththr(Q, \T)\cong \Z_p$ and $X^3(Q, \T)=\lrbrace{0}$.

For the proof, we need some preparation. First, if $n=1$, then our assertion is trivially true for any abelian group $Q$ with no assumption on torsion. We then assume that our assertion is true for cocycle dimension  $1, \cdots, n-1$ with $n\in\N$ fixed and for any torsion free abelian group $Q$.  With this induction hypothesis, we prepare a couple of lemmas for cocycle dimension $n$.

\proclaim{Lemma 2.8} {\rm i)} If $M$ is an abelian group such that a cocycle $c\in \tZunlr{M, \T}$ is a coboundary if and only if $\AS c=1$, then the same is true for the product group $Q=M\times \Z$.

{\rm ii)} If $M$ is an abelian group such that the asymmetrization $\AS c$ of each cocycle $c\in\tZunlr{M, \T}$ is a multi-character, then the same is true for the product group $Q=M\times \Z$.
\endproclaim
\demo{Proof} Let $\zdz$ denote the distinguished element of $Q$ associated with the product decomposition $Q=M\times \Z$ so that every element $q\in Q$ is written uniquely in the form
$q=m\zdzuk, m\in M, k\in\Z$.

 i)  The triviality of the asymmetrization of a coboundary was proven in Lemma 2.2. Thus we prove the converse.  
Suppose that 
$\AS c=1, c\in \tZ^n(Q, \T)$. By Lemma 2.5 the cocycle $c$ is cohomologous to a cocycle  $c_\txs$ 
of standard form and $\AS c_\txs=\AS c=1$
by Lemma 2.2.
So we may and do assume that $c$ is standard:
$$
c(\tp_1, \tp_2, \cdots,\tp_n)=d_c(p_2,p_3,\cdots,p_n )^{\ell_1}c_M(p_1, p_2,\cdots,p_n)
$$ 
where $
\tp_i=p_iz_0^{\ell_i}\in Q=M\times \Z$. As $Q$ does not act on $\T$, the $d$-part $d_c$ is a cocycle in $\tZunmonelr{Q, \T}$.

We look at  the  asymmetrization of $c$:
$$\aligned
(\AS c)\lr{\tp_1,\tp_2,\cdots,\tp_n}
&=
\prod_{\sigma\in S_n}\Big(d_c\lr{p_{\sigma(2)},
p_{\sigma(3)},\cdots,p_{\sigma(n)}}^{\ell_{\sigma(1)}}\\
&\hskip.5in\times
c_M\lr{p_{\sigma(1)},p_{\sigma(2)},\cdots,p_{\sigma(n)}}
\Big)^{\sign\ \sigma}\\
&=
\prod_{\sigma\in S_n}d_c\lr{p_{\sigma(2)},p_{\sigma(3)},\cdots,p_{\sigma(n)}}^{\ell_{\sigma(1)}\sign\ \sigma}\\
&\hskip.5in\times
\prod_{\sigma\in S_n}c_M\lr{p_{\sigma(1)},p_{\sigma(2)},\cdots,p_{\sigma(n)}}^{\sign\ \sigma},
\endaligned
$$
i.e.,
$$\aligned
(\AS c)\lr{\tp_1,\tp_2,\cdots,\tp_n}&=\prod_{\sigma\in S_n}d_c\lr{p_{\sigma(2)},p_{\sigma(3)},\cdots,p_{\sigma(n)}}^{\ell_{\sigma(1)}\sign\ \sigma}\\&\hskip.5in\times
\lr{\AS c_M}\lr{p_1,p_2,\cdots,p_n}.\\
\endaligned\tag2.12
$$
To compute the first term of the above expression, we take a closer look at the permutation group $S_n$. In particular, we have to pay attention to the fact that the first term in the variables of $d_c$ is mission. To this end, we fix $k, 1\leq k\leq n$, which represents the missing term in $d_c$, and
 consider the cyclic permutation: 
$$\eightpoint
 S_{n-1}(k)=\lr{1,2,\cdots, k-1,k+1,\cdots, n}\in \Pi\Big(\lrbrace{1,2,\cdots, k-1,k+1,\cdots, n}\Big).
$$
For $
\sigma=(k,\sigma(2),\sigma(3),\cdots,\sigma(n)) \in S_n
$, define $\rho,\ \tilde\rho$ and $\tilde\sigma$ as follows:
$$\eightpoint\aligned
\rho&=S^{(n-k+1)}\sigma\\
&=
\pmatrix
1&2&\cdots&k-1&k&k+1&\cdots&n\\
\sigma({n-k+2})&\sigma({n-k+3})&\cdots&\sigma(n)&k&\sigma(2)&\cdots&\sigma({n-k+1})
\endpmatrix;\\
\tilde \rho&=\pmatrix
1&2&\cdots&k-1&k+1&\cdots&n\\
\sigma({n-k+2})&\sigma({n-k+3})&\cdots&\sigma(n)&\sigma(2)&\cdots&\sigma({n-k+1})
\endpmatrix;\\
\tilde\sigma&={S_{n-1}(k)}^{k-1}\tilde\rho\\
&=
\pmatrix
1&2&\cdots&k-1&k+1&\cdots&n\\
\sigma(2)&\sigma(3)&\cdots&\sigma({k})&\sigma({k+1})&\cdots&\sigma(n)&
\endpmatrix\\
&=\big(\sigma(2),\sigma(3),\cdots,\sigma(n)\big).
\endaligned
$$
Then observing $\sign\ \tilde\rho=\sign\ \rho$, we compute
$$\aligned
&\sign\  \sigma
=\sign\ S^{k-1}\sign\ \rho=(-1)^{(n-1)(k-1)} 
\sign\ \tilde\rho\\
&= 
(-1)^{(n-1)(k-1)}\sign\lr {S_{n-1}(k)^{n-k}}\sign\ \tilde\sigma\\
&=(-1)^{(n-1)(k-1)+(n-2)(n-k)}\sign\ \tilde\sigma
=(-1)^{k-1}\sign\ \tilde\sigma.
\endaligned
$$
Hence  the first term of (2.12) becomes the following:
$$\allowdisplaybreaks\eightpoint\align
\prod_{\sigma\in S_n}&\left\{d_c(p_{\sigma(2)},p_{\sigma(3)},\cdots,p_{\sigma(n)})\right\}^{\ell_{\sigma(1)}\sign\ \sigma}\\
&=\prod_{k=1}^n\left\{\prod_{\tilde\sigma\in S_{n-1}(k)}
\left\{d_c(p_{\tilde\sigma(1)},p_{\tilde\sigma(2)},\cdots,p_{\tilde\sigma({n-1})})\right\}^{\sign\ \tilde\sigma}\right\}^{\ell_k(-1)^{k-1}}\\
&=\prod_{k=1}^n\left\{(\AS d_c)(p_1,p_2,\cdots, \overset{\smallsmile}\to{p_k},\cdots, p_n)\right\}^{\ell_k(-1)^{k-1}}
\endalign
$$
where the notation $\smallsmile$ stands for removing the corresponding variable. Thus 
(2.12) is replaced by the following:
$$\aligned
(\AS c)&(\tp_1,\tp_2,\cdots,\tp_n)\\
&=\prod_{k=1}^n\left\{(\AS d_c)(p_1,p_2,\cdots, \overset{\smallsmile}\to{p_k},\cdots, p_n)\right\}^{\ell_k(-1)^{k-1}}\\
&\hskip3cm\times
(\AS c_M)(p_1,p_2,\cdots,p_n).
\endaligned\tag2.12$'$
$$
The condition $\AS c=1$
yields that:
$$\gathered
\AS c_M=1\quad \text{with } 
\elldk=0, k=1,\dots, n;\\
\AS d_c=1 \quad \text{with } \elldone=1, \elldk=0, k=2,\dots, n,\ p_1=e.
\endgathered
$$
Hence $c_M$ and $d_c$ are both coboundaries by
the induction hypothesis. 
Choose $b\in \tC^{n-1}(M, \T)$ and
$a\in\tC^{n-2}(M, \T)$ such that
$$
c_M=\partm b\quad \text{and}\quad d_c=\partm a.
$$
Then the cocycle  $c$ has the form:
$$\aligned
c(\tp_1&, \tp_2,\cdots,\tp_n)=
d_c(p_2,p_3,\cdots, p_n)^{\ell_1}
c(p_1, p_2,\cdots, p_n)\\
&=\left((\partm a)(p_2,p_3,\cdots, p_n)\right)^{\ell_1}
(\partm b)(p_1, p_2,\cdots, p_n).
\endaligned
$$
Setting
$$
f(\tp_1, \tp_2,\cdots, \tp_{n-1})
=a(p_2,p_3,\cdots, p_{n-1})^{-\ell_1}
b(p_1, p_2,\cdots, p_{n-1}),
$$
for $\tp_i=\zdzu\elldi\pdi\in Q, i=1,\cdots, n-1$, 
we compute
$$\allowdisplaybreaks\align
(\partq f)&(\tp_1, \tp_2,\cdots, \tp_n)\\
&=f(\tp_2,\tp_3,\cdots, \tp_n)
\times\prod_{k=1}^{n-1}
f(\tp_1,\cdots,\tp_k\tp_{k+1},\cdots, \tp_n)^{(-1)^k} \\
&\hskip5cm\times
f(\tp_1,\tp_2,\cdots, \tp_{n-1})^{(-1)^n}
\\
&=a(p_3,\cdots, p_n)^{-\ell_2}
a(p_3,\cdots, p_n)^{\ell_1+\ell_2}\\
&\hskip3cm\times\prod_{k=2}^{n-1} 
a(p_2,\cdots,p_kp_{k+1},\cdots,p_n)^{-\ell_1(-1)^k}\\
&\hskip5cm\times 
a(p_2,p_3,\cdots, p_{n-1})^{-\ell_1(-1)^n}
\\
&\hskip1cm\times b(p_2,p_3,\cdots, p_{n})
\times\prod_{k=1}^{n-1} 
b(p_1,\cdots,p_kp_{k+1},\cdots,p_n)^{(-1)^k}
\\
&\hskip5cm\times 
b(p_1,p_3,\cdots, p_{n})^{(-1)^n}
\\
\endalign
$$
$$\eightpoint\aligned&=a(p_3,\cdots, p_n)^{\ell_1}\prod_{k=2}^{n-1} 
a(p_2,\cdots,p_kp_{k+1},\cdots,p_n)^{-\ell_1(-1)^k}\\
&\hskip2cm\times 
a(p_2,p_3,\cdots, p_{n-1})^{-\ell_1(-1)^n}
\times (\partm b)(p_1, p_2,\cdots, p_n)\\
&=\left((\partm a)(p_2,p_3,\cdots, p_{n})\right)^{\ell_1}(\partm b)(p_1, p_2,\cdots, p_n)\\
&=c(\tp_1, \tp_2,\cdots,\tp_n).
\endaligned
$$
Therefore $c$ is a coboundary. This completes the proof of the assertion (i).

ii)  Fix a standard cocycle $c\in \tZ^n(Q, \T)$:
$$
c(\tp_1, \tp_2,\cdots,\tp_n)=
d_c(p_2,p_3,\cdots, p_n)^{\ell_1}
c(p_1, p_2,\cdots, p_n)
$$
with $d_c\in \tZ^{n-1}(M, \T)$ and $c_M\in 
\tZ^n(M, \T)$. Observing that $\AS c_M$ and $\AS d_c$ are both multi-characters by the assumptions, we compute with (2.13), for $\tq_1=q_1z_0^{k_1}$,
$$\eightpoint\aligned
(\AS c)&(\tp_1\tq_1,\tp_2,\cdots,\tp_n)\\
&=
(\AS d_c)(p_2,\cdots,p_n)^{\ell_1+k_1}\\
&\hskip.2in\times
\prod_{j=2}^n\left\{(\AS d_c)(p_1q_1,p_2,\cdots, \overset{\smallsmile}\to{p_j},\cdots, p_n)\right\}^{\ell_j(-1)^{j-1}}\\
&\hskip.4in\times
(\AS c_M)(p_1q_1,p_2,\cdots,p_n)\\
&=
(\AS d_c)(p_2,\cdots,p_n)^{\ell_1}\\
&\hskip.2in\times
\prod_{j=2}^n\left\{(\AS d_c)(p_1,p_2,\cdots, \overset{\smallsmile}\to{p_j},\cdots, p_n)\right\}^{\ell_j(-1)^{j-1}}\\
&\hskip0.5cm\times(\AS d_c)(p_2,\cdots,p_n)^{k_1}\\
&\hskip.4in\times
\prod_{j=2}^n\left\{(\AS d_c)(q_1,p_2,\cdots, \overset{\smallsmile}\to{p_j},\cdots, p_n)\right\}^{\ell_j(-1)^{j-1}}\\
&\hskip2cm\times
(\AS c_M)(p_1,p_2,\cdots,p_n)
(\AS c_M)(q_1,p_2,\cdots,p_n)
\\
&=(\AS c)(\tp_1,\tp_2,\cdots,\tp_n)
(\AS c)(\tq_1,\tp_2,\cdots,\tp_n).
\endaligned
$$
Thus $\AS c$ is indeed multiplicative on the first variable, so that it is an asymmetric multi-character 
of $Q=M\times \Z$.
\QED
\enddemo

\proclaim{Lemma 2.9} Suppose  that $c\in \tZ^n(Q, \T)$ has a trivial asymmetrization, i.e., $\AS
c=1$.  Assume the following\:
\roster
\item"a)" $M$ is a finitely generated subgroup of $Q$\;
\item"b)" $a_0\in Q$ but not in $M$\;
\item"c)" $f\in\tC^{n-1}(M, \T)$ cobounds the restriction $c_M$ of $c$ to $M$, i.e.,
$$
\partm f=c_M.
$$
\endroster
Then the cochain $f$ has an extension to the subgroup $N=\langle M, \adz\rangle$ generated by
$M$ and $\adz$ such
that 
$$
\partn f=c_N,
$$
where $c_N$ is the restriction of $c$ to the subgroup $N$.
\endproclaim
\demo{Proof} To apply the structure theory of abelian groups, we use the additive group operation in the group $Q$.
From the general theory of abelian groups, it follows that $M$ and $N$ are both
free abelian groups and there exists a free basis $\{\zdone, \zdtwo, \cdots ,\zdm\}$ of $N$ and non-negative
integers $\{\pdone, \pdtwo, \cdots, \pdr\}\i \Z_+, 1\leq r\leq m$, such that
$$
N=\langle \zdone, \zdtwo, \cdots ,\zdm\rangle, \quad M=\langle \pdone\zdone, \cdots, \pdr\zdr\rangle.
$$
With the assumption for $n-1$, every $(n-1)$-cocycle $\mu\in \tZ^{n-1}(M, \T)$ is cohomologous to an
asymmetric multi-character $\mu_{\text a}$, i.e., there exist   $a_{i_1,i_2,\cdots,i_r}\in \R$ such that
$$\aligned
\mu_{\text a}&\lr{g_1, g_2,\cdots,g_{n-1}}\\
&=\exp\Big[\twopii \sum\Sb i_j\in\{1,2,\cdots, r\}\\1\leq i_1<i_2<\cdots<i_{n-1}\leq n-1\endSb
a_{i_1,i_2,\cdots,i_{n-1}} \\
&\hskip.5cm\times
\Big(e_{i_1,M}\wedge e_{i_2,M}\wedge\cdots \wedge e_{i_{n-1},M}\Big)(g_1, g_2,
\cdots,  g_{n-1})\Big]
\endaligned
$$
where $\lrbrace{\edim: 1\leq i\leq r}$ is the coordinate sysem of $M$ relative to the basis
$\lrbrace{\pdone\zdone, \cdots ,\pdr\zdr}$. 
Setting 
$$\aligned
\nu_a&(g_1, g_2,\cdots,g_{n-1})\\
&=\exp\Big\{\twopii \sum\Sb i_j\in\{1,2,\cdots, r\}\\1\leq i_1<i_2<\cdots<i_{n-1}\leq n-1\endSb
\dfrac{a_{i_1,i_2,\cdots,i_{n-1}}}{p_{i_1}p_{i_2}\cdots p_{i_{n-1}}} \\
&\hskip3cm\times
\Big(e_{i_1}\wedge e_{i_2}\wedge\cdots \wedge e_{i_{n-1}}\Big)(g_1, g_2,
\cdots,  g_{n-1})\Big\}
\endaligned
$$
where $\lrbrace{\edi: 1\leq i\leq m}$ is the coordinate system of $N$ relative to the basis
$\lrbrace{\zdone, \cdots ,\zdm}$, we obtain an extension $\nu$ of $\mu_{\text
a}$. Choose $\xi\in \tC^{n-2}(M, \T)$ so that $\mu=(\partm \xi)\mudtxa$ and extends $\xi$ to a cochain
$\xi\in \tC^{n-2}(N, \T)$. Then the second cocycle $(\partn \xi)\nu$ gives an extension of the original
$(n-1)$-cocycle $\mu\in \tZ^{n-1}(M, \T)$. Thus we obtain the surjectivity of the restriction map $\res:
\mu\in\tZ^{n-1}(N, \T)\mapsto \mum\in \tZ^{n-1}(M, \T)$, 
i.e. the exactness of the sequence:
$$\CD
\tZ^{n-1}(N, \T)@>\res>>\tZ^{n-1}(M, \T)@>>>1.
\endCD
$$

By induction on generators, Lemma 2.9 yields that the restriction $c_N$ of $c$ to $N$ is a coboundary. Hence there exists $\xi\in \tC^{n-1}(N, \T)$ such that
$c_N=\partn \xi$. Then we have $\partm f=c_M=\partm \xi_M$, so that we obtain
$\mum=\xi_M^{-1}f\in \tZ^{n-1}(M, \T)$. By the first arguments, we can extend $\mum$ to an element
$\nu\in \tZ^{n-1}(N, \T)$. Set 
$$
f=\nu\xi\in \tC^{n-1}(N, \T),
$$
and  the newly defined cochain $f$ on $N$ extends the original $f\in \tC^{n-1}(M,\allowmathbreak \T)$
and cobounds the cocycle $c_N$:
$$
\partn f=(\partn \nu)( \partn \xi)=\partn \xi=c_N.\nopagebreak
$$
This completes the proof.
\QED
\enddemo

We are now ready to complete the proof of Theorem 2.6, proceeding from cocycle dimension  $1, \cdots, n-1$ to the 
cocycle dimension $n$.

\demo{Proof of {\rm Theorem 2.6}}  Suppose that $c\in \tZ^n(Q, \T)$ and $\AS c=1$. Let $\{\zdk:
k\in
\N\}$ be a sequence of generators of $Q$ and let
$$
M_m=\langle \zdone, \zdtwo, \cdots, \zdm\rangle, \quad m\in \N.
$$ 
The sequence $\{M_m\}$ is then increasing and $Q=\bigcup M_m$.
The triviality assumption $\AS c=1$ and Lemma 2.8 (i) yield  that the restriction $c_m$ of the cocycle $c$ to
each $M_m$ is a coboundary, so that there exists $f_m\in 
\tC^{n-1}(M_m, \T)$ such that
$$
c_m=\part_{M_m} f_m.
$$
The last lemma however allows us to choose the sequence $\{f_m\}$ in such a way that
each $f_m$ is an extension of the previous $f_{m-1}$. Hence the sequence $\{f_m\}$ gives a 
cochain $f\in \tC^{n-1}(Q, \T)$ such that $f|_{M_m}=f_m, m\in\N$, and therefore
$$
\partq f=c.
$$
Thus we conclude that $\Ker(\AS)\i \tB^n(Q, \T)$. The  inclusion,
$\Ker(\AS)\j \tB^n(Q, \T)$, was proven in
Lemma 2.2. Hence $\Ker(\AS)=\tB^n(Q, \T)$.  

Lemma 2.8 (ii) 
for  $\{M_m\}_{m\in\N}$ yields that the asymmetrization $\AS c$ of every $c\in \tZ^n(Q,\T)$ is a multi-character.

Set $\cdtxa=\AS c$ for an arbitrary cocycle $c\in \tZ^n(Q, \T)$. Then $c_{\text a}\in
X^n(Q, \T)$. Since
$Q$ is torsion free, the group $X^n(Q, \T)$ is indefinitely divisible. So
the $n!$-th power mapping: $\xi\in X^n(Q, \T)\mapsto \xi^{n!}\in X^n(Q, \T)$ is surjective. But
the asymmetrization $\AS$ on $X^n(Q, \T)$ is precisely the $n!$-th power. Hence there exists $\xi\in
X^n(Q, \T)$ such that $\AS \xi=\xi^{n!}=c_{\text a}$. Now we have $\AS (\xi\inv
c)=\xi^{-n!}\cdtxa=1$. Thus
$\xi\inv c \in \tB^n(Q, \T)$. Consequently, we conclude
$$\aligned
&\hskip.2in \tZ^n(Q, \T)=X^n(Q, \T)\tB^n(Q, \T);\\
 X^n&(Q, \T)\cap \tB^n(Q, \T)= X^n(Q, \T)\cap\Ker(\AS) \\
&=\{c\in X^n(Q, \T): c^{n!}=1\}.
\endaligned
$$
This completes the proof.
\QED
\enddemo

\proclaim{Corollary 2.10} If $G$ is a discrete abelian group, then the
asymmetrization of every $n$-cocycle $c\in\tZ^n(G, \T)$ is a multi-character, i.e., $\AS c\in
X^n(G, \T)$.
\endproclaim
\demo{Proof} Let $F$ be a large enough free abelian group so that there exists a
surjective homomorphism $\pi\!: F\mapsto G$. Consider the pullback $\pi^*(c)$ and its
asymmetrization, $\AS\pi^*(c)=\pi^*(\AS c)$. It follows from Theorem 2.6 that the
pull back $\pi^*(\AS c)$ is a multi-character of $F$, consequently the original
asymmetrization $\AS c$ is a multi-character of $G$. 
\QED
\enddemo

\head{\bf \S3. Universal Resolution for a Countable Discrete Abelian Group.}
\endhead

We discuss a universal resolution group for a countable discrete abelian group. We consider only the case that the abelian group under consideration has an infinitely many generators since the finitely generated case can be covered by the infinite generator case. 
 Let $G=\Zinf$ be the
free abelian group of a finite sequences of integers, i.e., every element $g\in G$ is of the form:
$$
g=(\ggdone, \ggdtwo, \cdots, \ggdi, \cdots, g_\ell, 0, 0, \cdots ), \quad \ggdi\in\Z,
$$
with $\ell=\ell(g)\in\N$, the last non-zero term of $g\in\Zinf$. With
$$
\adi=(0, 0,\cdots,0, \overset{\underset{\downarrow}\to{i}}\to {1}, 0, 0\cdots), \tag3.1
$$
every element $g\in\Zinf$ is written uniquely
$$
g=\sumd{i\in\N}\edi(g)\adi.\tag3.2
$$
We call $\{\adi: i\in\N\}$ the {\it standard basis} of $\Zinf$.
We also fix a subgroup $N$ of $G$ which is generated by a sequence $\lrbrace{\pdi\adi: i\in\N}$ with $\pdi\in\Z_+, i\in\N$. We will use the matrix:
$$
P=\pmatrix \pdone&0&0&\cdots\\
0&\pdtwo&0&\cdots\\
0&0&\pdthr&\cdots\\
\vdots&\vdots&\vdots&\ddots\\
\endpmatrix, \quad
N=P\Zinf.
$$

Let $M$ be the additive group of integer coefficient upper triangular matrices:
$$
M=\left\{m=\pmatrix 0&\mdonetwo&\mdonethr&\mdonefour&\cdots\\
0&0&\mdtwothr&\mdtwofour&\cdots\\
0&
0&0&\ddots&\cdots\\
\vdots&\vdots&\vdots&\ddots&\cdots
 \endpmatrix: \mdij\in\Z
\right\}
$$
and set 
$$
\edjk(m)=\mdjk,\quad j<k,\ m\in M.
$$
Let $\adi\wedge \adj, i<j,$ be the element of $M$ such that
$$
\edkell(\adi\wedge \adj)=\ddik\ddjell,
$$
i.e., the matrix with only $(i,j)$-component 1 and all others $0$, equivalently $\adi\wedge\adj, i<j$, is the $(i,j)$-matrix unit of $M$. 
Let
$\fndm$ be the $M$-valued second cocycle of $G$ defined by:
$$\aligned
\edjk&\lr{\fndm(g;h) }=\edj(g)\edk(h), \quad\ g, h\in G,\ 1\leq j<k;\\
\fndm(g; h)&=\pmatrix
0&\edone(g)\edtwo(h)&\edone(g)\edthr(h)&\edone(g)\edfour(h)&\cdots\\
0&0&\edtwo(g)\edthr(h)&\edtwo(g)\edfour(h)&\cdots\\
0&0&\ddots&\edthr(g)\edfour(h)&\cdots\\
\vdots&\vdots&\ddots&\vdots&\cdots
\endpmatrix.
\endaligned\tag3.3
$$
Let $H$ be the group extension of $G$ associated with 
$\fndm\in\tztw(G, M)$: 
$$\aligned
H&=M\times_\fndm G \quad \text{and}\quad L=M\times_\fndm N.
\endaligned
$$
The group operation in $H$ is given by:
$$\aligned
(m, g)(n, h)&=(m+n+\fndm(g; h), g+h),\quad (m, g), (n, h)\in H.
\endaligned
$$
The inverse $(m, g)\inv$ is given by:
$$
(m, g)\inv=(-m+\fndm(g, -g), -g)
$$
because
$$\aligned
(0, 0)&=(m, g)(m', g')=(m+m'+\fndm(g; g'), g+g')\\
&\quad g'=-g, \quad m'=-m+\fndm(g; g).
\endaligned
$$
To determine the commutater subgroup $[H, H]$, we compute the commutator:
$$\aligned
(m, g)&(n, h)(m, g)\inv(n, h)\inv, \hskip.5in  (m, g), (n, h)\in H,\\
&=(m, g)(n, h)(-m+\fndm(g; g), -g)(-n+\fndm(h; h); -h)\\
&=(m+n+\fndm(g, h), g+h)\\
&\hskip.5in\times
(-m-n+\fndm(g; g)+\fndm(h; h)+\fndm(g; h), -g-h)\\
&=\lr{\fndm(g; h)+\fndm(g; g)+\fndm(h; h)+\fndm(g; h)+\fndm(g+h; -(g+h)), 0}\\
&=\lr{\fndm(g; h)-\fndm(h;g), 0}\\
&=\lr{\sumd{j<k}(\edj(g)\edk(h)-\edj(h)\edk(g))\lr{\adj\wedge\adk}, 0}.
\endaligned
$$
This shows immediately the following:
\proclaim{Lemma 3.1} The commutator subgroup $[H, H]$ of $H$ is the center $M$ of $H$.
\endproclaim
\demo{Proof} From the computation above, it follows that for each pair $j<k$
$$
\fsdh(\adj)\fsdh(\adk)\fsdh(\adj)\inv\fsdh(\adk)\inv=\adj\wedge\adk,
$$
with $\fsdh$ the cross-section of $\pidz: (m, g)\in H\mapsto g\in G$ given by
$$
\fsdh(g)=(0, g)\in H, \quad g\in G,
$$
so that the commutator subgroup $[H, H]$ contains the generators $\adj\wedge \adk, j<k$, of $M$. Thus our
assertion follows.
\QED
\enddemo

\proclaim{Theorem 3.2} The pair $\{H, \pidz\}$ is a universal resolution of the third cocycle group $\tzthr(G,
\T)$ of $G$. Consequently, if $K$ is a countable discrete abelian group, then for any surjective
homomorphism $\pi\!\!: \Zinf\mapsto K$, the composed map $\pi_K=\pi\scirc \pidz\!\!: H\mapsto K$ 
makes the pair $\{H, \pi_K\}$ a universal resolution of the third cocycle group $\tzthr(K, \T)$.
\endproclaim
\demo{Proof} Since $\Zinf$ is a free abelian group on countably infinite generators, there exists a surjective
homomorphism from $G$ to any countable abelian group $G$. So it is sufficient to prove that 
$$
\pidzs\lr{\tzthr(G, \T)}\i \tbthr(H, \T).
$$
For each triplet $\xi, \eta, \z \in \Hom(G, \R)$, we define a multi-homomorphism, called the 
{\it tensor product} and denoted by $\xi\ot\eta\ot\z\in\tcthr(G, \R)$, as follows:
$$
\lr{\xi\ot\eta\ot\z}(g; h; k)=\xi(g)\eta(h)\z(k), \quad g, h, k\in G.
$$
Then the tensor product $\xi\ot\eta\ot\z$  generate the 
third cocycle group $\tzthr(G, \R)$ up to coboundary, i.e., 
$$
\Big\langle \Big\{\xi\ot\eta\ot\z: \xi, \eta, \z\in\Hom(G, \R)\Big\}\Big\rangle+\tbthr(G, \R)=\tzthr(G, \R).
$$
Now for each pair $\eta, \z\in\Hom(G, \R)$ we define a cochain $\Bdetaz\in\tcone(H, \R)$
$$
\Bdetaz(g)=\sumd{j< k}\eta(\adj)\z(\adk)\edjk(\mdz(g)), \quad g=(\mdz(g), \pidz(g))\in H.\tag3.4
$$
Then we have
$$\eightpoint\aligned
\Big(\parth &(\pidzs\xi\ot\Bdetaz)\Big)(\ggdone; \ggdtwo;
\ggdthree)=\xi(\pidz(\ggdtwo))\Bdetaz(\ggdthree)-\xi(\pidz(\ggdone)+\pidz(\ggdtwo))\Bdetaz(\ggdthree)\\
&\hskip1in
+\xi(\pidz(\ggdone))\Bdetaz(\ggdtwo\ggdthree)-\xi(\pidz(\ggdone))\Bdetaz(\ggdtwo)\\
&=-\xi(\pidz(\ggdone))\Bdetaz(\ggdthree)
+\xi(\pidz(\ggdone))\lr{\sumd{j<k}\eta(\adj)\z(\adk)\edjk(\mdz(\ggdtwo\ggdthree))}\\
&\hskip1.5in
-\xi(\pidz(\ggdone))\Bdetaz(\ggdtwo)\\
&=-\xi(\pidz(\ggdone))\Bdetaz(\ggdthree)\\
&\hskip.2in
+\xi(\pidz(\ggdone))\lr{\sumd{j<k}\eta(\adj)\z(\adk)\lr{\edjk(\mdz(\ggdtwo)+\mdz(\ggdthree)+\pidzs\fndm(\ggdtwo; \ggdthree)}}\\ 
&\hskip1.5in
-\xi(\pidz(\ggdone)\Bdetaz(\ggdtwo)\\
&=\xi(\pidz(\ggdone))\lr{\sumd{j<k}
\eta(\adj)\z(\adk)\edj(\pidz(\ggdtwo))\edk(\pidz(\ggdthree))}.
\endaligned
$$
Choosing $\xi, \eta, \z\in \Hom(G, \T)$ to be $\xi=\edi, \eta=\edj$ and $\z=\edk$ for $i<j<k$, we
obtain
$$\aligned
\pidzs(\edi\ot\edj\ot\edk)&=\parth \lr{\pidzs\edi\ot B_{\edj, \edk}}.
\endaligned
$$
Every third cocycle in $\tzthr(G, \T)$ is cohomologous to a cocycle $\cda\in \tzthr(G, \T)$ of the form:
$$
\cda(\ggdone; \ggdtwo; \ggdthree)=\explrtwopii{\sumdijk
\aaijk\edi(\ggdone)\edj(\ggdtwo)\edk(\ggdthree)}.\tag3.5
$$
So with $\bda\in\tctw(H, \T)$ defined by:
$$\aligned
\bda(\ggdone; \ggdtwo)&=\explrtwopii{\sumdijk \aaijk \edi(\pidz(\ggdone))B_{\edj,\edk}(\ggdthree)},
\endaligned\tag3.6
$$
we have
$$
\pidzs\cda=\parth \bda.\tag3.7
$$
Hence we get
$$
\pidzs\lr{\tzthr(G, \T)}\i \tbthr(H, \T),
$$
which concludes that the pair $\{H,\pidz\}$ is a universal resolution of $\tzthr(G, \T)$ and completes the proof.
\QED
\enddemo

\proclaim{Corollary 3.3} The $\mu$-part of every characteristic cocycle $(\la, \mu)\in \tZ(H, M,\allowmathbreak \T)$ is trivial.
\endproclaim
\demo{Proof} Since $M\triangleleft H$ is central, $\la$ is a bicharacter of $M\times H$, in particular $\la(m,
\cdot)$ is a character of $H$ for every $m\in M$. Hence it must vanish on the commutator subgroup, i.e.,
$\la(m, n)=1$ for ever $m, n\in M$. Thus $\mu\in \tztw(M, \T)$ is a coboundary.
\QED
\enddemo

Consider $(\la, \mu)\in \tZ(H, L,  \T)$ with $L=M\times_\fndm N$.  We
may and do assume the triviality $\mum=1$  of the restriction of $\mu$ to $M$. We then have the
corresponding crossed extension:
$$\CD
1@>>>\T@>>>E@>j>\underset{u}\to\longleftarrow>L@>>>1
\endCD
$$
The triviality of $\mum$ means that the cross-section $u$ is multiplicative on $M$, i.e., $u(mn)=u(m)u(n), m, n\in M$. Here we use the multiplicative group operation as $M$ sits in the noncommutative group $H$.

\proclaim{Lemma 3.4} If $\fsdh$ is a cross-section of the quotient map $\pidz: H\mapsto \Zinf=H/M$ with
$\fndm=\part \fsdh\in \tztw(\Zinf, M)$, then each characteristic cocycle in $\tZ(H, L, M, \T)$ is cohomologous to the one $(\la, \mu)\in \tZ(H, L, M, \T)$ such that\:
$$\gathered
\la(m; n\fsdh(h))=\la(m; \fsdh(h)), \quad m, n\in M,\ h\in \Zinf; \\
\mu(m\fsdh(g); n\fsdh(h))=\la(n; \fsdh(g))\mu\lr{\fsdh(g); \fsdh(h)}, \quad m,n\in M,\ g, h\in N.
\endgathered
$$

\endproclaim
\demo{Proof} In the crossed extension $E\in \X(\Hdm, L, M, \T)$ associated with $(\la, \mu)\in\tZ(\Hdm, L, M, \T)$:
$$\CD
1@>>>\T@>>>E@>>>L@>>>1,
\endCD
$$
we redefine the cross-section $u$ in the following way:
$$
u(m\fsdh(g))=u(m)u(\fsdh(g)), \quad m\in M, g\in N, 
$$
so that  $\mu(m; \fsdh(g))=1, m\in M, g\in N$. We now compute, for $m,n\in M,  h\in \Zinf$:
$$\aligned
\la(m; n\fsdh(h))u(m)&=\a_{n\fsdh(h)}\lr{u(m)}
=u(n)\a_{\fsdh(h)}\lr{u(m)}u(n)\inv\\
&=\la\lr{m; \fsdh(h)}u(mn)u(n)\inv\\
&=\la\lr{m; \fsdh(h)}u(m),\\
\endaligned
$$
for $g, h\in N$, we continue the computation:
$$\aligned
\mu\lr{m\fsdh(g); n\fsdh(h)}&u\lr{m\fsdh(g)n\fsdh(h)}\\
&=u\lr{m\fsdh(g)} u\lr{n\fsdh(h)}\\
&=u(m)u\lr{\fsdh(g)}u(n)u\lr{\fsdh(h))}\\
&=u(m)\a_{\fsdh(g)}\lr{u(n)}u\lr{\fsdh(g)}u\lr{\fsdh(h)}\\
&=\la\lr{n; \fsdh(g)}u(m)u(n)\mu\lr{\fsdh(g);\fsdh(h)}u\lr{\fsdh(g)\fsdh(h)}\\
&=\la\lr{n; \fsdh(g)}u(mn)\mu\lr{\fsdh(g); \fsdh(h)}
u\lr{\fsdh(g)\fsdh(h)}\\
&=\la\lr{n; \fsdh(g)}\mu\lr{\fsdh(g); \fsdh(h)}
u\lr{m\fsdh(g)n\fsdh(h)},
\endaligned
$$
and complete the proof.
\QED

\enddemo

\subhead\nofrills{\bf Groups $\pmb{G,\ \Hdm,\ \Gdm,\ \text{ and }\ \Qdm}$:\quad}
\endsubhead
First, we fix notations.
To work on the quotient group
$\Z/p\Z=\Z_p, p\in\N, p\geq 2$, we set
$${\eightpoint\aligned
\bracketp{i}&=i+p\Z\in \Zp;\quad
i=np+\bracep{i}, \ 0\leq \bracep{i}<p,\\
\etapp\left(\bracketp{i}, \bracketp{j}\right)&=\bracep{i}+\bracep{j}-\bracep{i+j}
=\cases
0\quad\text{if } \bracep{i}+\bracep{j}<p;\\
p\quad \text{if }  \bracep{i}+\bracep{j}\geq p.
\endcases\\
\endaligned}\tag3.8
$$
We shall call the $p\Z$-valued cocycle $\etapp\in \tztw(\Zp, p\Z)$ the {\it Gauss
cocycle}, which can be written in the following way:
$$
\etapp\left(\bracketp{i},
\bracketp{j}\right)=p\left(\left[\frac{i+j}p\right]-\left[\frac ip\right]-\left[\frac
jp\right]\right),\tag3.8$'$
$$
where $[x], x\in \R,$ is the Gauss symbol, i.e., the largest integer less than or equal to $x$.

Given a homomorphism $\txm$ of the group $G$ to  $\rt'z$ such that $\Ker(\txm)\j N$, we consider the group extension:
$$\gathered
\Gdm=\lrbrace{(g, s)\in G\times \R: \dstp=s+T'\Z=\txm(g)\in\rt'z},\\
\CD
0@>>>\Z@>n\to\zdzun>>\Gdm@>\pim>>G@>>>1,
\endCD
\endgathered
$$
where 
$$
\zdz=(0, T')\in\Gdm.
$$
Identifying $\txm$ with $\txm\scirc \pidz\in\Hom(H, \rt'z)$, we also form a group extension:
$$\gathered\aligned
\Hdm&=\lrbrace{(h, s)\in H\times \R: \txm(h)=\dstp\in\rt'z}\\
&=\lrbrace{(m, g, s)\in M\times G\times \R: \txm(g)=\dstp\in\rt'z},
\endaligned\\
\CD
0@>>>\Z@>n\to \zdzun>>\Hdm@>>>H@>>>1,
\endCD
\endgathered
$$
where the central element
$$
\zdz=\lr{1, T'}\in\Hdm
$$
appears in both $\Gdm$ and $\Hdm$. We hope that this abuse use of the same notation for two distinct elements  in the different groups
will not cause a headhach later: just like the zero elements in the ring theory.

By the assumption $N\i \Ker(\txm)$, the homomophism $\txm$ factors through the quotient group $Q=G/N$, so that it is also
viewed as a homomorphism of $Q\mapsto \rt'z$ and therefore we can form the group extension $\Qdm$ as before,
which sits on the following commutative diagram of  exact sequences:
$$\CD
@.@.1@.1\\
@.@.@VVV@VVV\\
@.0@>>>\fsm (N)@>>>N@>>>0\\
@.@VVV@VVV@VVV\\
0@>>>\Z@>>>\Gdm@>\pim>\underset{\fsdm}\to{\longleftarrow}>G@>>>1\\
@.@|@V\piq V\uparrow\fs V@V\piq V\uparrow\fs V\\
0@>>>\Z@>>>\Qdm@>\dpim>\underset{\dfsm}\to{\longleftarrow}>Q@>>>1\\
@.@VVV@VVV@VVV\\
@.0@.1@.1
\endCD
$$

From the assumption $\Ker\lr{\txm}\j N$, it follows that
 $\txm(\pdi\adi)=0$, so that there exists an integer $\qdi\in \Z, 0\leq \qdi<\pdi$
such that
$$\gathered
\txmdi=\bracett'{\txm(\adi)}=\frac{\qdi T'}\pdi\in \lr{\frac{T'}\pdi\Z},\\
\txm(\adi)=\dottxmdi=\txmdi+\T'\Z\in\rt'z.
\endgathered\tag3.9
$$
We set
$$\left.\aligned
\zdi&=(\adi, \txmdi)\in \Gdm, i\neq 0, \quad \zdz=(0, T')\in \Gdm,\\
 \fsdm(g)&=\sumd{i\in\N}\edi(g)\zdi
=\lr{g, \sumd{i\in\N}\edi(g)\txmdi}=\lr{g, \txn(g)},\\
&\txn(g)=\sumd{i\in\N}\edi(g)\txmdi, \quad g\in G.
\endaligned\right\}\tag3.10
$$
Then $\Gdm$ decomposes in the following way:
$$\left.\aligned
&\Gdm=\Z\zdz\oplus \fsdm(G)=\sumd{i\in\Ndz}{}^\oplus \Z\zdi,\ \text{where}\ \Ndz=\N\cup\lrbrace{0}, \\
\tilg&=\tedz(\tilg)\zdz+\sumd{i\in\N}\tedi(\tilg)\zdi\in \Gdm;\\
\tilg&=(g, s)=(0,\tedz(g, s)T')+\lr{\sumd{i\in\N}\tedi(\tilg)\adi,
\sumd{i\in\N}\tedi(\tilg)\txmdi}\\
&=\lr{0, \tedz(g, s)T'}+\sumd{i\in\N}\tedi(\tilg)\zdi;\\
\tedz(g, s)&=\frac{s-\txn(g)}{T'}\in\Z,\quad 
\tedi(g, s)=\edi(g), \quad i\in\N.
\endaligned\right\}\tag3.11
$$
In particular, if $g\in N$, we have
$$
g=(g, 0)=-\frac{\txn(g)}{T'}\zdz+\sumd{i\in\N}\edi(g)\zdi,
$$
so that
$$
\tedz(g)=-\frac{\txn(g)}{T'}\neq 0 \quad \text{unless}
\quad \txn(g)=\sumd{i\in\N}\edi(g)\txmdi=0.
$$
We then have 
$$
\txm(g)=\brackett'{\txn(g)}\in\rt'z.
$$
Setting
$$
\bdj=\pdj\adj,\quad j\in\N,
$$
we write every $g\in N$ uniquely in the form:
$$
g=\sumd{j\in\N}\frac{\edj(g)}{\pdj}\bdj
=\sumd{j\in\N}\edjn(g)\bdj,\tag3.12
$$
where 
$$
\edjn(g)=\frac{\edj(g)}{\pdj},
$$
and also in $\Hdm$ we have
$$
\bdj=\pdj\zdj-\pdj\mdj\zdz=\pdj\zdj-\qdj\zdz.\tag3.12$'$
$$

\vskip.1in
{\smc Remark.} The element $(\adi, 0)$ is NOT a member of $\Gdm$.
\vskip.1in

Next we define a cross-section $\dfsm: Q\mapsto \Qdm$ in such a way that
the following diagram commutes:
$$\CD
\Gdm@<\fsdm<<G\\
@AA\fs A@AA\fs A\\
\Qdm@<\dfsm<<Q
\endCD
$$
First, we set
$$\left.\aligned
&\dot g=g+N\in Q=G/N, \quad g\in G;\\
&\fs(q)=\sumd{i\in\N}\bracedpi{\edi(q)}\adi, \quad q\in Q;\\
&\dotadi=\piqm(\adi),\quad \dotzdi=(\dotadi, \txmdi)\\
\dfsm(q)&=\sumd{i\in\N}\bracedpi{\edi(q)}\dotzdi
=\sumd{i\in\N}\bracedpi{\edi(q)}\Big(\dotadi,\txmdi\Big)\\
&=\lr{q, \sumd{i\in\N}\bracedpi{\edi(q)}\txmdi} ;\\
&\fs(q, s)=(\fs(q), s)\in \Gdm,\quad (q, s)\in \Qdm.
\endaligned\right\}\tag3.13
$$
The cross-section $\fs: \Qdm\mapsto \Gdm$ gives rise to an $N$-valued
cocycle:
$$
\fnn=\partq\fs\in\tztw(\Qdm, N),\tag3.14
$$
which is given by:
$$\aligned
\fnn(\tqdone; \tqdtwo)&=\fs(\qdone, \sdone)+\fs(\qdtwo,
\sdtwo)-\fs(\qdone+\qdtwo, \sdone+\sdtwo)\\
&=(\fs(\qdone)+\fs(\qdtwo)-\fs(\qdone+\qdtwo), 0)\\
&=\lr{\sumd{i\in\N}{\edietan{\qdone}{\qdtwo}\adi, 0}}\\
&=\sumd{i\in\N}\lr{{\edietan{\qdone}{\qdtwo}\adi, 0}}\in N =N\times\{0\}.
\endaligned\nopagebreak
$$
for each pair $\tqdone=(\qdone, \sdone), \tqdtwo=(\qdtwo, \sdtwo)\in\Qdm$.

For each element 
$$
h=(m, g)\in H, \quad m\in M, g\in G,
$$
we write $m=\mdz(h)$ and $g=\pig(h)$. Then we have
$$
L=\piginv(N)
$$
and
$$
\mdz(gh)=\mdz(g)+\mdz(h)+\fndm(\pig(g); \pig(h)), \quad g, h\in H.
$$
For short, we write:
$$
\edij(\tilg)=\edij(\mdz(g))\ \text{for }\tilg=(\mdz(g), g,s)\in \Hdm, i, j\in\N.
$$
With 
$$
\sh(g)=(0, g)\in H\quad \text{for each}\ g\in G,
$$
we have
$$
\fndm(g; h)=\sh(g)+\sh(h)-\sh(g+h)=\partg \sh(g; h), \ g, h\in G.
$$
With
$$
\dot\fs=\sh\scirc \fs,
$$
we obtain a cross-section $\dot\fs$ of $\piq\scirc \pig: H\mapsto Q=H/L$,
which gives rise to an $L$-valued second cocycle $\fnl\in\tztw(Q, L)$:
$$\aligned
\fnl&(\qdone; \qdtwo)=\dfs(\qdone)\dfs(\qdtwo)\dfs(\qdone+\qdtwo)\inv, 
\quad \qdone, \qdtwo\in Q,\\
&=\sh\Big(\fs(\qdone)\Big)\sh\Big(\fs(\qdtwo)\Big)
\sh\Big(\fs(\qdone+\qdtwo)\Big)\inv\\
&=\fndm\Big(\fs(\qdone); \fs(\qdtwo)\Big)\sh\Big(\fs(\qdone)+\fs(\qdtwo)
\Big)
\sh\Big(\fs(\qdone+\qdtwo)\Big)\inv\\
&=\fndm\Big(\fs(\qdone);
\fs(\qdtwo)\Big)\sh\Big(\fnn(\qdone; \qdtwo)+\fs(\qdone+\qdtwo)\Big)
\sh\Big(\fs(\qdone+\qdtwo)\Big)\inv\\
&=\fndm\Big(\fs(\qdone);
\fs(\qdtwo)\Big)\fndm\Big(\fnn(\qdone; \qdtwo);
\fs(\qdone+\qdtwo)\Big)\inv \sh\Big(\fnn(\qdone; \qdtwo)\Big).
\endaligned\tag3.15
$$
We further compute the $j,k$-components and $k$-components:
$$\left.\aligned
\edjk&\Big(\fndm(\fs(\qdone); \fs(\qdtwo))\Big)
=\edj(\fs(\qdone))\edk(\fs(\qdtwo))\\
&\hskip1.05in
=\bracedpj{\edj(\qdone)}\bracedpk{\edk(\qdtwo)},\\
\edjk&\Big(\fndm(\fnn(\qdone; \qdtwo); \fs(\qdone+\qdtwo)\Big)
=\edj\lr{\fnn(\qdone; \qdtwo)}\edk\lr{\fs(\qdone+\qdtwo)}\\
&\hskip1.05in
=\edjetan{\qdone}{\qdtwo}\bracedpk{\edk(\qdone+\qdtwo)},\\
\edk&\Big(\sh\Big(\fnn(\qdone; \qdtwo)\Big)\Big)
=\edketan{\qdone}{\qdtwo}.
\endaligned\right\}\tag3.16
$$
Since
$$\aligned
\Hdm&=M\times_{\pims(\fnm)} \lr{\sumd{i\in\N}{}^\oplus \Z \zdi\oplus \Z\zdz},
\endaligned
$$
for each $h=(m, g)\in H$, we set
$$\aligned
\fsdm(h)&=(m, \fsdm(g))=\lr{m, \sumd{i\in\N}\edi(g)\zdi}
=\lr{m, g, \sumd{i\in\N}\edi(g)\txmdi},
\endaligned\tag3.17
$$
and we identify $\ell=(m, Pg)\in L$ with $(m, Pg, 0)\in \Hdm$, so that
$L$ is a subgroup of $\Hdm$, while $H$ is not.

\head{{\pmb\S4}. \bf The Characteristic Cohomology Group $\pmb{\La(\Hdm, L, M, \T)}$.}
\endhead

Since $H$ is a universal resolution group for $G=\Zinf$, every third
cohomology class $[c]\in\ththr(G, \T)$ is of the form $[c]=\dhjr[\la, \mu]$
for some $[\la, \mu]\in \La(H, M, \T)$. So every outer action $\da$ of
$G$ on a factor $\sM$ of type \threel\  comes from an action $\a$ of $H$,
i.e., the outer action $\da$ is given by 
$$
\da_g=\a_{\sh(g)}, \quad g\in G.\tag4.1
$$
But the action $\a$ of $H$ does not
give rise to an action of $H$ on the reduced (discrete) core $\tMd$. Instead,
the action $\a$ of $H$ on $\sM$ gives rise naturally to an action, denoted by the
same notation $\a$, of  $\Hdm$ on $\tMd$ where 
$$
\txm(h)=\mod(\a_h) \in\rt'z,\ h\in H.
$$
If $N=\da\inv(\cntr(\sM)))\i G$, then $L=\a\inv(\cntr(\sM))$. We make a basic assumption on the subgroup $N$:
$$
N=PG=P\Zinf.
$$
In the case that $G$ is finitely generated free abelian group, the fundamental structure theorem for finitely generated abelian groups guarantees that every subgroup of $G$ is of this form.

We study first the characteristic cohomology group $\La(\Hdm, L, M, \T)$
and modified HJR-map $\d: \La(\Hdm, L, M, \T)\mapsto \thmsout(G, N, \T)$.

We introduce a series of notations first:
$$\left.\aligned
&\Ndz=\N\cup\{0\}=\Z_+;\\
\Deldz=\{(i,j,k)\in \Ndzuthr&: i<j<k\}\cup\{(i,i,k)\in\Ndzuthr; i<k\}\\
&
\cup\{(k, i,k)\in\Ndzuthr: i<k\},\\
&\Delta=\Deldz\cap \N^3.
\endaligned\right\}\tag4.2
$$
For each $g\in \Hdm$, let $\mdz(g)$ be the $M$-component of $g$, i.e.,
$$
\mdz(g)=g\sh(\pig(g))\inv\in M, \quad g \in \Hdm.\tag4.3
$$
We regard $\edi$ and $\edjk$ as functions defined on $\Hdm$  by fixing the
coordinate system:
$${\eightpoint\aligned
\tilg=\lr{\sumd{1\leq j<k}\edjk(g)(\adj\wedge \adk), \sumd{i\in\Ndz}
\tedi(\tilg)\zdi}\in\Hdm,\ \text{with}\ g=\pim(\tilg)\in H.
\endaligned}
\tag4.4
$$
We then introduce a cochain
$\Bdjk\in\tcone(\Hdm, \R)$ defined by the following:
 $$\aligned
\Bdjk(h)=\cases -\edjk(\mdz(h))\quad&\text{for }j<k;\\
-\dfrac{(\edj\edj)(h)}2\quad&\text{for }j=k;\\
\edkj(\mdz(h))-\lr{\edj\edk}(h)\quad &\text{for } j>k,
\endcases\quad h\in\Hdm.
\endaligned\tag4.5
$$
The cochain enjoys the property:
$$
\parth \Bdjk= \pidzs(\edj\ot\edk)\quad \text{for } j, k\in\N.\tag4.6
$$
We continue to  define the following cochains
for each $a\in \R^\Ndzuthr$:
$$\allowdisplaybreaks\eightpoint\align
\xxdaijk&=\aaijk\edjk\ot\edi+\aajik\edik\ot\edj+\aakij\edij\ot\edk,\\
\xxdaik&=\aaiik\edik\otimes\edi+\aakik\edik\otimes\edk;\\
\yydaijk&=\aaijk\Big(\Bdij\ot\edk+\edk\ot\Bdji-
\Bdik\ot\edj-\edj\ot\Bdki\Big)\\
&\hskip.2in
+\aajik\Big(\Bdji\ot\edk+\edk\ot\Bdij-\Bdjk\ot\edi-\edi\ot\Bdkj\Big)\\
&\hskip.2in
+\aakij\Big(\Bdki\ot\edj+\edj\ot\Bdik-\Bdkj\ot\edi-\edi\ot\Bdjk\Big),\\
\yydaik&=\aaiik\lr{\Bdii\otimes
\edk+\edk\otimes\Bdii-\Bdik\otimes\edi-\edi\otimes
\Bdki}\\
&\hskip.2in
+\aakik\lr{\Bdki\otimes \edk+\edk\otimes \Bdik- \Bdkk\otimes
\edi-\edi\otimes \Bdkk},\\
&Z(\cdots)(g; h)=Y(\cdots)(\mdz(h); g);\\
\zzdaijk&=\aaijk\Big(\edj\ot\edik-\edk\ot\edij\Big)\\
&\hskip-.2in
+\aajik\Big(\edk\ot\edij+\edi\ot\edjk\Big)
+\aakij\Big(\edj\ot\edik-\edi\ot\edjk\Big),\\
\zzdaik&=\aaiik\edi\ot\edik+\aakik\edk\otimes\edik;\\
\ffdijk&=2(\edi\edj)\ot\edk-3\edi\ot(\edj\edk)+\edj\ot(\edi\edk)\\
&\hskip1.5in
-2(\edi\edk)\ot\edj
-\edk\ot(\edi\edj),\\
\uudaijk&=\frac16\Big( \aaijk\ffdijk+\aajik\ffdjik+\aakij\ffdkij\\
&\hskip1in
-\Asaijk\ffdijk\Big),\\
\uudaik&=-{\aaiik}\Bdii\otimes \edk+\aakik\lr{
\Bdkk\otimes\edi-\edk\otimes(\edi\edk)},\\
&\vvdaijk=\zzdaijk+\pigs\uudaijk,\\
&\hskip.1in
\vvdaik=\zzdaik+\pigs\uudaik.
\endalign
$$

The infinite summations:
$$\left.\aligned
\xxda&=\sumdijk\xxdaijk+\sumd{i<k}\xxdaik;\\
\yyda&=\sumdijk\yydaijk+\sumd{i<k}\yydaik;\\
\uuda&=\sumdijk\uudaijk+\sumd{i<k}\uudaik;\\
\vvda&=\sumdijk\vvdaijk+\sumd{i<k}\vvdaik;\\
\zzda&=\sumdijk\zzdaijk+\sumd{i<k}\zzdaik\\
\endaligned\right\}\tag4.7
$$
will become all finite sums  as soon as variables from $M$ or $\Hdm$ are fed in.
So  no divergence problem in the infinite sums will occur. 

The cochain $\ffdijk$ relates basic cocycles $\edi\ot\edj\ot\edk$ and the asymmetric tri-character:
$$\aligned
\detdijk&=\lr{\edi\ot\edj\ot\edk+\edj\ot\edk\ot\edi+\edk\ot\edi\ot\edj}\\
&\hskip.5in-
\lr{\edj\ot\edi\ot\edk+\edi\ot\edk\ot\edj+\edk\ot\edj\ot\edi}\\
&=\edi\wedge\edj\wedge\edk
\endaligned
$$
in the following way:
$$
\detdijk=\partl \ffdijk+6\edi\ot\edj\ot\edk, \quad i<j<k,
\tag4.8
$$
which can be confirmed by a direct computation.

Let  $Z$ be the set of all pairs $(a, b)$ of functions $a$: $ (i,j,k)\in
\N^3\mapsto\aaijk\in \R$ and $b$: $(i,j)\in\Ndzutw\mapsto \bbij\in\R$
satisfying the following requirements:
\roster
\item"a)" The requirements on the parameter $a$ is given by:
$$\left.\gathered
\aaijk=0\text{ for } j\geq k \text{ and }  a(0,j,k)=0\quad \text{for every } j,
k\in \Ndz,\\
\Asaijk=\aaijk-\aajik+\aakij\\
\in \lr{\frac1{\gcd(\pdi, \pdj, \pdk)}\Z}.
\endgathered\right\}\tag4.9Z-a
$$
\item"b)" The requirements for the parameter $b$ is given by: 
$$\left.\aligned
\bbij\pdj&- \bbiz\qdj \in \Z\quad \text{for }\quad i,
j\in\N, i<j,\\
&\bbzj=0,\quad j\in\Ndz.
\endaligned\right\}\tag4.9Z-b
$$
\endroster
Let $Z_a$ be the set of $a\in \R^{\N^3}$ satisfying the above requirement (4.9Z-a)  and $Z_b$ be the set of all
$b\in \R^{\Ndzutw}$ with the properties of (4.9Z-b). So we have
$$
Z=\Zda\oplus \Zdb.
$$

Let $B$ be the subgroup of $Z$ consisting of all those $(a, b)\in Z$ such
that 
\roster
\item"a)" The coboundary condition on the parameter $a$ is given by:
$$\left.\gathered
\aaijk, \aakij, \aajik\in \Z\quad \text{if } i<j<k,\\
\aaiik\in{2\Z}\quad \text{if } i<k; \quad 
\aakik\in{2\Z}\quad \text{if } i<k,\\
\endgathered\right\}\tag4.9B-a
$$
\item"b)" The coboundary condition on the parameter $b$ is given by:
$$\left.\gathered
\frac\bbij\pdi+\frac\bbji\pdj\in
\lr{\frac1\pdi\Z}+\lr{\frac1\pdj\Z}=\lr{\frac1{\lcm(\pdi,
\pdj)}\Z}, i<j,\\
\bbiz\in\Z,\quad \bbii\in \Z, \quad i\in\N.
\endgathered\right\}\tag4.9B-b
$$
\endroster
Respectively, let $\Bda$ (resp. $\Bdb$) be the set of
all those $b\in \R^\Ndzutw$ satisfying the requirement of (4.9B-a) (resp. (4.9B-b)). Thus we have
$$
B=\Bda\oplus \Bdb.
$$
and set
$$\gathered
\La=\Lada\oplus \Ladb, \quad \Lada=\tzda/\tbbda, \quad
\Ladb=\tzdb/\tbbdb,\\
H_a=\Zda/\Bda, \quad H_b=\Zdb/\Bdb.
\endgathered
$$
With $D(i,j,k)=\gcdlr{\pdi, \pdj, \pdk}$ for each triplet $i<j<k, i,j,k\in\N$, we set
$$\gathered
\tZ_a(i,j,k)=\lrbrace{(u, v, w)\in\R^3: u-v+w\in \lr{\frac1{D(i,j,k)}\Z}},\\
\tB_a(i,j,k)=\Z\oplus\Z\oplus\Z,\\
\text{where } u=\aaijk,\ v=\aajik,\ w=\aakij.
\endgathered
$$
For a pair $i<k, i, k\in\N$, we set
$$
\tZ_a(i,k)=\lrbrace{(x, y)\in \R^2}=\R\oplus \R, \quad \tB_a(i, k)=(2\Z)\oplus (2\Z),
$$
where $x=\aaiik$ and $y=\aakik$. We then naturally define:
$$\gathered\aligned
\Lada(i,j,k)&=\tZ_a(i,j,k)/\tB_a(i,j,k)\\
&\cong \lr{\left.\lr{\frac1{D(i,j,k)}\Z}\right/\Z}\oplus \R/\Z\oplus \R/\Z, \quad i<j<k;\\
\endaligned\\
\Lada(i,k)=\tZ_a(i,k)/\tB_a(i,k)=\R/\lr{2\Z}\oplus
\R/\lr{2\Z}, \quad i<k.
\endgathered
$$
Here the above second isomorphism can be seen easily by considering the matrix:
$$
A=\pmatrix1&-1&1\\0&1&0\\0&0&1\endpmatrix
\in\SL(3, \Z).
$$
For each ordered pair $i<j, i, j\in\N$, we define
$${\eightpoint\left.\gathered
\tzzdbij=\lrbrace{(x, u, y, v)\in\R^4: \pdj x-\qdj u\in\Z, \pdi y-\qdi v\in \Z};\\
\tbbdbij=
\lrbrace{(x, u, y, v)\in\tzdbij: \pdj x+\pdi y\in \Ddij\Z,
 u, v\in \Z},\\
\tzzdbii=\lrbrace{z=(x, u)\in \R^2: \pdi x -\qdi u\in\Z},\quad
\tbbdbii=\Z\oplus \Z,\\
\Ladb(i,j)=\tzzdbij/\tbbdbij, \quad \Ladb(i,i)=\tzzdbii/\tbbdbii,
\endgathered\right\}}\tag4.10
$$
where
$$
\Ddij=\gcd(\pdi,\pdj).
$$

{\smc Definition 4.1.} To each $(a, b)\in Z$ we associate a cochain $\lr{\ladab, \muda}$ defined by the following:
$$\left.\aligned
\ladab(g; h)&=\explrtwopii{\lr{\yyda+X_{\AS a}}(g; h)}\\
&\hskip.5in\times
\explrtwopii{\!\!\sumd{i\in\N, j\in\Ndz}  {\bbij\edin(g)\tedj(h)}},\\
\etada(g; h)&=\explrtwopii{\yyda(g; h)},\\
\muda(g; h)&=\exp\Big(\twopii\vvda(g; h)\Big)\\
&=\ladab(\mdz(h); g)\exp\Big(\twopii\uuda(\pig(g); \pig(h))\Big),
\endaligned\right\}\tag4.11
$$
for each $(g, h)\in L\times \Hdm$. In the case that $b=0$ (resp. $a=0$) we denote the corresponding cochains by $(\lada, \muda)$ (resp.
$\ladb$).

\proclaim{Theorem 4.2} {\rm a)} The cochain $\lr{\lada, \muda}$ is a characteristic cocycle in $\tZ(\Hdm, \allowmathbreak L, M, \T)$
and the correspondence $a\in Z_a\mapsto (\lada, \muda)\in \tzzda$ gives 
 the following
commutative diagram of exact sequences\:
$${\eightpoint\CD
@.@.@.0\\
@.@.@.@VVV\\
0@>>>\Bda@>>>a\in \Zda@>>>\lrbracket{a}\in H_a@>>>0\\
@.@VVV@VVV@VVV\\
0@>>>\tbbda@>>>(\lada, \muda)\in\tzzda@>>>\lrbracket{\lada, \muda }\in \Lada@>>>1\\
@.@.@.@VVV\\
@.@.@.0
\endCD}\tag4.12-a
$$

{\rm b)} The correspondence $b\in Z_b\mapsto \ladb\in \tZ(\Hdm, L, M, \T)$ gives the following commutative diagram of exact sequences\:
$$\CD
@.@.@.0\\
@.@.@.@VVV\\
0@>>>\Bdb@>>>b\in \Zdb@>>>\lrbracket{b}\in H_b@>>>0\\
@.@VVV@VVV@VVV\\
0@>>>\tbbdb@>>>(\ladb, 1)\in\tzzdb@>>>\lrbracket{\ladb}\in\Ladb@>>>1\\
@.@.@.@VVV\\
@.@.@.0
\endCD\tag4.12-b
$$
 
{\rm c)} The characteristic cohomology group 
$$
\La(\Hdm, L, M, \T)=\Lada\oplus \Ladb
$$ 
has further fine structure\:
\roster
\item"i)" The group $\Lada$ has the Cartesian product decomposition\:
$$\left.\gathered
\Lada=\prod_{i<j<k}\Lada(i,j,k)\oplus\prod_{i<j}\Lada(i,j);\\
\Lada\lr{i,j,k}\cong\Z_{D(i,j,k)}\oplus \R/\Z\oplus\R/\Z, \\
D(i,j,k)=\gcdlr{\pdi,\pdj, \pdk};\\
\Lada(i,j)\cong\R/\lr{2\Z}\oplus \R/\lr{2\Z}.\\
\endgathered\right\}\tag4.13-a
$$
\item"ii)" The group $\Ladb$ has the fiber product decomposition into the family  $\lrbrace{\Ladb(i,j): i,j\in\N}$ and each group $\Ladb(i,j)$ is described as follows\:
$$\left.\gathered
\Ladb(i,j)\cong\Z/\lr{\gcdlr{\pdi,\pdj,\qdi, \qdj}\Z}\oplus
\lr{\R/\Z}\oplus \lr{\R/\Z}, i<j,\\
\Ladb(i,i)\cong \Z/\lr{\gcdlr{\pdi, \qdi}\Z}\oplus \lr{\R/\Z}.
\endgathered\right\}\tag4.13-b
$$
The group $\Ladb(i,j)$ and  {\rm(}resp. $\Ladb(i,i)${\rm)} is equipped with three {\rm(}resp. one{\rm)} homomorphisms\: 
$$\left.\gathered
\pidij: \Ladb(i,j)\mapsto \left.\lr{\frac 1{D(i,j)}\Z}\right/\Z,\\
\pidijui: \Ladb(i,j)\mapsto \R/\Z,  \quad \pidijuj: 
\Ladb(i, j)\mapsto \R/\Z,\\
\pidiui: \Ladb(i,i)\mapsto \R/\Z,
\endgathered\right\}\tag4.14
$$
such that for each $z=(x,u,y,v)\in \tzdbij$
$$\left.\gathered
\pidij\lr{\lrbracket{\la_z}}=\lrbracket{\mdij\lr{x\rdji+y\rdij}-\ndij\lr{u\sdji+v\sdij}}_\Z,\\
\pidijui\lr{\lrbracket{\la_z}}=\lrbracket{u}_\Z,\quad
\pidijuj\lr{\lrbracket{\la_z}}=\lrbracket{v}_\Z,\\
\pidii\lr{\lrbracket{\la_z}}=\lrbracketd{\pdi x-\qdi u}{\Z}, \quad
\pidiui\lr{\lrbracket{\la_z}}=\lrbracketd{u}{\Z},
\endgathered\right\}\tag4.15
$$
where
$$\left.\gathered
D(i,j)=\gcdlr{\pdi, \pdj, \qdi, \qdj},\\  \Ddij
=\gcdlr{\pdi, \pdj}, \quad \Edij=\gcdlr{\qdi, \qdj},\\
\rdij=\frac\pdi\Ddij, \quad \rdji=\frac\pdj\Ddij\quad
\sdij=\frac \qdi\Edij, \quad \sdji=\frac\qdj\Edij,\\
\mdij=\frac \Ddij {D(i,j)}, \quad\ndij=\frac \Edij {D(i,j)},\\
\qdi\wdij+\qdj\wdji=\Edij, \quad 
\xdij \Ddij+\ydij \Edij=D(i,j).
\endgathered\right\}\tag4.16
$$
The group $\Ladb$ is the fiber product of $\lrbrace{\Ladb(i,j): i,j\in\N}$ relative to the maps
$\lrbrace{\pidijui, \pidijuj, \pidiui: i,j\in\N}$ in the sense that $\Ladb$ is the group of all those 
$\ladb\in \prod_{(i,j)\in\N^2}\Ladb(i,j)$ such that
$$\aligned
\pidijui\lrbracket{\ladb(i,j)}&=\pidiui\lrbracket{\ladb(i,i)}=\pidkiui\lrbracket{\ladb(k,i)}, \quad i,j,k\in\N.
\endaligned\tag4.17
$$
\endroster
\endproclaim

We will prove the theorem in several steps.

First, we  observe that the asymmetrization of $\ffdijk$ is given by:
$$\aligned
\AS\ffdijk&=2(\edi\edj)\wedge\edk-3\edi\wedge(\edj\edk)
+\edj\wedge(\edi\edk)\\
&\hskip1.5in
 -2(\edi\edk)\wedge\edj -\edk\wedge(\edi\edj)\\
&=3\Big((\edj\edk)\wedge \edi-(\edi\edk)\wedge
\edj+(\edi\edj)\wedge \edk\Big).
\endaligned\tag4.18
$$

\proclaim{Lemma 4.3} {\rm i)} The difference of $\xxda$ and $\yyda$ on
$M\times \Hdm$ is given by\:
$$\aligned
\xxda-\yyda=X_{\AS a}\quad \text{ on}\quad M\times \Hdm.
\endaligned
$$
In particular if the following integers
$$
\edij(m)\edk(g),\ \edjk(m)\edi(g),\ \edik(m)\edj(g),\ \edjk(m)\edi(g)
$$
are all divisible by $\gcdpijk$, then we have for each $a\in Z$
$$\aligned
\yydaijk(m; g)\equiv\xxdaijk(m; g)\quad \mod \Z, \quad m \in M, g\in \Hdm.
\endaligned
$$
Therefore, if either $g\in L$ or $m\in L\wedge \Hdm$, then the following
congruence holds\:
$$\aligned
\xxdaijk(m; g)&\equiv \yydaijk(m; g) \quad \mod \Z;\\
\xxdaijk(\hdone\wedge g; \hdtwo)&\equiv \yydaijk(\hdone\wedge g; \hdtwo)
\quad \mod \Z\ 
\endaligned\tag4.19
$$
for each $ \hdone, \hdtwo\in\Hdm.$

{\rm ii)} For every $m\in M$ and $g\in \Hdm$ and $i<k$ we have
$$
\xxdaik(m; g)=\yydaik(m; g).\tag4.20
$$
\endproclaim
\demo\nofrills{Proof. \quad}  i) We simply compute for $i<j<k$:
$$\aligned
\lr{\xxdaijk-\yydaijk}&(m; g)\\
&\hskip-1in
=\aaijk\edjk(m)\edi(g)+\aajik\edik(m)\edj(g)\\
&\hskip.5in
+\aakij\edij(m)\edk(g)\\
&\hskip-.7in
-\aaijk\Big(\edik(m)\edj(g)-\edij(m)\edk(g)\Big)\\
&\hskip-.5in
-\aajik\Big(\edij(m)\edk(g)+\edjk(m)\edi(g)\Big)\\
&\hskip-.3in
-\aakij\Big(\edik(m)\edj(g)-\edjk(m)\edi(g)\Big)\\
&\hskip-1in
=\Big(\aaijk-\aajik+\aakij\Big)\\
&\hskip-.7in\times
\Big(\edjk(m)\edi(g)-\edik(m)\edj(g)+\edij(m)\edk(g)\Big).\\
\endaligned
$$
Thus we conclude
$$\aligned
\lr{\xxda-\yyda}(m; g)=X_{\AS a}(m; g), \quad m\in M, g \in \Hdm.
\endaligned
$$

ii)  The assertion follows from an easy direct computation.
\QED
\enddemo

\proclaim{Lemma 4.4} If $a\in \RuDelt$ is asymmetric modulo
$\lr{\dfrac1{\pdi\pdj\pdk}\Z}$ in the sense that\:
$$
\Asaijk=\aaijk-\aajik+\aakij\in\lr{\frac1{\pdi\pdj\pdk}\Z}\tag4.21
$$
for each triplet $i<j<k$, then the cochain $\muda$ of {\rm (4.11)}, i.e.,  
$$
\muda(g; h)=\explrtwopii{\vvda(g; h)}, \quad g, h\in L,
$$
is a second cocycle $\muda\in\tztw(L, \T)$.
\endproclaim
\demo\nofrills{Proof.} \quad  Observing
$$
\lr{\partl \muda}(\ggdone; \ggdtwo; \ggdthree)=\explrtwopii{\partl \vvda(\ggdone; \ggdtwo; \ggdthree)}, \quad \ggdone, \ggdtwo,
\ggdthree\in L,
$$
we compute the coboundary of $\vvda$:
$$\aligned
\partl&\vvdaijk=\partl \zzdaijk+\partl\uudaijk\\
&=\aaijk\lr{\edj\ot\edi\ot\edk-\edk\ot\edi\ot\edj}\\
&\hskip.5in
+\aajik\lr{\edk\ot\edi\ot\edj+\edi\ot\edj\ot\edk}\\
&\hskip1in
+\aakij\lr{\edj\ot\edi\ot\edk-\edi\ot\edj\ot\edk}\\
&\hskip.5in
+\frac16\partl\Big( \aaijk\ffdijk+\aajik\ffdjik+\aakij\ffdkij\\
&\hskip1.5in
-\Asaijk\ffdijk\Big)\\
\endaligned
$$
$$\aligned
&=\aaijk\lr{\edj\ot\edi\ot\edk-\edk\ot\edi\ot\edj}\\
&\hskip.5in
+\aajik\lr{\edk\ot\edi\ot\edj+\edi\ot\edj\ot\edk}\\
&\hskip1in
+\aakij\lr{\edj\ot\edi\ot\edk-\edi\ot\edj\ot\edk}\\
&\hskip.5in
+\frac16\Big(\aaijk(\detdijk-6\edi\ot\edj\ot\edk)\\
&\hskip.7in
+\aajik(\detdjik-6\edj\ot\edi\ot\edk)\\
&\hskip.9in
+\aakij(\detdkij-6\edk\ot\edi\ot\edj)\\
&\hskip1.1in
-\Asaijk(\detdijk-6\edi\ot\edj\ot\edk)\Big)\\
&\equiv
-\Asaijk\Big(\edi\ot\edj\ot\edk-\edj\ot\edi\ot\edk+\edk\ot\edi\ot\edj\Big)\\
&\equiv 0\quad \mod \quad \Z\quad\text{on }L\times L\times L,
\endaligned
$$
since $\edi\ot\edj\ot\edk$ takes values in $\pdi\pdj\pdk\Z$ on $L\times
L\times L$. Also we have
$$\aligned
\partl &\vvdaik=\partl  \zzdaik+\partl\uudaik\\
&=\aaiik\edi\ot\edi\ot\edk+\aakik\edk\ot\edi\ot\edk
-\aaiik\edi\ot\edi\ot\edk\\ &\hskip.5in
+\aakik\Big(\edk\ot\edk\ot\edi-\edk\ot(\edi\ot\edk+\edk\ot\edi)\Big)\\
&=0.
\endaligned
$$
Hence $\muda$ is a second cocycle on $L$.
\QED
\enddemo

\proclaim{Lemma 4.5} 
{\rm i)} For every $(a, b)\in\tZ$, the pair 
$\{\ladab, \muda\}$  is a characteristic cocycles
in $\tZ(\Hdm, L, M, \T)$.

{\rm ii)} Every characteristic cocycle $(\la, \mu)\in
\tZ(\Hdm, L, M, \T)$ is cohomologous to some $(\ladab, \muda)$.

{\rm iii)} The characteristic cocycle $\{\ladab, \muda\}\in \tZ(\Hdm, L, M,
\T)$ is a coboundary if and only if $(a, b)\in \tB$.

\endproclaim
\demo\nofrills{Proof.}\quad  i) We first check the cocycle identities for each
$g, \ggdone, \ggdtwo \in L$ and $h, \hdone, \hdtwo\in \Hdm$:
$$\align
\Big(\Big(\partl \ot\id\Big)\ladab\Big)(\ggdone; \ggdtwo;
h)&=\frac{\muda(h\inv
\ggdone h; h\inv\ggdtwo h)}{\muda(\ggdone; \ggdtwo)}\\
&=\ladab(\ggdtwo\wedge h;  \ggdone);\tag a\\
\Big(\Big(\id\ot\parthm\Big)\ladab\Big)(g; \hdone; \hdtwo)
&=\frac1 {\ladab(g\wedge \hdone; \hdtwo)}\\
&=\ladab(\hdone\wedge g; \hdtwo);
\tag b\\
&\hskip-1in
\ladab(g; h)=\frac{\muda(h; h\inv g h)}{\muda(g; h)}, \quad g, h\in L.
\tag c
\endalign
$$
Second, we compute for $\ggdone, \ggdtwo\in L$ and $h\in\Hdm$:
$$\aligned
\xxdaijk&(\ggdtwo\wedge h; \ggdone)\\
&=\aaijk\edjk(\ggdtwo\wedge h)\edi(\ggdone)
+\aajik\edik(\ggdtwo\wedge h)\edj(\ggdone)\\
&\hskip1in
+\aakij\edij(\ggdtwo\wedge h)\edk(\ggdone)\\
\endaligned
$$
$$\aligned&=\aaijk\edi(\ggdone)\Big(\edj(\ggdtwo)\edk(h)-\edk(\ggdtwo)\edj(h)\Big)\\
&\hskip.5in
+\aajik\edj(\ggdone)\Big(\edi(\ggdtwo)\edi(h)-\edk(\ggdtwo)\edi(h)\Big)\\
&\hskip.5in
+\aakij\edk(\ggdone)\Big(\edi(\ggdtwo)\edj(h)-\edj(\ggdtwo)\edi(h)\Big)\\
&=\Big[\aaijk\edi\ot\Big(\edj\ot\edk-\edk\ot\edj\Big)
+\aajik\edj\ot\Big(\edi\ot\edi-\edk\ot\edi\Big)\\
&\hskip.5in
+\aakij\edk\ot\Big(\edi\ot\edj-\edj\ot\edi\Big)
\Big]\Big(\ggdone; \ggdtwo; h\Big).
\endaligned
$$
On the other hand, we have
$$\left.\aligned
\Big(\Big(\partl \ot \id\Big)\yydaijk\Big)
&=\aaijk\Big(\edi\ot\edj\ot\edk
-\edi\ot\edk\ot\edj\Big)\\
&\hskip.2in+\aajik\Big(\edj\ot\edi\ot\edk
-\edj\ot\edk\ot\edi\Big)\\
&\hskip.4in
+\aakij\Big(\edk\ot\edi\ot\edj-\edk\ot\edj\ot\edi\Big).
\endaligned\right\}\tag4.22
$$
Since
$$
X_{\AS a}(i,j,k)\lr{\ggdtwo\wedge h; \ggdone}\equiv
0\quad \mod \Z,
$$
Lemma 4.3 yields, for each $\ggdone, \ggdtwo\in L, h\in
\Hdm$, the following:
$$\aligned
\Big((\partl\ot\id)&\yydaijk\Big)(\ggdone; \ggdtwo; h)
=\xxdaijk(\ggdtwo\wedge h; \ggdone)\\
&\equiv\yydaijk(\ggdtwo\wedge h; \ggdone) \ \mod \Z.
\endaligned
$$
Similarly, we have
$$
\lr{\lr{\partl \ot \id}\yydaik}\lr{\ggdone, \ggdtwo; h}\equiv\yydaik\lr{\ggdtwo\wedge h; \ggdone}\quad \mod \Z.
$$
Next, we have
$$\left.\eightpoint{\aligned
\xxdaijk&(\hdone\wedge g; \hdtwo)
=\aaijk\edjk(\hdone\wedge g)\edi(\hdtwo)\\
&
+\aajik\edik(\hdone\wedge g)\edj(\hdtwo)
+\aakij\edij(\hdone\wedge g)\edk(\hdtwo)\\
&\hskip-.5in
=\Big[\aaijk\Big(\edk\ot\edj-\edj\ot\edk\Big)\ot\edi
+\aajik\Big(\edk\ot\edi-\edi\ot\edk\Big)\ot\edj\\
&
+\aakij\Big(\edj\ot\edi-\edi\ot\edj\Big)\ot\edk\Big]\Big(g;
\hdone;\hdtwo\Big);\\
\Big(\id\ot\parthm\Big)&\yydaijk(g; \hdone; \hdtwo)
=\Big[\aaijk\Big(\edk\ot\edj\ot\edi- \edj\ot\edk\ot\edi\Big)\\
&\hskip.5in
+\aajik\Big(\edk\ot\edi\ot\edj-\edi\ot\edk\ot\edj\Big)\\
&\hskip.5in
+\aakij\Big(\edj\ot\edi\ot\edk-\edi\ot\edj\ot\edk\Big)\Big] \Big(g;
\hdone;\hdtwo\Big)\\
&=\xxdaijk(\hdone\wedge g; \hdtwo)
\endaligned}\right\}\tag4.23
$$
and
$$\aligned
\Big(\id\ot\parthm\Big)\xxdasaijk=0.
\endaligned
$$
Hence Lemma 4.3 again yields, for each $g\in L, \hdone,\hdtwo\in\Hdm$, that:
$$\aligned
\Big(\id&\ot\parthm\Big)\Big(\yydaijk+\xxdasaijk\Big)(g; \hdone; \hdtwo)\\
&\equiv \Big(\yydaijk+\xxdasaijk\Big)
(\hdone\wedge g; \hdtwo)\ \mod\ \Z.
\endaligned
$$
Similarly, we get the following:
$$\gathered
\lr{\lr{\id\ot\parthm}\yydaik}\lr{g, \hdone; \hdtwo}=\yydaik\lr{\hdone\wedge g; \hdtwo},\quad g\in L, \hdone, \hdtwo\in\Hdm,\\
\xxdasa(i,k)=0.
\endgathered
$$
Thus so far we have established the formulae (a) and (b). 

Now we move on
to (c). Fixing $g, h\in L$, we compute the right hand side of (c):
$$\aligned
\frac{\muda(h; h\inv g h)}{\muda(g; h)}&=\frac{\muda(h; (g\wedge
h)g)}{\muda(g; h)}=\lada(g\wedge h; h)\frac{\muda(h; g)}{\muda(g;h)}\\
&\hskip-.5in
=\lada(g\wedge h; h)\frac{\muda\Big(\mdz(h)\sh(h); \mdz(g)\sh(g)\Big)}
{\muda\Big(\mdz(g)\sh(g); \mdz(h)\sh(h)\Big)}\\
&\hskip-.5in=\explrtwopii{\xxda(g\wedge h; h)}(\explrtwopii{\AS \vvda(h; 
g)})\\
&\hskip-.5in=\explrtwopii{\lr{\yyda+\xxdasa}(g\wedge h;
h)}(\explrtwopii{\AS
\vvda(h;  g)}.
\endaligned
$$
Next we prove the following:
$$\aligned
\lada(\sh(g); \sh(h))=\lada(g\wedge  h; h)
(\AS \muda)(\sh(h); \sh(g)), \quad g, h\in N.
\endaligned
$$
First we observe that 
$$
\xxdasa(g; h)\equiv 0\quad \mod \Z\quad \text{for }\quad g, h \in L.
$$
So for the proof of (c), the term $\xxdasa$ can be ignored. 
With this fact in mind, we compute:
$$\aligned
\xxdaijk&(g\wedge h; h)\\
&=\aaijk\edjk(g\wedge h)\edi(h)+\aajik\edik(g\wedge h)\edj(h)\\
&\hskip.5in
+\aakij\edij(g\wedge h)\edk(h)\\
&=\Big[\aaijk\Big(\edj\ot(\edk\edi)-\edk\ot(\edj\edi)\Big)\\
&\hskip.5in
+\aajik\Big(\edi\ot(\edk\edj)-\edk\ot(\edi\edj)\Big)\\
&\hskip1in
+\aakij\Big(\edi\ot(\edj\edk)-\edj\ot(\edi\edk)\Big)\Big]\Big(g; h\Big),
\endaligned
$$
and also
$$\aligned
\xxdaik&(g\wedge h; h)=\aaiik\Big(\edi(g)\edk(h)-\edk(g)\edi(h)\Big)\edi(h)\\
&\hskip.7in
+\aakik\Big(\edi(g)\edk(h)-\edk(g)\edi(h)\Big)\edk(h)\\
&=\aaiik\Big(\edi\ot(\edi\edk)-\edk\ot e_i^2\Big)\Big(g; h\Big)\\
&\hskip1in
+\aakik\Big(\edi\ot e_k^2-\edk\ot (\edi\edk)\Big)\Big(g; h\Big).
\endaligned
$$
Next we determine the asymmetrization of $\uudaijk$ based on (4.18):
$$\allowdisplaybreaks\eightpoint\align
\AS&\uudaijk\\
&=\frac16\Big(
\aaijk\AS\ffdijk+\aajik\AS\ffdjik+\aakij\AS\ffdkij\\
&\hskip1in
-\Asaijk\AS\ffdijk\Big)\\
&=\frac12\Big[\aaijk\Big((\edj\edk)\wedge \edi-(\edi\edk)\wedge
\edj+(\edi\edj)\wedge\edk\Big)\\
&\hskip.5in
+\aajik\Big((\edi\edk)\wedge\edj-(\edj\edk)\wedge\edi
+(\edi\edj)\wedge\edk\Big)\\
&\hskip.5in
+\aakij\Big((\edi\edj)\wedge\edk-(\edj\edk)\wedge\edi
+(\edi\edk)\wedge\edj\Big)\\
&\hskip.5in
-\Big(\aaijk-\aajik+\aakij\Big)\\
&\hskip1in\times
\Big((\edj\edk)\wedge
\edi-(\edi\edk)\wedge\edj+(\edi\edj)\wedge \edk\Big)\Big]\\
&=\frac12\Big[\aajik\Big((\edi\edk)\wedge\edj-(\edj\edk)\wedge\edi
+(\edi\edj)\wedge\edk\Big)\\
&\hskip.5in
+\aakij\Big((\edi\edj)\wedge\edk-(\edj\edk)\wedge\edi
+(\edi\edk)\wedge\edj\Big)\\
&\hskip.5in
+\Big(\aajik-\aakij\Big)\\
&\hskip1in\times
\Big((\edj\edk)\wedge
\edi-(\edi\edk)\wedge\edj+(\edi\edj)\wedge \edk\Big)\Big]\\
&=-\aakij(\edj\edk)\wedge\edi +\aakij(\edi\edk)\wedge\edj\\
&\hskip.5in
+\aajik(\edi\edj)\wedge\edk.
\endalign
$$
Hence we get
$$\left.\eightpoint{\aligned
\AS\uudaijk&=-\aakij\Big((\edj\edk)\ot\edi-\edi\ot(\edj\edk)\Big)\\
&\hskip.5in
+\aakij\Big((\edi\edk)\ot\edj-\edj\ot(\edi\edk)\Big)\\
&\hskip1in
 +\aajik\Big((\edi\edj)\ot\edk-\edk\ot(\edi\edj)\Big).
\endaligned}\right\}\tag4.24
$$
We also check the asymmetrization of $\uudaik$:
$$\aligned
\AS
\uudaik&=\aaiik\edk\wedge\Bdii+\aakik\lr{\Bdkk\wedge\edi
-\edk\wedge(\edi\edk)}\\
&=\frac\aaiik2\Big(e_i^2\ot\edk-\edk\ot e_i^2\Big)
+\frac\aakik2\Big(\edi\ot e_k^2-e_k^2\ot\edi\Big)\\
&\hskip.5in
+\aakik\Big((\edi\edk)\ot\edk-\edk\ot(\edi\edk)\Big).
\endaligned
$$
We then combine these with the above computations for $\xxdaijk$, paying attention to the order of variables in the first term and the second term
\footnote{In the first term,  the variables $g$ and $h$ appear in this order, but in the second term they appear in the opposite order.}:
$$\allowdisplaybreaks\align
\xxdaijk&(g\wedge h; h)+\AS\uudaijk(\sh(h); \sh(g))\\
&\hskip-.5in
=\aaijk\Big(\edj\ot(\edk\edi)-\edk\ot(\edj\edi)\Big)\\
&\hskip.5in
+\aajik\Big(\edi\ot(\edk\edj)-\edk\ot(\edi\edj)\Big)\\
&\hskip1in
+\aakij\Big(\edi\ot(\edj\edk)-\edj\ot(\edi\edk)\Big)\\
&
+\Big(\aakij(\edj\edk)\wedge\edi
-\aakij(\edi\edk)\wedge\edj\\
&\hskip1in
 -\aajik(\edi\edj)\wedge\edk\Big)\\
&=\aaijk\Big(\edj\ot(\edk\edi)-\edk\ot(\edj\edi)\Big)\\
&\hskip.5in
+\aajik\Big(\edi\ot(\edk\edj)-\edk\ot(\edi\edj)-(\edi\edj)\wedge\edk\Big)\\
&\hskip.5in
+\aakij\Big(\edi\ot(\edj\edk)-\edj\ot(\edi\edk)+(\edj\edk)\wedge\edi\\
&\hskip2in
-(\edi\edk)\wedge\edj\Big)\\
&=\aaijk\Big(\edj\ot(\edk\edi)-\edk\ot(\edj\edi)\Big)\\
&\hskip.5in
+\aajik\Big(\edi\ot(\edk\edj)-\edk\ot(\edi\edj)\\
&\hskip1.5in
-(\edi\edj)\ot\edk+\edk\ot(\edi\edj)\Big)\\
&\hskip.5in
+\aakij\Big(\edi\ot(\edj\edk)-\edj\ot(\edi\edk)\\
&\hskip1.5in
+(\edj\edk)\ot\edi-\edi\ot(\edj\edk)\\
&\hskip1.5in
-(\edi\edk)\ot\edj+\edj\ot(\edi\edk)\Big)\\
&=\aaijk\Big(\edj\ot(\edk\edi)-\edk\ot(\edj\edi)\Big)\\
&\hskip.5in
+\aajik\Big(\edi\ot(\edk\edj)
-(\edi\edj)\ot\edk\Big)\\
&\hskip.5in
+\aakij\Big((\edj\edk)\ot\edi
-(\edi\edk)\ot\edj\Big).
\endalign
$$
and
$$\allowdisplaybreaks\align
\xxdaik&(g\wedge h; h)+\AS\uudaik(\sh(h); \sh(g))\\
&=\aaiik\Big(\edi\ot(\edi\edk)-\edk\ot e_i^2\Big)\\
&\hskip.2in
+\aakik\Big(\edi\ot e_k^2-\edk\ot (\edi\edk)\Big)\\
&\hskip.4in
+\frac\aaiik2\Big(\edk\ot e_i^2-e_i^2\ot\edk\Big)
+\frac\aakik2\Big(e_k^2\ot\edi-\edi\ot e_k^2\Big)\\
&\hskip.6in
+\aakik\Big(\edk\ot(\edi\edk)-(\edi\edk)\ot\edk\Big)\\
&=\aaiik\lr{\edi\ot(\edi\edk)-\left.\left.\frac12\right(\!\edk\ot
e_i^2+e_i^2\ot\edk\right)}\\
&\hskip.2in
+\aakik\lr{\left.\left.\frac12\right(\edi\ot
e_k^2+e_k^2\ot\edi\right)-(\edi\edk)\ot\edk}.
\endalign
$$
We now compare these with $\yydaijk$:
$$\aligned
\yydaijk&(\sh(g); \sh(h))\\
&=\aaijk\Big(\edj\ot(\edi\edk)-\edk\ot(\edi\edj)\Big)\\
&\hskip.5in
+\aajik\Big(\edi\ot(\edj\edk)-(\edi\edj)\ot\edk\Big)\\
&\hskip1in
+\aakij\Big((\edj\edk)\ot\edi-(\edk\edi)\ot\edj\Big)\\
\endaligned
$$
$$\aligned&=\xxdaijk(g\wedge h; h)+\AS\uudaijk(\sh(h); \sh(g))\\
&\equiv\yydaijk(g\wedge h; h)+\AS\uudaijk(\sh(h); \sh(g)),
\endaligned
$$
and 
$$\aligned
\yydaik&(\sh(g); \sh(h))\\
&=\Big(\aaiik\lr{\Bdii\otimes
\edk+\edk\otimes\Bdii-\Bdik\otimes\edi-\edi\otimes
\Bdki}\\
&\hskip.2in
+\aakik\lr{\Bdki\otimes \edk+\edk\otimes \Bdik- \Bdkk\otimes
\edi-\edi\otimes \Bdkk}\Big)\\
&=\aaiik\lr{\edi\ot(\edi\edk)-\left.\left.\frac12\right(\edk\ot
e_i^2+e_i^2\ot\edk\right)}\\
&\hskip.4in
+\aakik\lr{\left.\left.\frac12\right(\edi\ot
e_k^2+e_k^2\ot\edi\right)-(\edi\edk)\ot\edk}\\
&=\xxdaik(g\wedge h; h)+\AS\uudaik(\sh(h); \sh(g))\\
&=\yydaik(g\wedge h; h)+\AS\uudaik(\sh(h); \sh(g)).
\endaligned
$$
Therefore, we have
$$\aligned
\ladab(\sh(g); \sh(h))&=\ladab(g\wedge h; h)\frac{\muda(\sh(h);
\sh(g))} {\muda(\sh(g); \sh(h))}\\
\endaligned
$$
Since we have
$$\aligned
\yyda&(mg; nh)=\yyda(m; h)+\yyda(g; n)+\yyda(g; h)
\endaligned
$$
for every $m,n\in M$ and $g, h\in \Hdm$,  we get, for each $m, n\in M$ and
$g, h\in N$,
$$\aligned
\ladab&(m\sh(g); n\sh(h))\\
&=\ladab(m; \sh(h))\ladab(\sh(g); n)
\ladab(\sh(g); \sh(h))\\
&=\frac{\ladab(m; \sh(h))}{\ladab(n; \sh(g))}\ladab(g\wedge h; h)
\frac{\muda(\sh(h); \sh(g))}{\muda(\sh(g); \sh(h))}\\
&=\frac{\muda(n\sh(h); (n\sh(h))\inv m\sh(g)n\sh(h))} {\muda(m\sh(g);
n\sh(h))}.
\endaligned
$$
This proves the cocycle identity (c).  Consequently $\{\ladab, \muda\}$ is  a characteristic cocycle in $\tZ(\Hdm, L, M, \T)$.

ii) Suppose that $(\la, \mu)\in \tZ(\Hdm, L, M, \T)$.
Since $M$ is central in
$\Hdm$, the $\la$-part is a bicharacter on $M\times \Hdm$, so that there exists
$a=\{\aaijk\}\in\R^\Delta$ such that
$$
\la(m; h)=\explrtwopii{\sumddijk\aaijk\edjk(m)\edi(h)}, \quad m\in M, h\in\Hdm.
$$
As $[\Hdm, \Hdm]=M$, for each fixed $m\in M$ the character $\la(m; \cdot)$  on
$\Hdm$ must vanish on $M$, i.e.,
$$
\la(m; n)=1, \quad m, n\in M.
$$
Thus the restriction $\mum$ of the second cocycle $\mu$ to $M$ is a coboundary.
Hence, replacing $\mu$ by a cohomologous cocycle if necessary, we may and do
assume that $\mum=1$. Now consider the corresponding $E\in\X(\Hdm, L, M,
\T)$:
$$\CD
1@>>>\T@>>>E@>j>\underset{\sj}\to\longleftarrow >L@>>>1.
\endCD
$$
Redefining the cross-section $\sj$ as
$$
\sj(m\sh(g))=\sj(m)\sj(\sh(g)), \quad m\in M, \ g\in N,
$$
we may and do assume that $\mu(m; g)=1, m\in M, g \in L$. Now we compute
the second cocycle $\mu$ with $m, n\in M$ and $g, h\in L$:
$$\aligned
\mu(mg; nh)&\sj(mgnh)=\sj(mg)\sj(ng)=\sj(m)\sj(g)\sj(n)\sj(h)\\
&=\sj(m)\la(n; g)\sj(n)\sj(g)\sj(h)\\
&=\la(n; g)\mu(g; h)\sj(m)\sj(n)\sj(gh)\\
&=\la(n; g)\mu(g; h)\sj(mngh)=\la(n; g)\mu(g; h)\sj(mgnh),
\endaligned
$$
which gives
$$
\mu(mg; nh)=\la(n; g)\mu(g; h), \quad m,n\in M, g, h\in L.
$$
In particular, we have
$$
\mu(g; h)=\la(\mdz(h); g)\mu(\sh(\pig(g)); \sh(\pig(h)), \quad g, h\in L,
$$
where 
$$
\mdz(h)=h\sh(\pig(h))\inv\in M.
$$
Now with $\ggdone, \ggdtwo, \ggdthree\in N$, we compute the coboundary:
$$\left.\aligned
1&=(\partl \mu)(\sh(\ggdone); \sh(\ggdtwo); \sh(\ggdthree))\\
&=\frac{\mu(\sh(\ggdtwo); \sh(\ggdthree))\mu(\sh(\ggdone);
\sh(\ggdtwo)\sh(\ggdthree))} 
{\mu(\sh(\ggdone)\sh(\ggdtwo); \sh(\ggdthree))\mu(\sh(\ggdone);
\sh(\ggdtwo))}\\
&=\frac{\mu(\sh(\ggdtwo); \sh(\ggdthree))\mu(\sh(\ggdone); \fndm(\ggdtwo;
\ggdthree) \sh(\ggdtwo+\ggdthree))} 
{\mu(\fndm\lr{\ggdone; \ggdtwo}\sh(\ggdone+\ggdtwo); \sh(\ggdthree))\mu(\sh(\ggdone);
\sh(\ggdtwo))}\\
&=\la(\fndm(\ggdtwo;
\ggdthree); \sh(\ggdone)) \frac{\mu(\sh(\ggdtwo);
\sh(\ggdthree))\mu(\sh(\ggdone); \sh(\ggdtwo+\ggdthree))} 
{\mu(\sh(\ggdone+\ggdtwo);
\sh(\ggdthree))\mu(\sh(\ggdone);
\sh(\ggdtwo))}.
\endaligned\right\}\tag4.25
$$
Thus the cocycle $\cda\in\tzthr(N, \T)$ given by:
$$\aligned
\cda(\ggdone; \ggdtwo; \ggdthree)&=\la(\fndm(\ggdtwo; \ggdthree); \ggdone)\\
&=\explrtwopii{\sumddijk \aaijk\edjk(\fndm(\ggdtwo; \ggdthree))\edi(\ggdone)}\\
&=\explrtwopii{\sumddijk\aaijk\edi(\ggdone)\edj(\ggdtwo)\edk(\ggdthree)}
\endaligned
$$
is a coboundary in $\tzthr(N, \T)$. Thus we get, for every 
$\ggdonedotsthree\in N$,
$$\aligned
1&=\Big(\AS\cda\Big)(\ggdonedotsthree)\\
&=\explrtwopii{\sumddijk\aaijk\sumd{\sig\in\Pi(i,j,k)}\sign(\sig)
\edi\lr{\ggd{\sig(i)}}
\edj\lr{\ggd{\sig(j)}}\edk\lr{\ggd{\sig(k)}}}\\
&=\explrtwopii{\sumddijk \aaijk\detdijk(\ggdone; \ggdtwo; \ggdthree)}\\
&=\explrtwopii{\sumd{(i,j,k)\in\Delta}\Asaijk\detdijk(\ggdonedotsthree)}.
\endaligned
$$
Therefore the coefficient $a=\{\aaijk\}\in\R^\Delta$  is asymmetric in the sense of Lemma 4.4, so that it gives the second cocycle
$\muda\in\tztw(L, \T)$:
$$\aligned
\muda=\explr{\twopii\vvda}.
\endaligned
$$
Then the cocycle
$\mu\mudainv\in\tztw(L, \T)$ falls in the subgroup $\pigs\lr{\tztw(N,
\T)}\mathbreak\i
\tbtw(L, \T)$ because
$$\aligned
\mu(m\sh(g); n\sh(h))&=\la(n; \sh(g))\mu(\sh(g); \sh(h))\\
&=\frac{\muda(m\sh(g); n\sh(h))}{\muda(\sh(g); \sh(h))}
\mu(\sh(g); \sh(h))\\
&=\frac{\mu(\sh(g); \sh(h))}{\muda(\sh(g); \sh(h))}\muda(m\sh(g);
n\sh(h));\\
\mudainv\mu&=\pigs\scirc\shs(\mu\mudainv)\in\pigs\Big(\tztw(N, \T)\Big).
\endaligned
$$
 Thus there exists a cochain $f\in\tcone(L, \T)$ such that
$$
\muda(g; h)=\mu(g; h)\frac{f(g)f(h)}{f(gh)}, \quad g, h\in L.
$$
Since $1=\mu(m; h)=\muda(m; h), m\in M, h\in L$, we have
$$
f(mh)=f(m)f(h). 
$$
Since $(\partdone f)(m; h)=1, m\in M, h\in \Hdm$, we have
$$
\part f(\la, \mu)=(\la, \muda).
$$

Next we look at one of the cocycle identities: 
$$\aligned
\la(\ggdone\ggdtwo; h)&=\la(\ggdone; h)\la(\ggdtwo; h)\frac{\muda(\ggdone;
\ggdtwo)}{\muda(h\inv \ggdone h; h\inv \ggdtwo h) }\\
&=\frac1{\la(\ggdtwo \wedge h; \ggdone)}\la(\ggdone; h)\la(\ggdtwo; h),
\quad \ggdone, \ggdtwo \in L, h\in\Hdm,\\
&=\la(h\wedge \ggdtwo; \ggdone)\la(\ggdone; h)\la(\ggdtwo; h)\\
&=\explrtwopii{\sumddijk\aaijk\edi(\ggdone)\edjk(h\wedge \ggdtwo)}\la(\ggdone;
h)\la(\ggdtwo; h),
\endaligned
$$
which gives the following partial coboundary condition:
$$\aligned
(\partl\ot\id)\la=\explrtwopii{\sumddijk\aaijk\edi\ot(\edj\ot\edk-\edk\ot\edj)}.
\endaligned
$$
Another cocycle identity:
$$\aligned
\la(g; \hdone\hdtwo)&=\la(g; \hdone)\la(\hdoneinv g\hdone; \hdtwo),\quad g \in L,
\hdone, \hdtwo\in\Hdm,\\
&=\la(g\wedge \hdone; \hdtwo)\la(g;\hdone)\la(g; \hdtwo)\\
&=\explrtwopii{\sumddijk\aaijk\edjk(g\wedge \hdone)\edi(\hdtwo)}\la(g; \hdone)
\la(g; \hdtwo)
\endaligned
$$
gives the second partial coboundary condition:
$$\aligned
(\id\ot\parthm)\la=\explrtwopii{\sumddijk\aaijk (\edk\ot\edj-\edj\ot\edk)\ot\edi}.
\endaligned
$$
Setting
$$\aligned
\etada=\explrtwopii{\yyda},
\endaligned
$$
we obtain, by (4.22) and (4.23), 
$$\aligned
(\partl\ot\id)\la=(\partl\ot\id)\etada; \quad (\id\ot\parthm)\la=(\id\ot\parthm)\etada.
\endaligned
$$
Therefore the cochain $\overline\eta_a\la=\chi$ is a bicharacter on $L\times\Hdm$.
Since $M=[\Hdm, \Hdm]$, the bicharacter $\chi$ vanishes on $L\times M$, i.e., $\la(m; g)=\etada(m; g), m\in M, g\in L$. Thus we get
$$\aligned
1&=\la(m; g)\overline\eta_a(m; g)
=\explrtwopii{\xxda(m; g)-\yyda(m; g)}\\
&=\explrtwopii{\xxdasa(m; g)}=\ladasa(m; g),
\endaligned
$$
which is equivalent to the following fact:
$$\aligned
\Asaijk\in\lr{\frac1\gcdpijk\Z}.
\endaligned
$$
Thus we conclude the cocycle condition (4.9Z-a) on the parameter $\lrbrace{\aaijk}$.  Therefore the coefficient $a\in \R^\Delta$ satisfies the
requirement for the element $(a, 0) \in Z$. Therefore it follows from (i) that $(\ladaz,
\muda)\in \tZ(\Hdm, L, M, \T)$. Then 
the cocycle identity (c) for $(\ladaz, \muda)$ yields that
$$\aligned
\la(g; h)&=\frac{\muda(h; h\inv g h)}{\muda(g; h)}=\ladaz(g; h)=\etada(g; h),\quad g, h\in L.
\endaligned
$$
Thus the bicharacter $\chi=\overline\eta_a\la$ on $L\times
\Hdm$ vanishes on $L\times L$. Since Lemma 4.3(i) yields for each $m\in
M, h\in \Hdm$ that
$$\aligned
\chi(m; h)&=\la(m; h)\overline\eta_a(m; h)=\lada(m; h)\overline\eta_a(m; h)\\
&=\lad{\AS a}(m; h)=\explrtwopii{\sumddijk \Asaijk\edjk(m)\edi(h)},
\endaligned
$$
we conclude that $\chi$ is of the form:
$$\aligned
\chi(g; h)&=\chidz(\pig(g); \pig(h))
\explrtwopii{\sumdijk\Asaijk\edjk(g)\edi(h)}
\endaligned
$$
for $g\in L, h\in \Hdm$ where $\chidz$ is a bicharacter on $N\times \Gdm$
and $\pig\!: \Hdm \mapsto \Gdm$ the quotient map with $M=\Ker(\pig)$. 
We choose $\bbij\in\R$ so that
$$
\explrtwopii{\bbij}=\chidz(\bdi; \zdj), i\in\N, j\in\Ndz.
$$
Then we must have
$$\aligned
1&=\chidz(\bdi; \bdj)=\chidz(\bdi;
\pdj\zdj-\qdj\zdz)=\explrtwopii{\bbij\pdj-\bbiz\qdj},
\endaligned
$$
so that $\bbij\in\R, i\in\N, j\in\Ndz$, satisfies the following condition:
$$\gathered
\bbij\pdj\equiv \bbiz\qdj\quad \mod\ \Z, \quad i,j\in\N.\\
\endgathered
$$
Hence $\chidz$ is written in the form:
$$\eightpoint\aligned
\chidz(g; \tilh)&=\explrtwopii{{\sumd{i,j\in\N}\bbij
\edin(g)\edj(\tilh) +\sumd{i\in\N}b(i, 0)\edin(g)\tedz(\tilh)}}
\endaligned
$$
for each pair $g\in N$ and $\tilh\in\Hdm$, where the coefficients $\bbij$
satisfies the requirements:
$$\aligned
\bbij\pdj&-\bbiz\qdj\in \Z, \quad i,j\in \N, \quad \bbzi=0,\quad i\in\Ndz.\\
\endaligned
$$
Consequently the pair $(a, b)$ is a member of $Z$ and we conclude that
$(\la, \mu)$ is cohomologous to the characteristic cocycle $(\ladab, \muda) \in
\tZ(\Hdm, L, M, \T)$.

iii) Suppose $(\la, \mu)=(\ladab, \muda)=\part f$
with $f\in \tcone(L, \T)$. Since $\mum=1$ and $\muda(m; g)=1, m\in M,
g\in L$, we have
$$
f(mg)=f(m)f(g), \quad m\in M, g\in L,
$$
so that the restriction  $f|_M$ of $f$ to $M$ is of the form:
$$
\ffdc(m)=\explrtwopii{\sumd{1\leq i<j}\ccij\edij(m)}, \quad m\in M.
$$
Since $M$ is central in $\Hdm$, we have for every pair $(m, g)\in M\times
\Hdm$
$$\aligned
1&=\frac{\ffdc(g\inv mg)}{\ffdc(m)}=\la(m; g)=\explrtwopii{\sumddijk\aaijk\edjk(m)\edi(g)},
\endaligned
$$
which yields the integrality condition on $a$:
$$
\aaijk\in\Z \quad \text{for every}\ (i,j,k)\in\Delta.
$$
Hence  $\yydaijk$ and $\xxdaijk$ are both integer valued, so
that 
$$
\lada(m; h)=1,\quad m\in M, h \in\Hdm.
$$
Since $\chi=1$ on $L\times L$, 
for every
$g, h\in L$ we have
$$\aligned
1&=\la_{0,b}(g; h)=\lada(\mdz(g); h)\la_{0,b}(\sh(g); h)=\la(g; h)\\
&=\frac{f(h\inv gh)}{f(g)}=\ffdc(g\wedge h);\\
&\ccij\in\lr{\frac1{\pdi\pdj}\Z}, \quad i,j\in\N.
\endaligned
$$
This computation also shows that
$$
\la_{0,b}(g; h)=\ffdc(g\wedge h), \quad g\in L, h\in\Hdm.
$$
Furthermore, we have for each $m, n\in M$ and $g, h\in L$
$$\aligned
\muda(mg; nh)&=\ladab(n; g)\muda(g; h)=\muda(g; h),
\endaligned
$$
so that $\muda$ is of the form: $\muda=\pigs(\tilde \mu)$ with
$$\aligned
\tmuda(g; h)=\explrtwopii{\uuda(g; h)}, \quad g, h\in N.
\endaligned
$$
Since $\lada(\fndm(\ggdtwo; \ggdthree); \ggdone)=1, \
\ggdone, \ggdtwo, \ggdthree\in\Hdm$, we have
$\tmuda\in\tztw(N,\T)$ by (4.25). We first compute for each $g, h\in L$:
$$\aligned
(\AS \muda)(g; h)&=\frac{f(g)f(h)}{f(gh)}
\frac{f(hg)}{f(g)f(h)}=\frac{f(hgh\inv h)}{f(gh)}=\ffdc(h\wedge g)
=1.
\endaligned
$$
Since $\AS\uudaijk$ is also integer
valued, we have
$$\aligned
\AS\muda&=\explrtwopii{\sumd{i<k}{\AS\uudaik}}\\
&
=\explrtwopii{\sumd{i<k}\lr{\frac\aaiik2\Big(e_i^2\ot\edk-\edk\ot
e_i^2\Big)}}\\
&\hskip.5in\times
 \explrtwopii{\sumd{i<k}\frac\aakik2\Big(\edi\ot
e_k^2-e_k^2\ot\edi\Big)}\\
& =1.
\endaligned
$$
Thus we get
$$
\aaiik,\quad \aakik\in 2\Z, \quad \text{and }\quad
\uudaik\equiv 0\quad \mod\ \Z.
$$
Consequently, $\tmuda$ is a coboundary as a
member of $\tztw(N, \T)$. Hence there exists a cochain $\tilde
f\in \tcone(N, \T)$ such that
$$\aligned
\frac{f(g)f(h)}{f(gh)}&=\muda(g; h)=\tmuda(\pig(g); \pig(h))=\frac{\tilde
f(\pig(g))\tilde f(\pig(h))}{\tilde f(\pig(gh))}.
\endaligned
$$
Thus $f$ is of the form: 
$$\gathered
f(g)=\ffdc(\mdz(g))f(\sh(g))=\chi(g)\tilde f(\pig(g)),\quad g\in L;\\
\ffdc(m)=\chi(m), \quad m\in M.
\endgathered
$$
where $\chi\in\Hom(L, \T)$. Since 
$$\aligned
L/[L, L]&\cong M/PMP\oplus N,
\endaligned
$$
the homomorphism $\chi$ is of the form:
$$\aligned
\chi(g)&=\explrtwopii{{\sumd{j<k}\ccjk \edjk(g)+\sumd{k\in\Ndz}
\cck \tedk(g)}}, \quad g\in L,
\endaligned
$$
where 
$$
\ccij\in \lr{\dfrac 1{\pdi\pdj}\Z},\ i<j, \quad \text{and}\quad
\cck\in
\R.
$$
Since $\yyda$ is integer valued, the $\la$-part becomes the following:
$$\aligned
\la(g; h)&=\explrtwopii{{\sumd{j\in\N,k\in\Ndz} \bbjk \edjn(g)\tedk(h)}},
\quad g\in N,\ h\in\Hdm,\\
&=\frac{f(h\inv gh)}{f(g)}=\frac{f((g\wedge h)g)}{f(g)}
=\ffdc(g\wedge h)\\
\endaligned
$$
$$\aligned
&=\explrtwopii{\sumd{1\leq j<k}\ccjk\edjk(g\wedge h)}\\
&=\explrtwopii{\sumd{1\leq j<k}\ccjk\Big(\edj(g)\edk(h)
-\edk(g)\edj(h)\Big)}\\
&=\explrtwopii{\sumd{1\leq j<k}\ccjk\Big(\pdj\edjn(g)\edk(h)
-\pdk\edkn(g)\edj(h)\Big)}.\\
\endaligned
$$
Hence we conclude that for $j<k$ and $i\in\Ndz$
$$\gathered
\bbiz\in \Z,\quad \bbii\in \Z,\  i\in\N;\\
\bbjk\equiv \ccjk\pdj,\ \bbkj\equiv -\ccjk\pdk\
\mod\ \Z, \quad j<k.\\
\endgathered
$$
Thus we have for $i<j$
$$\gathered
\bbij=\ccij\pdi+\mdij\quad\text{with some
}\mdij\in\Z;\\
\bbji=-\ccij\pdj+\mdji\quad\text{with some
}\mdji\in\Z;\\
\frac\bbij\pdi+\frac\bbji\pdj=\frac{\mdij}\pdi+\frac{\mdji}\pdj\in\lr{\frac1\pdi\Z}
+\lr{\frac1\pdj\Z}=\lr{\frac1{\lcm(\pdi, \pdj)}\Z}.
\endgathered
$$

Conversely suppose $(a, b)\in \tB$, i.e., 
$$\gathered
\aaijk\in\Z\quad \text{for}\quad i<j<k;\\
\aaiik,\ \aakik\in 2\Z\quad \text{for}\quad i<k,
\endgathered
$$
and
$$\gathered
\frac\bbij\pdi+\frac\bbji\pdj\in \lr{\frac1{\lcm(\pdi,\pdj)}\Z}\quad
\text{and}\quad
\bbii\in \Z;\\
\bbiz \in\Z,\quad i\in \N.
\endgathered
$$
So we can write
$$
\frac\bbij\pdi+\frac\bbji\pdj
=\frac\mdij\pdi+\frac\mdji\pdj\quad \text{with}\quad
\mdij,\mdji\in \Z.
$$
Set
$$\aligned
\ccij=\frac\bbij\pdi-\frac\mdij\pdi\quad \text{for } i<j,\quad \quad
\ccii=\bbii,
\endaligned
$$
so that
$$\aligned
\frac\bbji\pdj&=-\ccij+\frac{\mdji}\pdj.
\endaligned
$$
Then we have 
$$\aligned
\sumd{i,j\in\N}&\bbij\edin(g)\edj(h)\\
&\equiv 
\sumd{i<j}\ccij\Big(\pdi\edin(g)\edj(h)-\pdj\edjn(g)\edi(h)\Big) 
\mod \ \Z,\\
&=\sumd{i<j}\ccij\edij(g\wedge h).
\endaligned
$$
Thus with 
$$
\ffdc(g)=\explrtwopii{\sumd{1\leq i<j}\ccij\edij(g)}, \quad g\in L,
$$
we have
$$
\explrtwopii{\sumdij\bbij\edin(g)\edj(h)}=\frac{\ffdc(h\inv
gh)}{\ffdc(g)}=\part_1\ffdc(g; h),
$$
where $\edin(g)$ means $\edin\scirc \pig$.
We then compute the coboundary of $\ffdc$:
$$\aligned
\Big(\partl &\ffdc\Big)(g;
h)=\frac{\ffdc(g)\ffdc(h)}{\ffdc(gh)}, \quad g, h\in L,\\
&=\explrtwopii{\sumd{i<j}\ccij\lr{\edij(g)+\edij(h)-\edij(gh)}}\\
&=\explrtwopiim{\sumd{i<j}\ccij\edi(g)\edj(h)}=1,
\endaligned
$$
because $\edi(g)\in \pdi \Z$ and $\edj(h)\in \pdj\Z$ if $g, h\in
L$ and 
$$\aligned
\pdi\ccij\pdj&=\bbij\pdj-\mdij\pdj
\equiv \bbiz\qdj\equiv 0\quad \mod \ \Z.
\endaligned
$$
As $\aaijk\in\Z$ for every triplet $(i,j,k)\in\Delta$,  we get trivially
$$\aligned
\la_{a,0}=1,\quad \tmuda=\shs\muda\in \tztw(N, \T), \quad \text{and}\quad
\muda=\pigs(\tmuda).
\endaligned
$$
Since $\partn \uudaijk, i<j<k,$ is integer valued, the cochain
$$
\tmu_a^{ijk}=\explrtwopii{\uudaijk}
$$
belongs to $\tztw(N, \T)$. Since $\AS\uudaijk$ is integer
valued by (4.24), $\AS\tmu_a^{ijk}=1$ and therefore $\tmu_a^{ijk}\in
\tbtw(N, \T)$.
In view of the fact that
$$
\tmu_a^{ik}=\explrtwopii{\uudaik}=1,\quad i<k,
$$
we conclude that $\tmuda\in\tbtw(N, \T)$. Thus there exists a
cochain $\tilde f\in\tcone(N,  \T)$ such that
$$
\tmuda=\partn \tilde f.
$$
Define a cochain $f\in\tcone(L, \T)$ by
$$\aligned
f=(\pigs \tilde f)\ffdc.
\endaligned
$$
Then we get for each pair $g\in L, h\in\Hdm$
$$\aligned
(\partdone f)(g; h)&=\frac {f(h\inv
gh)}{f(g)}=\frac{\tilde f(\pig(h\inv gh))\ffdc(h\inv g
h)}{\tilde f(\pig(g))\ffdc(g)}\\
&=\frac{\ffdc(h\inv gh)}{\ffdc(g)}=\ladab(g; h);\\
(\partdtw f)(g; h)&=\frac{\tilde f(\pig(g))\ffdc(g)\tilde
f(\pig(h))\ffdc(h)} 
{\tilde f(\pig(gh))\ffdc(gh)}, \quad g, h\in L, \\
&=\partl \ffdc(g; h)\lr{\partn \tilde f}(\pig(g); \pig(h))\\
&=\tmuda(\pig(g); \pig(h))=\muda(g; h).
\endaligned
$$
Therefore we conclude
$$\aligned
\part f=\{\ladab, \muda\}\in \tB(\Hdm, L, M, \T).
\endaligned
$$
This completes the proof.
\QED
\enddemo
\proclaim{Lemma 4.6} The cocycle $\ladb$ corresponding to $b\in \Zdb$ does not depend on the $M$-component, i.e.,
$$
\ladb(mg; n\tilh)=\ladb(g; \tilh), \quad m, n\in M, g\in L, \tilh\in \Hdm.
$$
we will view $\ladb$ as a bicharacter on $N\times \Gdm$ rather than on $L\times \Hdm$.

{\rm i)} For $i\in\Z$, set
$$\aligned
\tzzdbii&=\lrbrace{z=(x, u)\in \R^2: \pdi x -\qdi u\in\Z},\quad
\tbbdbii=\Z\oplus \Z.
\endaligned
$$
The bicharacter $\la_z^{i,i}$ on $N\times \Gdm$ determined by\:
$$\aligned
\la_z^{i,i}(g; h)=\explrtwopii{x\edin(g)\tedi(h)+u\edin(g)\tedz(h)}, g\in N, h\in \Gdm,
\endaligned
$$
gives a characteristic cocycle of $\tZ(\Hdm, L, M, \T)$. It is a coboundary if and only
if $z$ is in $\tbbdbii$. The corresponding cohomology class $[\la_z^{i,i}]\in\Ladbii$
is given by\:
$$\aligned
\lrbracket{\la_z^{i,i}}=\lr{\lrbracket{\pdi x -\qdi u}_{\gcd(\pdi, \qdi)}, \lrbracket
{-\vdi x +\udi u}_\Z }\in \Z_{\gcd(\pdi, \qdi)}\oplus
(\R/\Z),
\endaligned
$$
where the integers $\udi, \vdi$ are determined by\:
$$
\pdi\udi-\qdi \vdi = \gcdlr{\pdi, \qdi}
$$
through the Euclid algorithm.

{\rm ii)}  Fix a pair $ i, j\in\N$ of indices and set
$$\gathered
\tzzdbij=\lrbrace{(x, u, y, v)\in\R^4: \pdj x-\qdj u\in\Z, \pdi y-\qdi v\in \Z};\\
\tbbdbij=
\lrbrace{(x, u, y, v)\in\tzdbij: \pdj x+\pdi y\in\gcd(\pdi,
\pdj)\Z, u, v\in \Z}.
\endgathered
$$
To each element $z=(x, u, y, v)\in \tzzdbij$, there corresponds a bicharacter $\la_z$ on $N\times
\Gdm$ determined by\:
$$\aligned
\la_z^{i,j}(g;
h)&=\explrtwopii{x{\edin(g)\tedj(h)}+y{\edjn(g)\tedi(h)}}\\
&\hskip.2in\times
\explrtwopii{u{\edin(g)\tedz(h)}+v{\edjn(g)\tedz(h)}},
\endaligned\ g\in N, h\in\Gdm,
$$
which is a characteristic cocycle in $\tZ(\Hdm, L, M, \T)$. It is a coboundary if and only if
$z\in\tbbdbij$. The cohomology class $[\la_z^{i,j}]\in\Ladbij$ of $\la_z$ corresponds to the parameter class\:
$$\aligned
\lrbracket{z}&=\pmatrix\lrbracket{\mdij\lr{x\rdji+y\rdij}-
\ndij\lr{u\sdji+v\sdij}}_\Z\\
 \lrbracket{\ydij\lr{x\rdji +y\rdij }+\xdij\lr{u\sdji+v\sdij}}_{\Z}\\
\lrbracket{-u\wdij+v\wdji}_\Z\endpmatrix
\in \pmatrix \left.\lr{\dfrac1{D(i,j)}\Z}\right/\Z\\\R/\Z\\\R/\Z\endpmatrix,
\endaligned
$$
where   $D(i,j), \cdots, \wdji$ are given in {\rm (4.17)} of {\rm Theorem 4.2}.
\endproclaim

\demo{Proof} i) Set 
$$
\Ddi=\gcd(\pdi, \qdi), \quad \rdi=\frac \pdi\Ddi, \quad \sdi=\frac\qdi\Ddi,
$$
and choose integers $\udi, \vdi\in \Z$  so that
$$
\rdi\udi-\sdi\vdi=1,
$$
where such a pair $(\udi, \vdi)\in \Z^2$ can be determined through the Euclid algorithm. Next we set 
$$\gathered
\edone=\lr{1,0}, \quad \edtwo=\lr{0,1};\\
\fdone=\udi\edone+\vdi\edtwo,\quad \edone=\rdi\fdone-\vdi\fdtwo,\\
\fdtwo=\sdi\edone+\rdi\edtwo,\quad\edtwo=-\sdi\fdone+\udi\fdtwo.\\
\endgathered
$$
Then
$\tzzdbii$ is given by the following: 
$$\aligned
\tzzdbii&=\lr{\frac1\Ddi\Z}\fdone+\R\fdtwo,
\endaligned
$$
and
$$
\tbbdbii=\Z\edone+\Z\edtwo=\Z\fdone+\Z\fdtwo,
$$
so that
$$\aligned
&\La_b(i,i)=\tzzdbii/\tbbdbii\cong \lr{\left.\frac 1\Ddi\Z\right/{ \Z}}\dot\fdone\oplus \lr{\R/\Z}\dot\fdtwo,
\endaligned
$$
where the dotted elements indicate the corresponding elements in the quotient group $\Ladbii$.
Now we chase the parameter:
$$\aligned
z&=x\edone+u\edtwo=x\lr{\rdi\fdone-\vdi\fdtwo}+u\lr{-\sdi\fdone+\udi\fdtwo}\\
&=\lr{\rdi x-\sdi u}\fdone+\lr{-\vdi x +\udi u}\fdtwo;\\
\dot z&=\lrbracket{\rdi x -\sdi u}_\Z \dot\fdone
+\lrbracket{-\vdi x +\udi u}_\Z\dot\fdtwo.
\endaligned
$$

$$\aligned
\la_{z}^{i,i}(g; \tilh)&=\explrtwopii{\lr{x\edin(g)\edi(\tilh)+u\edin(g)\edz(\tilh)}}
\endaligned
$$
for each pair $g\in N$ and $\tilh\in \Gdm$.

 ii)  First we fix the standard basis $\lrbrace{\edone,\cdots, \edfour}$ of $\R^4$ and set
$$\gathered
\ggdz=\rdij\edone-\rdji\edthr,\quad \ggdone=\udji\edone+\udij\edthr,
\endgathered
$$
where we choose $\udij,\udji\in\Z$ so that
$$
\rdij\udij+\rdji\udji=1.
$$
Since 
$$
\edone=\udij\ggdz+\rdji\ggdone, \quad 
\edtwo=-\udji\ggdz+\rdij\ggdone,
$$
we have 
$$
\Z\edone+\Z\edthr=\Z\ggdz+\Z\ggdone.
$$
Also we have
$$
\tbbdbij+\R\ggdz=\R\ggdz+\Z\ggdone+\Z\edtwo+\Z\edfour.
$$

Consider an integer $3\times 4$-matrix: 
$$\gathered
T=\pmatrix \mdij\rdji&-\ndij\sdji&\mdij\rdij&-\ndij\sdij\\\ydij\rdji&\xdij\sdji&\ydij\rdij&\xdij\sdij\\
0&-\wdij&0&\wdji\endpmatrix.
\endgathered
$$
We claim that 
$$
T\lr{\tzdbij+\R\ggdz}=\lr{\frac1\ddij\Z}\oplus \R\oplus \R.
$$
To prove the claim, for each vector
$$
z=x\edone+u\edtwo+y\edthr+v\edfour\in \R^4,
$$
we simply compute, 
$$\aligned
&\hskip1in
 T\ggdz=0,\\
Tz
&=\pmatrix \mdij\rdji&-\ndij\sdji&\mdij\rdij&-\ndij\sdij\\\ydij\rdji&\xdij\sdji&\ydij\rdij&\xdij\sdij\\
0&-\wdij&0&\wdji\endpmatrix
\pmatrix x\\u\\y\\v
\endpmatrix\\
&=\pmatrix\mdij\lr{x\rdji+y\rdij}-\ndij\lr{u\sdji+v\sdij}\\
 {\ydij\lr{x\rdji +y\rdij}+\xdij\lr{u\sdji+v\sdij}}\\
{-u\wdij+v\wdji}\endpmatrix.
\endaligned
$$
Suppose
$$\aligned
\frac k \ddij&=\mdij\lr{x\rdji+y\rdij}-\ndij\lr{u\sdji+v\sdij}\in \lr{\frac1\ddij}\Z.
\endaligned
$$
Then we have
$$\aligned
k&=\lr{\mdij\lr{x\rdji+y\rdij}-\ndij\lr{u\sdji+v\sdij}}\ddij \\
&=\lr{x\pdj-u\qdj}+\lr{y\pdi-v\qdi}\\
&=\lr{\lr{x+t\rdij}\pdj-u\qdj}+\lr{\lr{y-t\rdji}\pdi-v\qdi}.
\endaligned
$$
A choice of $t\in \R$, such that $\lr{x+t\rdij}\pdj-u\qdj$ is an integer, yields the integrality of the other term
$\lr{y-t\rdji}\pdi-v\qdi$, so that
$$\aligned
z+t\ggdz\in \tzdbij.
\endaligned
$$
Now we prove that
$$
T\inv \Z^3=\tbbdbij+\R\ggdz.
$$
Since $T$ is a matrix with integer coefficients and the generators $\ggdone, \edtwo, \edfour$ are all integer vectors, 
we have $T(\tbbdbij)\i \Z^3$. Conversely suppose that
$Tz\in \Z^3$. Then we have
$$\gathered
k=\mdij\lr{x\rdji+y\rdij}-\ndij\lr{u\sdji+v\sdij}\in\Z,\\
\ell=\ydij\lr{x\rdji +y\rdij}+\xdij\lr{u\sdji+v\sdij}\in\Z,\\
m=-u\wdij+v\wdji\in\Z.
\endgathered
$$
Hence we get
$$\gathered
x\rdji+y\rdij=\xdij k+\ndij\ell\in\Z,\quad 
n=u\sdji+v\sdij=-\ydij k+\mdij \ell\in\Z,\\
u=n\wdji-m\sdij\in\Z, \quad v=n\wdij+m\sdji\in\Z,\\
x\pdj+y\pdi=\lr{x\rdji+y\rdij}\Ddij\in \Ddij\Z.
\endgathered
$$
Therefore $z\in \tbbdbij+\R\ggdz$.

Consequently, we conclude
$$
\Ladbij\cong \tzdbij/\tbbdbij\cong \lr{\left.\lr{\frac1\ddij\Z}\right/\Z}\oplus\lr{\R/\Z}\oplus\lr{\R/\Z},
$$
in the sense that
 the cohomology class $\lrbracket{\la_z^{i,j}}\in \Ladbij$ corresponds to the following:
$$\aligned
\lrbracket{z}&=\pmatrix\lrbracket{\mdij\lr{x\rdji+y\rdij}-\ndij\lr{u\sdji+v\sdij}}_{\Z}\\
 \lrbracket{\ydij\lr{x\rdji +y\rdij }+\xdij\lr{u\sdji+v\sdij}}_{\Z}\\
\lrbracket{-u\wdij+v\wdji}_\Z\endpmatrix
\in \pmatrix\left. \lr{\frac1\ddij\Z}\right/\Z\\\R/\Z\\\R/\Z\endpmatrix.
\endaligned\nopagebreak
$$

For each $i, j\in\N$, define  maps $\pidiui: \Ladb(i,i)\mapsto \R/\Z$,
$\pidijui: \Ladb(i,j)\mapsto \R/\Z$, $\pidijuj: \Ladb(i,j)\mapsto \R/\Z$ and $\pidij: \Ladb(i, j)\mapsto \left.\lr{\frac 1\ddij\Z}\right/\Z$ by
the following:
$$\gathered
\pidiui([\la_z^{i,i}])=[u]_\Z\in \R/\Z,\\
\pidii\lr{\lrbracket{\la_z^{i,i}}}=\lrbracket{x\rdi-u\sdi}_\Z\in \left.\lr{\frac 1\Ddi\Z}\right/\Z\\
\endgathered
$$
for each $ z=(x, u)\in \tzzdbii$, and
$$
\gathered
\pidijui([\la_z^{i,j}])=[u]_\Z\in \R /\Z,\quad
\pidijuj([\la_z^{i,j}])=[v]_\Z\in\R/\Z,\\
\pidij\lr{\lrbracket{\la_z^{i,j}}}=\lrbracket{\mdij\lr{x\rdji+y\rdij}-\ndij\lr{u\sdji+v\sdij}}_{\Z}\in \left. \lr{\frac1\ddij\Z}\right/\Z
\endgathered
$$ 
for each $ z=(x, u, y, v)\in \tzzdbij$. The above maps $\pidijui$ and $\pidijuj$ are both well-defined because the coboundary condition on $z$ implies the integrality of $u$ and $v$.

Let $\Ladb$ be the set of all those
$$
\la_b=\{\ladb(i,i), \ladbij\}\in \prod_{i\in\N}\Ladbii\times \prod_{\underset{i,j\in\N}\to{i<j}}\Ladbij
$$
such that
$$\aligned
\pidiui(\ladb(i,i))=\pidijui(\ladbij)=\pidkiui(\ladbki)\quad 
\text{for all }\ i,j,k\in\N.
\endaligned
$$
Finally we have
$$
\La(\Hdm, L, M, \T)=\Lada\oplus \Ladb.
$$
This completes the proof.
\QED
\enddemo

{\smc Remark 4.6.} The direct sum homomorphism $\pidij\oplus\pidijui\oplus\pidijuj$ is a homomorphism of $\Lada(i,j)$ onto the direct sum group:
$$\CD
\Ladb(i,j)@>\pidij\oplus \pidijui\oplus\pidijuj>>\lr{\left.
\lr{\dfrac 1\ddij\Z}\right/\Z}\oplus \lr{\R/\Z}\oplus 
\lr{\R/\Z}.
\endCD
$$

i) By multiplying $\Ddi$ to $\pidii\lr{\la_z}$ we get
$$
\Ddi\pidii\lr{\lrbracket{\la_z}}
=\lrbracket{x\pdi-u\qdi}_{\Ddi\Z}\in \Z/\lr{\Ddi\Z};
$$
Similarly, we have
$$
D(i,j)\pidij(\la_z)=\lrbracket{\lr{x\pdj+y\pdi}-\lr{u\qdj+v\qdi}}_{D(i,j)}\in\Z/\lr{D(i,j)\Z}.
$$

ii) The kernel of $\pidij\oplus \pidijui\oplus\pidijuj$ is given by the following:
$$
\Ker\lr{\pidij\oplus \pidijui\oplus\pidijuj}
=\pmatrix \lrbrace{0}\\\left.\lr{\dfrac 1\mdij\Z}\right/\Z\\\lrbrace{0}\endpmatrix.
$$
At the parameter level, the kernel is described as follows:
$$
\lrbracket{\la_z}\in \Ker\lr{\pidij\oplus \pidijui\oplus\pidijuj}\Leftrightarrow x\pdj+y\pdi\in \ddij\Z, u, v \in\Z.
$$

\head{\bf{\pmb\S5}. The {\Rmhjrsq}.}
\endhead

We are now going to investigate the reduced modified \hjr:
$$\CD
\vdots@.\vdots\\
@VVV@VVV\\
\thtw(H,\T)@=\thtw(H,\T)\\
@V\Res VV@V\res VV\\
\La(\Hdm,L ,M, \T)@>\res>>\La(H, M, \T)\\
@V\d  VV@V\dhjr VV\\
\thmsout(G, N, \T)@>\partqm>>\ththr(G, \T)\\
@V\Inf VV@V\inf VV\\
\ththr(H, \T)@=\ththr(H, \T)
\endCD\tag5.1
$$
We refer to \cite{KtT3: page 116} for detail. 
So we first discuss  the second cohomology group $\tztw(H, \T)$ and the
restriction map $\Res$. Each second cocycle $\mu\in \tztw(H, \T)$ gives rise to a
group extension equipped with a cross-section:
$$\CD
1@>>>\T@>>>E@>j>\underset{\fsj}\to\longleftarrow>H@>>>1
\endCD
$$
such that
$$
\sj(g)\sj(h)=\mu(g; h)\sj(gh), \quad g, h\in H.
$$
With
$$
\ladmu(g; h)=\frac{\mu(h; h\inv gh)}{\mu(g; h)}, \quad g, h\in H,
$$
we obtain a characteristic cocycle $(\ladmu, \mu)\in  \tZ(H, H, \T)$. This
corresponds to the case that $P=1$ in the previous section. So we set
$$\gathered
\tztw=\left\{a\in \R^{\N^3}: \aaijk=0  \ \text{if } j\geq k,
\Asaijk\in\Z\right\};\\
\tbtw=\left\{a\in \tztw: \aaijk\in\Z, \aaiik,\aakik\in2\Z\right\}.
\endgathered\tag5.2
$$
\proclaim{Theorem 5.1} {\rm i)} Each element $a\in\tztw$ gives rise to a cocycle $\muda\in\tztw(H, \T)$\:
$$
\muda=\explr{\twopii\vvda}\in\tztw(H, \T), \quad a\in\tztw,
\tag5.3
$$
and the following diagram describes the second cohomology 
$\thtw(H, \T)$\:
$$\eightpoint\CD
1@>>>\tbtw@>>> a\in\tztw@>>>\lrbracket{a}\in\thtw@>>>1\\
@.@VVV@VVV@|\\
1@>>>\tbtw(H, \T)@>>>\muda\in \tztw(H, \T)@>>>\lrbracket{\lada}\in\thtw(H, \T)@>>>1
\endCD
$$
More precisely, with 
$$\left.\gathered
\tztw(i,j,k)=\{(x, y, z)\in\R^3: x-y+z\in\Z\},\quad \tbtw(i,j,k)=\Z^3,\\
\tztw(i, k)=\R^2, \quad \tbtw(i,k)=(2\Z)^2, \\
\thtw(i,j,k)=\tztw(i,j,k)/\tbtw(i,j,k), \quad
\thtw(i,k)=\tztw(i,k)/\tbtw(i,k)
\endgathered\right\}\tag5.4
$$
for each triplet $i<j<k$ {\rm(}resp. pair $i<k${\rm)} and
$$\gathered
\aaijk=x, \quad \aajik=y, \quad \aakij=z,\\
\lr{\text{\rm resp. }\quad \aaiik=x, \quad \aakik=y},
\endgathered
$$
we set
$$\gathered
\mudauijk=\explrtwopii{\vvdaijk}\in\tztw(H, \T);\\
\mudauik=\explrtwopii{\vvdaik}\in\tztw(H, \T).
\endgathered
$$
Then we have
$$\gathered
\tztw(H,\T)=\prod_{i<j<k}\tztw(i,j,k)\times
\prod_{i<k}\tztw(i,k),\\
\tbtw(H,\T)=\prod_{i<j<k}\tbtw(i,j,k)\times
 \prod_{i<k}\tbtw(i,k),\\
\muda=\lr{\prod_{i<j<k}\mudauijk}\lr{\prod_{i<k}\mudauik}
\in\tztw(H,\T),\\
\thtw(H,\T)\cong\prod_{i<j<k}\thtw(i,j,k)\times
\prod_{i<k}\thtw(i,k),\\
\lrbracket{\muda}=\lr{\lrbracket{\mudauijk}, \lrbracket{\mudauik}: i<j<k\ \text{\rm and}\ i<k}
\in \thtw(H,\T).
\endgathered
$$
Each $\thtw(i,j,k), i<j<k$, {\rm (}resp. $\thtw(i,k), i<k${\rm )}, is 
given by\:
$$\gathered
\thtw(i,j,k)\cong \lr{\R/\Z}\oplus\lr{\R/\Z},\\
\lr{\text{\rm resp.}\quad \thtw(i,k)\cong \lr{\R/2\Z}\oplus\lr{\R/2\Z}}.
\endgathered
$$
\endproclaim
\demo\nofrills{Proof.}\quad Most of the claims have been proved already except the claim for the structure of $\thtw(i,j,k)$. To prove the assertion on $\thtw(i,j,k)$, it is convenient to introduce a matrix $A\in \SL(3, \Z)$:
$$
A=\pmatrix1&-1&1\\0&1&0\\0&0&1\endpmatrix
\in \SL(3, \Z), \quad A\inv =\pmatrix1&1&-1\\0&1&0\\0&0&1\endpmatrix.
$$
We then observe that
$$
A\tztw(i,j, k)=\lr{\Z\oplus \R\oplus\R}, \quad 
A \tbtw=\Z^3,
$$
and conclude
$$\gathered
\thtw(i,j,k)\cong \lrbrace{0}\oplus \lr{\R/\Z}\oplus\lr{\R/\Z}.
\endgathered
$$
This completes the proof. 
\QED
\enddemo
\proclaim{Theorem 5.2} {\rm i)} Each second cocycle  $\muda\in\tztw(H, \T),\ a\in\tztw$, gives the corresponding characteristic 
cocycle\:
$$
\Res(\muda)=\lr{\lada, \muda}=\pims\lr{\lada|_{L\times\Hdm}, \muda|_L}\in\tZ(\Hdm, L, M, \T).
$$
The image $\Res\lr{\tztw(H, \T)}$ is therefore given by:
$$\gathered
\Res\lr{\tztw(H, \T)}=\lrbrace{\lr{\lada, \muda}: a\in \tzda,
\Asaijk\in\Z, i<j<k}.
\endgathered
$$
The $\lr{i,j,k}$-component $\Res(i,j,k)$ of the restriction map
$\Res$ gives rise to the following commutative diagram of short exact sequences\:
$$\eightpoint\CD
1@.1\\
@VVV@VVV\\
\tbtw(i,j, k)=\Z^3@>\xxdaijk\longrightarrow\xxdaijk>>\tbbbdaijk=\Z^3\\
@VVV@VVV\\
\tztw(i,j,k)=A\inv(\Z\oplus
\R^2)@>\xxdaijk\longrightarrow\xxdaijk>>\tzzdaijk
=A\inv\lr{\dfrac1D\Z\oplus \R^2}\\
@VVV@VVV\\
\thtw(i,j,k)=\{0\}\oplus\T^2@>\Res(i,j,k)>>\Lada(i,j,k)=\Z_D\oplus \T^2\\
@VVV@VVV\\
1@.0
\endCD
$$
where $D=D(i,j,k)=\gcd\lr{\pdi, \pdj, \pdk}$.
Also the restriction map $\Res_a(i,k): \thtw(i,k)\mapsto \Lada(i,k)$ is given
by
$$\CD
@.1@.1\\
@.@VVV@VVV\\
@.\tbtw(i,k)=(2\Z)^2@>\xxdaik\longrightarrow
\xxdaik>>\tbbbda(i,k)=(2\Z)^2\\
@.@VVV@VVV\\
@.\tztw(i,k)=\R^2@>\xxdaik\longrightarrow \xxdaik>>\tzzdaik=\R^2\\
@.@VVV@VVV\\
@. \thtw(i,k)=\lr{\R/2\Z}^2@>\Res(i,k)>>\Lada(i,k)=\lr{\R/2\Z}^2\\
@.@VVV@VVV\\
@.1@.1
\endCD
$$
Consequently, we get
$$\gathered
\Lada(i,j,k)/\Res(i,j,k)\lr{\thtw(i,j,k)}\cong\Z/\lr{D\Z},\\
\Lada(i,k)/\Res(i,k)\lr{\thtw(i,k)}\cong \lrbrace{0}.
\endgathered
$$

{\rm ii)}  The modified {\rm HJR}-map $\d: \La(\Hdm, L, M, \T)\mapsto \thmsout\lr{G, N, \T}$ enjoys the following properties\:
\roster
\item"a)" The $(i,j,k)$-component and $(i,k)$-component of
$\Ker(\d)$ are given by:
$$\gathered
\Ker(\d)_{ijk}=\lrbrace{0}\oplus \lr{\R/\Z}\oplus\lr{\R/\Z},\\
\Ker(\d)_{ik}=\lr{\R/2\Z}\oplus\lr{\R/2\Z}=\Lada(i,k).
\endgathered
$$
\item "b)" The image $\d\lr{\lrbracket{\lada, \muda}}\in\thmsout\lr{G, N, \T}, a \in \tzda,$ depends
only on the asymmetrization $\AS a$, i.e.,
$$\gathered
\d\lr{\lrbracket{\lada, \muda}}=\d\lr{\lrbracket{\la_\whaa, 1}}
\endgathered
$$
where
$$\gathered
\whaaijk=\Asaijk\in \lr{\frac1 D\Z},\quad i<j<k;\\
\whaajik=\whaakij=\whaaiik=\whaaijj=\whaakik=0.
\endgathered\tag5.5
$$
\item"c)" Set
$$
\tZ_\whaa=\lrbrace{a\in \tzda: a\text{ \rm  satisfies the requirement (5.5)}}.
$$
If $a\in\tzdwha$, then the image
$\cda=\d(\lada, 1)\in\tzout(\Gdm, N, \T)$ under the modified {\rm HJR}-map $\d$ is in the pull back $\pims\lr{\ththr(Q, \T)}$ and given by\:
$$\eightpoint\aligned
\cda&(\tqdone, \tqdtwo, \tqdthr)=\cda(\qdone,\qdtwo, \qdthr)\\
&=\explrtwopii{\sumdijk \aaijk \bracedpi{\edi(\qdone)}
\bracedpj{\edj(\qdtwo)} \bracedpk{\edk(\qdthr)}}
\endaligned\tag5.6
$$
for each
$\tqdone=\lr{\qdone, \sdone},\ \tqdtwo=\lr{\qdtwo, \sdtwo},
\ \tqdthr=\lr{\qdthr, \sdthr}\in\Qdm. $
\item"d)" The modified {\rm HJR}-map $\dhjr$ is injective on $\Ladb$ and $\Ker(\d)$ is precisely the connected component of $\La(\Hdm, L, M, \T)$. 
If $b\in\tzdb$, then 
$$
\lrbracket{\cdb, \nudb}=\d\lr{\ladb, 1}\in \tzmsout(G, N, \T)
$$
is given by\:
$${\eightpoint\aligned
\cdb&\lr{\tqdone, \tqdtwo, \tqdthr}
=\explrtwopii{\sumd{ i\in\N,j\in\Ndz}\bbij\edin\lr{\fnn(\tqdtwo; \tqdthr)}
\tedj\lr{\fs(\tqdone)}}
\endaligned}\tag5.7
$$
where
$$\gathered
\edin\lr{\fnn(\tqdtwo; \tqdthr)}=\frac{\edietan{\qdtwo}{\qdthr}}\pdi,\\
\tedi\lr{\fs(\tqdone)}=\bracedpi{\edi\lr{\qdone}}\ \text{for } i\geq 1, \quad \tedz\lr{\fs(\tqdone)}=\tedz\lr{\qdone} .\\
\endgathered\tag5.8
$$
The $d$-part $d_\cdb$ of $\cdb$ is given by $\nudb$:
$$\left.\aligned
d_\cdb\lr{\qdtwo; \qdthr}&=
\explrtwopii{\sumd{j\in\N}\bbjz\frac{\edjetan{\qdtwo}{\qdthr}}\pdj}\\
&=\explrtwopii{\frac{\bracettt{\nudb\lr{\fnn\lr{\qdtwo; \qdthr}}}}T},\\
\nudb\lr{g}&=\pi_T\lr{T\sumd{j\in\N}b(j,0)\edjn(g)}\in \rtz, \quad
g\in N,
\endaligned\right\}\tag5.9
$$
where $\pi_T: s\in\R\mapsto s_T=s+T\Z\in\rtz$ is the quotient map.
\endroster
The modular obstruction group $\thmsout(G, N, \T)$ looks like the following\:
$$\left.{\eightpoint\gathered
\thmsout(G, N, \T)=\tH_a^\out\oplus\tH_b^\out,
\quad \tH_b^\out\cong\Ladb,\\
\d\lr{\lrbracket{\lada, \muda}}=[c_{\AS a}]\in \prod_{i<j<k}\lr{\left.\lr{\frac1{\gcd\lr{\pdi,\pdj,\pdk}}\Z}\right/\Z},\ a\in \tZ_a,\\
\lrbracket{\cdb, \nudb}=\d\lr{\lrbracket{\ladb, 1}},\quad \nudb\in \Hom\lr{N, \rtz},\\
\lrbracket{\cdbuii}=\lr{\lrbracket{\pdi\bbii-\qdi\bbiz}_{\Ddi\Z}, \lrbracketd{-\vdi\bbii+\udi\bbiz}{\Z}}\\
\in \Z\Big/\lr{\Ddi\Z}
\oplus \R/\Z,\\
\lrbracket{\cdbuij}=\pmatrix\lrbracket{\mdij\lr{\bbij\rdji+\bbji\rdij}-
\ndij\lr{\bbiz\sdji+\bbjz\sdij}}_\Z\\
 \lrbracket{\ydij\lr{\bbij\rdji +\bbji\rdij }+\xdij\lr{\bbiz\sdji+\bbjz\sdij}}_{\Z}\\
\lrbracket{-\bbiz\wdij+\bbjz\wdji}_\Z\endpmatrix\\
\in \pmatrix \left.\lr{\dfrac1{D(i,j)}\Z}\right/\Z\\\R/\Z\\\R/\Z\endpmatrix,\quad
 D(i,j)=\gcdlr{\pdi,\pdj,\qdi,\qdj}.\\
\endgathered}\right\}\tag5.10
$$

{\rm iii)} The map $\partqm: \thmsout(G, N, \T)\mapsto \ththr(G, \T)$ in the \mhjr\ is given by\:
$$\left.\gathered
\partqm\lr{\lrbracket{c_{\whaa}} \lrbracket{\cdb \nudb}}=
\lrbracket{c_\whaa^G}\in\ththr(G, \T)=X^3(G, \T), \whaa\in \tZ_\whaa\\
c_\whaa^G=
\explrtwopii{\sumdijk\Asaijk\edi\ot\edj\ot\edk},\\
\partqm\lr{\thmsout(G, N, \T)}=\piqs\lr{\ththr(Q, \T)}.
\endgathered\right\}\tag5.11
$$
\endproclaim
\demo\nofrills{Proof.} \quad i) The assertion has been already proven.

ii) For each $i<j<k$, let
$$
D(i,j,k)=\gcdlr{\pdi,\pdj, \pdk}\in\Z.
$$
Fix $a\in\tzda$, i.e., $a\in\R^{\Delta}$ such that
$$\gathered
\Asaijk=\aaijk-\aajik+\aakij\in\lr{\frac1{D(i,j,k)}\Z}\\
\aaijk=0\quad \text{if}\quad j\geq k.
\endgathered
$$
Set
$$
\zdaijk=\pmatrix\aaijk\\\aajik\\\aakij\endpmatrix\in \tzda=A\inv
\pmatrix\lr{\frac1{D(i,j,k)}\Z}\\\R\\\R
\endpmatrix.
$$
Then we get
$$\gathered
A\zdaijk=\pmatrix \Asaijk\\\aajik\\\aakij
\endpmatrix
\in \pmatrix \frac1{D(i,j,k)}\Z\\\R\\\R\endpmatrix\\
A\tbbda(i,j,k)=\Z^3,
\endgathered
$$
so that
$$
\lrbracket{\ladauijk, \mudauijk}\sim\pmatrix \lrbracket{\Asaijk}_\Z\\
\lrbracket{\aajik}_\Z\\\lrbracket{\aakij}_\Z\endpmatrix
\in\pmatrix \lr{\frac1{D(i,j,k)}\Z}\\\R/\Z\\\R/\Z\endpmatrix.
$$
If $\Asaijk\in\Z$, the second cocycle $\mudauijk$ extends to a second cocycle on $H$ which gives $\lr{\ladauijk, \mu_a^{i,j,k}}=\Res\lr{\mu_a^{i,j,k}}$. Since $\text{\rm Range}(\Res)=\Ker(\d)$, the image $\d(\ladauijk, \mu_a^{i,j,k})$ depends only on the first term $\Asaijk$
of $A\zdaijk$. Hence we conclude $\d\lr{\lrbracket{\lada,\muda}}=\d\lr{\lrbracket{\la_\whaa}, 1}$. For $\Lada(i,k)$, we have
$$
\Lada(i,k)=\Res(i,k)\lr{\thtw(i,k)},
$$
so that the map $\d$ kills the entire $\Lada(i,k)$. This completes the proof of (iia) and (iib).

iic) Set $\cda=\d(\lada, \muda)$ with $a\in \tzdwha$.
We then look at the 
crossed extension $E_{\lada, \muda}\in\X\lr{\Hdm, L, M, \T}$
$$\CD
1@>>>\T@>>>E@>j>\underset{\sj}\to\longleftarrow>L@>>>1.
\endCD
$$
As 
$$\aaijk\in\lr{ \frac1{\gcdlr{\pdi,\pdj, \pdk}}\Z} 
\quad \text{and}\quad
 \edi(g)\in \pdi\Z,\quad g\in L,
$$
we have 
$\muda=1$. Hence observing that $\lada(g; \tilh)=1$ for every $g\in L\wedge \Hdm$ and $\tilh\in\Hdm$,  we get from (3.15) and (3.16) the following:
$$\aligned
\cda\lr{\tqdone, \tqdtwo, \tqdthr}
&=\a_{\fs(\tqdone)}\lr{\sj(\fnl(\tqdtwo; \tqdthr))} 
\sj\lr{\fnl\lr{\tqdone; \tqdtwo\tqdthr}}\\
&\hskip.5in\times
\lrbrace{\sj\lr{\fnl(\tqdone; \tqdtwo)}\sj\lr{\fnl(\tqdone\tqdtwo;\tqdthr)}}\inv\\
&=\lada\lr{\fs(\tqdone)\fnl\lr{\tqdtwo; \tqdthr}\fs(\tqdone)\inv; \fs(\tqdone)}\\
&=\lada\lr{\lr{\fs(\tqdone)\wedge\fnl(\tqdtwo; \tqdthr)}\fnl(\tqdtwo; \tqdthr); \fs(\tqdone)}\\
\endaligned
$$
$$\aligned
&=\lada\lr{\fs(\tqdone)\wedge\fnl(\tqdtwo; \tqdthr); \fs(\tqdone)}\lada\lr{\fnl(\tqdtwo; \tqdthr); \fs\lr{\tqdone}}\\
&=\lada\lr{\fnl(\tqdtwo; \tqdthr); \fs\lr{\tqdone}}\\
&=\explrtwopii{\sumdijk\aaijk\edjk\lr{\fnl(\tqdtwo; \tqdthr)}\edi\lr{\fs(\tqdone)}}\\
&=\explrtwopii{\sumdijk\aaijk \bracedpi{\edi(\tqdone)}\bracedpj{\edj(\tqdtwo)}\bracedpk{\edk(\tqdthr)}}\\
&=\explrtwopii{\sumdijk\aaijk\bracedpi{\edi(\qdone)}\bracedpj{\edj(\qdtwo)}\bracedpk{\edk(\qdthr)}}\\
&=\cda\lr{\qdone; \qdtwo; \qdthr}
\endaligned
$$
for each $\tqdone=\lr{\qdone,\sdone}, \tqdtwo=\lr{\tqdtwo, \sdtwo}, \tqdthr=\lr{\qdthr, \sdthr}\in \Qdm$. Thus the assertion (iic) follows.

iid) Since $\Res\lr{\thtw(H, \T)}\cap \Ladb=\lrbrace{0}$, the modified HJR-map $\d$ is injective on $\Ladb$. Now fix $b\in \tzdb$. Since $\mudb=1$ and $\ladb(m; \tilh)=\lrbrace{1}$ for every pair $m\in M, \tilh\in\Hdm$, we have, as in (iic), the following:
$$\aligned
\cdb(\tqdone; \tqdtwo;\tqdthr)&=\ladb(\fnn(\qdtwo; \qdthr); \fs(\tqdone))\\
&=\explrtwopii{\sumd{i\in\N,j\in\Ndz}\bbij\edin\lr{\fnn(\qdtwo; \qdthr)}\tedj\lr{\fs(\tqdone)}}\\
&=\explrtwopii{\sumd{i,j\in\N}\bbij\edin\lr{\fnn(\qdtwo; \qdthr)}\edj\lr{\fs(\qdone)}}\\
&\qquad\times
\explrtwopii{\sumd{i\in\N}\bbiz\edin\lr{\fnn(\qdtwo; \qdthr)}\tedz\lr{\tqdone}}
\endaligned
$$
where $\edin\lr{\fnn(\qdtwo; \qdthr)}$ is given by (5.8). Also we compute
$$\aligned
d_\cdb(\qdtwo; \qdthr)&=\ladb\lr{\fnn(\qdtwo; \qdthr); \zdz}=\explrtwopii{\frac{\nudb\lr{\fnn(\qdtwo; \qdthr)}}{T}}\\
&=\explrtwopii{\sumd{i\in\N}\bbiz \edin\lr{\fnn(\qdtwo; \qdthr)}},\\
\nudb(g)&=\pi_T\lr{T\sumd{i\in\N}\bbiz\edin(g)}\in\rtz, \quad g\in N,
\endaligned
$$ 
with $\pi_T: s\in \R\mapsto s_T=s+T\Z\in\rtz$ the quotient map. 

The last assertion, (5.10), on $\thmsout(G, N, \T)$ follows almost automatically from the above computations and Lemma 4.6 in the last section.

iii) We now compute the map 
$$
\part_\pim: \thmsout(G, N, \T)\mapsto\ththr(G, \T).
$$
 We continue to work on the cocycle $(\ladab, 1)$ for
$a\in \tzdwha$ whose
restriction to $\{\Hdm, K\}$ gives rise to the crossed extension
$U\in\X(\Hdm, K, \T)$:
$$\CD
1@>>>\T@>>>U@>j>\underset{\sj }\to\longleftarrow>K@>>>1
\endCD
$$
where the group $K$ is given by the following:
$$\aligned
K&=\Ker(\nudb\scirc\pig)=\lrbrace{g\in L: \sumd{i\in\N}\bbjz\edjn(g)\in \Z}.
\endaligned
$$
Then the following third  cocycle $c_G\in\tzthr(G, \T)$:
$$\aligned
c_G(\ggdone; \ggdtwo;
\ggdthree)&=\a_{\sh(\ggdone)}\Big(\sj(\fndm(\ggdtwo;
\ggdthree)\Big)\sj\lr{\fndm(\ggdone; \ggdtwo\ggdthree)}\\
&\hskip,5in\times
\Big(\sj(\fndm(\ggdone; \ggdtwo))\sj(\fndm(\ggdone \ggdtwo; \ggdthree))
\Big)^{-1}\\
&=\ladab(\fndm(\ggdtwo; \ggdthree); \ggdone)
=\lada(\fndm(\ggdtwo; \ggdthree); \ggdone)\\
&=\explrtwopii{\sumdijk\aaijk\edi(\ggdone)\edj(\ggdtwo)\edk(\ggdthree)}\\
&=c_a^G(\ggdone;\ggdtwo; \ggdthree), \quad \ggdone, \ggdtwo, \ggdthree\in
G,
\endaligned
$$
is precisely the image $\part_\pim \scirc \d(\ladab, 1)$. 
\QED
\enddemo

\head{\bf\S6. Concluding Remark.}
\endhead
The history of cocycle (resp. outer) conjugacy analysis of group actions and group outer actions on an AFD factor goes back to the grand work of Connes, \cite{Cnn3, 4}, in the mid 1970's. Since then, the steady progress was accomplished by several hands through the three decades following Connes work,  the works of  V.F.R. Jone and A. Ocneanu are noteworthy, \cite{Jn, Ocn}. 

We have now computed the invariants, which determine the outer conjugacy class, of an outer action of a {\cdabg} on an AFD factor of type \threel, $0<\la<1$. The reduction of  outer conjugacy analysis of an outer action of a {\cdag} on an AFD factor of type {\threel} down to the associated complete invariants was successfully carried out in our previous work, \cite{KtT1,2,3}. As we have demonstrated in this paper, the computation of invariants is doable as soon as the group in question is specified, except the case of type \threez.

\subhead\nofrills{Toward One Parameter Automorphism Group:\quad}
\endsubhead
After the completion of cocycle (resp. outer) conjugacy classification of  {\cdag} (resp. outer) actions on an AFD factor, it is only natural to consider the same problem for a continuous group. The first step to this goal is obviously the study of one parameter automorphism group $\lrbrace{\a_t: t\in \R}$ of an {\afdf} $\sRdz$ of type \twoone. The first steps were already taken by Y. Kawahigashi, \cite{Kw1, 2, 3, 4}, who classified, up to cocycle (or stable) conjugacy, the most of one parameter automorphism groups of $\sRdz$ constructed from concrete data, which was extended to the case of type \threee\ by U-K. Hui, \cite{Hu}. But the general ones with full Connes spectrum are left untouched.
One of difficulties is the lack of technique which allows us to create a one cocycle $\lrbrace{u_s: s\in\R}$ for a projection 
$p\in \Proj(\sRdz)$ so that the perturbed one parameter automorphism group $\lrbrace{\Ad(u_t)\scirc \a_t: t\in\R}$ leaves the projection $p$ invariant which allows us to localize the analysis of the action. If a projection $p\in \Proj(\sRdz)$ is differentiable relative to $\a$, then the derivation $\d_\a$ associated with $\a$ generates a desired cocycle. But we don't know the answer to the following basic question:
\proclaim\nofrills{Question:} \quad Does the {\rm C$^*$}-algebra\: 
$$
A=\lrbrace{x\in \sRdz: \lim_{t\to0}\|x-\a_t(x)\|=0}
$$ 
contain a non-trivial projection?
\endproclaim
If $p\in \Proj(A)$, then for each smooth function $f\in C_c^\infty(\R)$ with compact support the element: 
$$
p(f)=\a_f(p)=\int_\R f(t)\a_t(p)\txd t
$$
is smooth and one can choose $f$ such a way that
$\|p-p(f)\|$ is arbitrarily small so that $\Sp(p(f))$ is concentrated on a neighborhood of the two points $\lrbrace{0, 1}$, which allows us to generate  a non trivial differentiable projection $q$ near $p$ via contour integral:
$$
q=\frac1\twopii \oint_{|z-1|=r}\lr{z-p(f)}\inv \text d z.
$$

On the other hand, thanks to the exponential functional calculus, one can generate plenty of differentiable unitaries. For example, if $h\in A_{\text {s.a}}$, then for a real valued smooth function $f$, we get a differentiable unitary element $\exp\lr{\txti f(h)}$ of $A$ which can stay near the unitary $\exp(\txti h)$ in norm. Hence the group of differentiable unitaries is $\sigus$-strongly dense in the unitary group $\sU(\sRdz)$.

\enddocument